\documentclass{amsart}
\usepackage{color}

\usepackage{amsmath}
\usepackage{amssymb}
\usepackage{mathrsfs}
\usepackage{mathtools}

\newtheorem{theorem}{Theorem}[section]
\newtheorem{claim}[theorem]{Claim}

\newtheorem{conclusion}[theorem]{Conclusion}
\newtheorem{observation}[theorem]{Observation}

\theoremstyle{definition}
\newtheorem{definition}[theorem]{Definition}

\theoremstyle{remark}
\newtheorem{remark}[theorem]{Remark}
\newtheorem{question}[theorem]{Question}
\newtheorem{notation}[theorem]{Notation}

\DeclareMathOperator{\Sep}{Sep}
\DeclareMathOperator{\Sol}{Sol}

\newcommand{\trp}{{\rm trp}}

\newcommand{\otp}{{\rm otp}}

\newcommand{\ZFC}{{\rm ZFC}}

\newcommand{\cd}{{\rm cd}}
\newcommand{\cf}{{\rm cf}}
\newcommand{\pr}{{\rm pr}}
\newcommand{\dom}{{\rm dom}}
\newcommand{\ran}{{\rm ran}}
\newcommand{\Dom}{{\rm Dom}}

\newcommand{\tr}{{\rm tr}}
\newcommand{\bd}{{\rm bd}}

\newcommand{\id}{{\rm id}}

\newcommand{\nacc}{{\rm nacc}}

\newcommand{\cov}{{\rm cov}}

\newcommand{\tcf}{{\rm tcf}}

\newcommand{\lqq }{{``}}

\newcommand{\dsize}{{\text{ }}}
\newcommand{\demo}{{\text{ }}}

\newcommand{\dbcu}{{\bigcup }}
\newcommand{\Bb}{{\rm BB}}
\newcommand{\midia}{{super BB}}

\newcommand{\rest}{{\restriction}}

\newcommand{\wilog}{{\rm without loss of generality}}

\newcommand{\varp}{{\varepsilon}}

\newcount\skewfactor
\def\mathunderaccent#1#2 {\let\theaccent#1\skewfactor#2
\mathpalette\putaccentunder}
\def\putaccentunder#1#2{\oalign{$#1#2$\crcr\hidewidth
\vbox to.2ex{\hbox{$#1\skew\skewfactor\theaccent{}$}\vss}\hidewidth}}

\newbox\noforkbox \newdimen\forklinewidth
\setbox0\hbox{$\textstyle\bigcup$}
\setbox1\hbox to \wd0{\hfil\vrule width \forklinewidth depth \dp0
                        height \ht0 \hfil}
\setbox\noforkbox\hbox{\box1\box0\relax}
\def\unionstick{\mathop{\copy\noforkbox}\limits}
\def\nonfork#1#2_#3{#1\unionstick_{\textstyle #3}#2}
\def\nonforkin#1#2_#3^#4{#1\unionstick_{\textstyle #3}^{\textstyle
    #4}#2}
\setbox0\hbox{$\textstyle\bigcup$}
\setbox1\hbox to \wd0{\hfil{\sl /\/}\hfil}
\setbox2\hbox to \wd0{\hfil\vrule height \ht0 depth \dp0 width
                                \forklinewidth\hfil}
\newbox\doesforkbox
\setbox\doesforkbox\hbox{\box1\box0\relax}
\def\nunionstick{\mathop{\copy\doesforkbox}\limits}

\def\fork#1#2_#3{#1\nunionstick_{\textstyle #3}#2}
\def\forkin#1#2_#3^#4{#1\nunionstick_{\textstyle #3}^{\textstyle
    #4}#2}

\newcommand{\stickT}{%
\setbox255=\hbox{\raise1ex\hbox{$\hspace{0.2pt}\,\bullet\,$}}
\mathord{\rlap{\hbox to\wd255{\hss\hbox{$|$}\hss}}
\box255}
}
\newcommand{\stickS}{%
\setbox255=\hbox{\raise0.6ex\hbox{$\scriptstyle\bullet$}}
\mathord{\rlap{\hbox to\wd255{\hss\hbox{$\scriptstyle|$}\hss}}
\box255}
}

\newenvironment{PROOF}[2][\proofname.]
   {\begin{proof}[#1]\renewcommand{\qedsymbol}{\textsquare$_{\rm #2}$}}
   {\end{proof}}

\numberwithin{equation}{section}

\usepackage[hidelinks]{hyperref}
\begin{document}
\makeatletter\def\shfiuwefootnote{\gdef\@thefnmark{}\@footnotetext}\makeatother\shfiuwefootnote{Version 2023-05-01\_2. See \url{https://shelah.logic.at/papers/775/} for possible updates.}

\parindent 0pt

\parskip 3pt

\title {Super Black Box, (ex. Middle Diamond) \\
775}  
\author {Saharon Shelah}
\address{Einstein Institute of Mathematics\\
Edmond J. Safra Campus, Givat Ram\\
The Hebrew University of Jerusalem\\
Jerusalem, 9190401, Israel\\
 and \\
 Department of Mathematics\\
 Hill Center - Busch Campus \\
 Rutgers, The State University of New Jersey \\
 110 Frelinghuysen Road \\
 Piscataway, NJ 08854-8019 USA}
\email{shelah@math.huji.ac.il}

\urladdr{http://shelah.logic.at}


\thanks{First typed: June 13, 2001. The author would like to thank ISF (Israel Science Foundation)    for
partially supporting this research by grant 1838(19). In earlier versions (up to 2019), the author thanks Alice Leonhardt for the beautiful typing. In later versions, the author would like to thank the typist for his work and is also grateful for the generous funding of typing services donated by a person who wishes to remain anonymous.  References like [Sh:950, Th0.2=Ly5] mean that the internal label of Th0.2 is y5 in Sh:950. The reader should note that the version in the author's website is usually more up-to-date than the one in arXiv. This is publication number 775 in Saharon Shelah's list.}

 

\makeatletter
   \@namedef{subjclassname@2020}{\textup{2020} Mathematics Subject Classification}
   \makeatother



\date{May 1, 2023}

\begin{abstract} This is a lightly corected version of an old work.
    
    Under certain cardinal arithmetic assumptions, we prove that for every large enough regular $\lambda$ cardinal, for many regular $\kappa < \lambda$,  many
    stationary subsets of $\lambda$ concentrating on cofinality $\kappa$  has the ``\midia".  In particular, we  have the \midia on $\{\delta < \lambda \colon \cf(\delta) = \kappa\}$. This is a strong negation of uniformization. 
    
    We have added some details to Case 2 in the  proof of Claim \ref{d.7}.  Works continuing it are \cite{Sh:898}    and \cite{Sh:1028}.  We thank Ari Brodski and Adi Jarden for their helpful comments.
    
    In this paper (and inside it) we had earlier used the notion \lqq middle diamond" 
    which is now replaced by ``super $\Bb$'', that is, ``super black box'', in order to be consistent with other papers (see \cite{Sh:898}).
\end{abstract}

\maketitle

\setcounter{section}{-1}

\newpage

\centerline{Annotated Content}

\S0 \quad  Introduction,  pg. \pageref{0}.

\begin{enumerate}
    \item[${{}}$] [We state concise versions of the main results: provably in ZFC there are many cases of the \midia, for  many triples $(\lambda,\Upsilon,\kappa)$ where $\kappa < \lambda$ are regular, $2 \le \Upsilon < \lambda$ we can give for any $\delta \in S^\lambda_\kappa$ a function $f_\delta$ from some subset  $C_\delta$ of $\delta$ of order type $\kappa$ unbounded in $\delta$ into $\Upsilon$ such that any $f \colon \lambda \rightarrow \Upsilon$ extends $f_\delta$ for stationarily many $\delta$'s.  We also review some history.]
\end{enumerate}

\S1 \quad Super black box: sufficient conditions for $\cf(2^{\mu})$,  pg.  \pageref{1}.

\begin{enumerate}
    \item[${{}}$]  [We  define families of relevant sequences $\bar C = \langle C_\delta \colon \delta \in S \rangle$ and properties needed for phrasing and proving our results (see Definition \ref{d.5}) and give easy implications; (later in Definition d.12A, \ref{d.14} we define other versions of \midia).  Then (in Claim \ref{d.6}) prove a version of the \midia for $(\lambda, \Upsilon, \theta, \kappa)$ guessing only the colour of $f \, {\rest} \, C_\delta$ when $\lambda = \cf(2^\mu)$, the number of possible colours, $\theta$, is not too large and $\kappa$ belongs to a not so small set of regular cardinals $\kappa < \mu$ (e.g., every large enough regular $\kappa <\beth_\omega$ if $\beth_\omega \le \mu$).  The demand on $\kappa$  says $2^\mu$ is not the number of $\kappa$-branches of a tree with $< 2^\mu$ nodes. We then prove that a requirement on $(\mu,\theta)$ is reasonable (see Claim \ref{d.7}, \ref{d.8}, \ref{d.8Z}).  In Claim \ref{d.12} we get the version of \midia giving $f \, {\rest} \, C_\delta$ from the apparently weaker one guessing only the colour of $f \, {\rest} \, C_\delta$.  We quote various needed results and definitions in Definition \ref{d.3}, Claim \ref{d.13}.]
\end{enumerate}

\S2 \quad Upgrading to larger cardinals,  pg. \pageref{2}.

\begin{enumerate}
    \item[${{}}$]   [Here we try to get such results for ``many" large enough regular $\lambda$.  We succeed to get it for every large enough not strongly inaccessible cardinal.  We prove upgrading: i.e., give sufficient conditions for deducing \midia on $\lambda_2$ from ones on $\lambda_1 < \lambda_2$ (see Claim \ref{x.20}) and Claim x.21A (many disjoint ones), Claim \ref{x.20.9}, \ref{x.20.10}). We deduce (in \ref{x.21}) results of the form ``if $\lambda$ is large enough to satisfy a weak requirement, then it has \midia for quite many $\kappa,\theta$", using upgrading claims on the results of \S1 and generalization of \cite{EK} (see Claims \ref{x.23}, \ref{x.23.9}, and  \ref{x.2.9}).  We deal with such theorem more than necessary just for the upgrading.]
\end{enumerate}

\newpage

\section{Introduction}\label{0}

Using the diamond (discovered by Jensen \cite{J}) we can  construct many objects; wonderful but there is a price - we must assume  $\mathbf V = \mathbf L$ or at least instances of G.C.H. (see Gregory
\cite{Gr}, Shelah \cite{Sh:108}).  In fact $\lambda = \mu^+ = 2^\mu > \beth_\omega \Rightarrow \diamondsuit_{\{\delta <
\lambda \colon  \cf(\delta) = \kappa\}}$ for every regular large enough $\kappa$ depending on $\mu$ (see \cite{Sh:460} and citations there). 
Naturally, we would like to get the benefit without paying the price. A way to do this is to use the weak diamond (see \cite{Sh:65}, \cite[AP,~\S1]{Sh:f}, \cite{Sh:208}, \cite{Sh:638} and lately
\cite{Sh:755}, see on the definable version in \cite[AP,~\S3]{Sh:f}, \cite{Sh:576}).  
In the diamond principle  for a stationary $S
\subseteq \lambda$, for $\delta \in S$ we guess $\eta \, {\rest} \, \delta$ for $\eta \in {}^\lambda 2$.  In the weak diamond, for stationary $S \subseteq \lambda$ for $\delta \in S$ we do not guess $\eta \, {\rest} \, \delta$ for $\eta \in {}^\lambda
2$, we just guess whether it satisfies a property - is black or white, i.e., the value of $\mathbf F(\eta \, {\rest} \, \delta) \in \{0,1\}$ for a pre-given $\mathbf F$.  Many times it is just enough in constructions using diamond to ``avoid undesirable company": whatever you  will do I will do otherwise.  This is, of course, a considerably weaker principle, but for $\aleph_1$ it follows from CH (which does not imply $\diamondsuit_{\aleph_1}$) 
and even WCH (i.e. $2^{\aleph_0} < 2^{\aleph_1}$).  In fact, it holds for
many cardinals: $2^\lambda > 2^{< \lambda} = 2^\theta$ suffices. Having $\mathbf F$ with two colours (= values) seemed  essential.  In \cite{Sh:638}, we get more colours for $\aleph_1$ if the club filter on $\omega_1$ is $\aleph_2$-saturated but this is a heavy assumption.  In \cite{Sh:755}, if $\lambda = \lambda^{< \lambda} = 2^\mu$, it is proved that for ``many'' regular $\kappa < \lambda$, if $S \subseteq \{\delta < \lambda \colon \text{cf}(\delta) = \kappa\}$ is stationary then we have the weak diamond for $S$.  Now what does ``many" mean?  For example,  if $\lambda \ge \beth_\omega$, ``many"  includes every large  enough regular $\kappa < \beth_\omega$; note that it is for ``every stationary $S$'' unlike the situation, e.g. for $\aleph_1$ under CH (see \cite[V]{Sh:f}).  We get there more  cases if the colouring $\mathbf F$ on ${}^\delta 2$ depends only on $\eta \, {\rest} \, C_\delta$,  where $C_\delta \subseteq \delta = \sup(C_\delta)$ and $C_\delta$ is
``small".  Here we get more than two colours and moreover: we guess $\eta \, {\rest} \, C_\delta, C_\delta \subseteq \delta = \sup(C_\delta)$ and so can choose $\langle M_\delta \colon \delta \in S \rangle$ such that $M_\delta$ is a model with universe $C_\delta$ and relations only 
(and e.g. vocabulary $\subseteq {\mathscr  H}(\aleph_0))$ such that if $M$ is a model with
universe $\lambda$ and relations only 
(and vocabulary $\subseteq {\mathscr  H}(\aleph_0))$ then for stationarily many $\delta \in S,M \, {\rest} \, C_\delta$ is
$M_\delta$. Of course, if $(\exists \mu < \lambda)[2^\mu > \lambda]$  then we cannot have $C_\delta = \delta$ so this is a very reasonable restriction. (It would be even better to have $M_\delta \prec M$ below but we do not know).  Also for ``almost" all $\lambda$ large enough  the result holds for many regular $\kappa$. 

The motivation has been trying to get better principles provably in $\ZFC$.  We use partial squares and $\hat{I}[\lambda]$, see \cite{Sh:108}, \cite{Sh:221}, \cite{Sh:237e}, \cite[\S1]{Sh:420} and \cite{Sh:562}. We prove (in fact in Theorem \ref{0.1} we can replace $S_{\kappa}^{\lambda}$ by for quite  ``many'' stationary $S \subseteq S^\lambda_\kappa$).

\begin{theorem}\label{0.1} 
    For every $\mu$ strong limit of
    uncountable cofinality,  for every regular $\lambda > \mu$ which is not strongly inaccessible, for every large enough regular $\kappa < \mu$, for every $\chi < \mu, S^\lambda_\kappa$ has the $\chi$-$\Bb$ 
    (called later ``$\lambda$ has the $(\kappa,\chi^+)$-MD$_2$-property"),
    which means:
    
    \begin{enumerate} 
        \item[$\boxplus_{\lambda,\kappa}$]
        for  some stationary $S \subseteq \{\delta < \lambda \colon {\text{\rm cf\/}}(\delta) = \kappa\}$ for every relational vocabulary
        $\tau$ of cardinality $< \chi$ (so $\tau$ has no function symbols) 
        we can find $\bar M = \langle M_\delta \colon \delta \in S \rangle$ which is
        a $\tau$-MD-sequence such that: 
    
        \begin{enumerate}
            \item[(a)]  $M_\delta$ is a $\tau$-model,
    
            \item[(b)] the universe of $M_\delta$ is an unbounded subset of $\delta$ of order type $\chi \times \kappa$ including $\chi,$
    
            \item[(c)] if $M$ is a $\tau$-model with universe $\lambda$ \underline{then} for stationarily many $\delta \in S$ we have $M_\delta \subseteq M$. 
        \end{enumerate}
    \end{enumerate}
\end{theorem} 

It seems that this is stronger in sufficient measure from the weak diamond to justify this name; there are many properties between the weak diamond and the property from Theorem \ref{0.1} and  we have to put the line somewhere.   We shall use the adjective middle diamond to describe even cases with coloring of ${}^{\alpha} 2 \ (\alpha < \lambda)$ or just functions with domain $\alpha$ for $\alpha < \lambda$ (justified by Claim \ref{d.12}). 

I thank Todd Eisworth and the referee and Adi Yarden  for their helpful comments.
Some improvements and elaborations are delayed to \cite{Sh:F611}. This work is continued in \cite{Sh:829}.

\begin{theorem}\label{0.2} 
    In \ref {0.1} we can have also $\bar f
    = \langle f_\delta \colon \delta \in S \rangle$ where  $f_\delta \colon u_\delta \rightarrow \delta, u_\delta \subseteq
    \delta, \ \lambda^{|u_\delta|} = \lambda$ such that: 
    
    \begin{enumerate} 
        \item[(d)]   if $M$ is a $\tau$-model with universe $\lambda$ and $f \colon \lambda \rightarrow \lambda$ and $u \subseteq \lambda,\lambda^{|u|} = \lambda$, 
    \end{enumerate} 
    
    \underline{then} for stationarily many $\delta \in S$ we have $M_\delta \subseteq M \, \wedge \, (f_\delta = f \, {\rest} \, u \subseteq f) \, \wedge \, u_\delta = u$.  
\end{theorem} 

Now, we prove Theorem \ref{0.1} and Theorem \ref{0.2}: 

\begin{proof}
    By Claim \ref{x.24}(1) if $\lambda > \mu$ is regular and not strongly inaccessible, then
    for every large enough regular $\kappa < \mu$ the following holds (see
    Definition \ref{d.12G}(4)):
    
    \begin{enumerate}  
        \item[$(*)_\kappa$]  if $\theta < \mu$ then some $(\lambda,\kappa,\kappa^+)$-parameter $\bar C$ has the  $\theta$-\Bb-property; let 
        $\bar C = \langle C_\delta \colon \delta \in S \rangle$.
    \end{enumerate} 
    
    Now \wilog \, $\delta \in S \Rightarrow (\chi \times \theta)^\omega \, \vert \, \delta$, and let $$C'_\delta = \{\gamma \colon \gamma < \chi \text{ or } (\exists \beta \in C_\delta)(\beta \le \gamma < \beta +\chi).$$ 
    
    Lastly, let $\theta = 2^{\chi + \kappa + |\tau|}$ and by Claim \ref{d.12} we get
    the conclusion of Theorem \ref{0.1}. As for Theorem \ref{0.2}, Clause (d) is easy by having a pairwise disjoint sequence of suitable $\langle S_i \colon i < \lambda \rangle$ which  holds by Claim \ref{x.24}(2).
\end{proof}
 
Note that having $M_\delta$'s is deduced in Claim \ref{d.12} from ``$\bar C$ has the $\theta$-\Bb-property" with $\theta$ large enough.  Note 
that Theorem \ref{0.1} has obvious applications in various constructions, e.g., the construction of
abelian groups $G$ with Ext$(G,K) \ne \emptyset$.

Note that in Theorem \ref{0.1}, Claim \ref{x.24} instead of using ``large enough $\kappa < \mu \ldots$", $\mu$ strong limit and quoting \cite{Sh:460} we can specifically write what we use on
$\kappa$ (and on $\theta$). Also in Theorem \ref{0.1}, if we omit cf$(\mu) > \aleph_0$ (but still $\mu$
is uncountable) the case $\lambda$ successor is not affected and in the case $\mu =
\beth_\delta$ when $\omega^2$ divides $\delta$, the statement instead ``every large enough $\kappa$" can be replaced by ``for arbitrarily large regular $\kappa < \mu$", or even:
``for arbitrarily large strong limit $\mu' < \mu$, for every large enough regular $\kappa < \mu'$".

A consequence is (see more on this principle in \cite{Sh:667}).

\begin{conclusion}\label{0.3} 
    1) If $\lambda,\kappa,\chi,S,\langle
    M_\delta \colon \delta \in S \rangle$ are as in Theorem \ref{0.1} and  $\gamma < \kappa \vee  \gamma < \chi$, \underline{then} we can
    find $\bar \eta = \langle \eta_\delta:\delta \in S \rangle,\eta_\delta$ is an increasing sequence of length $\kappa$ of ordinals $< \delta$ with limit $\delta$ such that:

    \begin{enumerate}
        \item[$\boxplus$] if $F \colon \lambda \rightarrow \gamma$, \underline{then} for stationarily many $\delta \in S$ we have $i < \kappa \Rightarrow F(\eta_\delta(2i)) = F(\eta_\delta(2i+1))$.
    \end{enumerate}
    
    2) In fact it is enough to find a $(\lambda,\kappa,\chi)$-\Bb-parameter
    $\langle C_\delta \colon \delta \in S \rangle$ which has the $(\gamma, 2^{|\gamma|^+})$-\Bb-property (see Definition \ref{d.5Y}(1)).
\end{conclusion}  

\begin{PROOF}{\ref{0.3}}
    1) By (2).
    
    2) Let $\tau = \{P_\zeta \colon \zeta < |\gamma|^+\},p_i$ a monadic predicate.  Let $\langle M_\delta \colon \delta \in S \rangle$ be as in Theorem \ref{0.1} for $\tau$. Let $F \colon \lambda \rightarrow \gamma$ and we let $M$ be the following $\tau$-model with universe $\lambda \colon P^M_\zeta =
    \{\alpha \colon \zeta = \pr(\zeta_0,\zeta_1), \zeta_0 < \zeta_1 <
    |\gamma|^+$ and $F(\alpha + \zeta_0) = F(\alpha + \zeta_1), \alpha \in
    M_\delta \Rightarrow \alpha + |\gamma|^+ < \min(M_\delta \backslash \alpha)$.  Without loss of generality $\otp(|M_\delta|) =
    \kappa,\delta \in S \Rightarrow |\gamma|^+/\delta$.  Let $(\eta_\delta(2,i),\eta_\delta(2i+1)) = (\alpha^\delta_i + \zeta^\delta_2,\alpha^\delta_i + \zeta^\delta_2)$ when $\alpha^\delta_i$ is the $i$-th member of $M_\delta,\alpha^\delta_i \in P$.
\end{PROOF} 

Note that if a naturally quite weak conjecture on pcf holds, then we can replace $\beth_\omega$ by $\aleph_\omega$ (and similar changes) in
most places but still need  $\lambda > \lambda_1 > \lambda_\delta,2^{\lambda_1} >
2^{\lambda_0} \ge 2^{\aleph_\omega}$ in the stepping-up lemma;  e.g. for every regular large
enough $\lambda$, for some $n$, we have \midia with $\lambda_0$ colours on $S^\lambda_{\aleph_n}$ - is helpful for abelian group theory.  (Under weak conditions we get this for $\aleph_1$, see \cite{Sh:F611}.

Even without this conjecture, it is quite hard to have the negation of this statement as it restricts $2^\mu$ and cf$(2^\mu)$ for every
$\mu$!

\begin{definition}\label{x.20.0109} 
    1) Let ``$S \subseteq \lambda$ has square" mean that $\lambda$ is regular uncountable, $S$ a stationary subset of $\lambda$ and: for some club $E$ of $\lambda$
    and $S', S \cap E \subseteq S' \subseteq \lambda$ there is a witness $\bar C = \langle C_\delta \colon \delta \in S' \rangle$ which means:
    
    \begin{enumerate}  
        \item[($\ast$)]    $\bar C = \langle C_\delta \colon \delta \in S' \rangle$ has  the \emph{square property} which means: 
        
        \begin{enumerate}  
            \item[(a)] $C_\delta \subseteq \delta$ is closed,
            
            \item[(b)]  $\delta \in S'$ limit $\Rightarrow \delta  = \sup(C_\delta)$,
            
            \item[(c)]   $\alpha \in C_\delta 
            \Rightarrow C_\alpha = \alpha \cap C_\delta$ and,
            
            \item[(d)]   $\delta \in S \Rightarrow \otp(C_\delta)  < \delta$. 
        \end{enumerate} 
    \end{enumerate}  
    
    1A) We say ``\emph{$S \subseteq \lambda$ has a strict square}" if we add: 

    \begin{enumerate}
        \item[(e)]  $\delta \in S \cap E \Rightarrow \otp(C_\delta) = \cf(\delta)$.
    \end{enumerate}
    
    2)  We say ``\emph{$S \subseteq \lambda$ has a $(\le \delta^*)$ square}" if we
    replace (d) above by:

    \begin{enumerate}
        \item[(d)] $\delta \in S \Rightarrow \otp(C_\delta) \le \delta^*$.
    \end{enumerate}

    If we replace square by $(\le \delta^*)$-square this means we add $\otp(C_\alpha) \le \delta^*$. 
    
    3) We say $(\le \kappa,\Theta)$-square if above the demand $\alpha \in C_\delta \Rightarrow C_\alpha = \alpha \cap C_\delta$ is restricted to the case cf$(\alpha) \in \Theta$. 
    
    4) $S^\lambda_\kappa = \{\delta < \lambda \colon \cf(\delta) =
    \kappa\}$.  
\end{definition} 

On $\hat{I}[\lambda]$ see \cite[\S1]{Sh:420}.

\begin{definition}\label{0.5}  
    1) For $\lambda$ regular
    uncountable let $\hat{I}[\lambda]$ be the family of sets $S \subseteq
    \lambda$ which have a witness $(E,\bar{\mathscr  P})$ for $S \in
    \hat{I}[\lambda]$, which means:

    \begin{enumerate}
        \item[$(\ast)$]  $E$ is a club of $\lambda, {\mathscr  P} = \langle {\mathscr  P}_\alpha \colon \alpha < \lambda \rangle,{\mathscr  P}_\alpha \subseteq {\mathscr  P}(\alpha), |{\mathscr  P}_\alpha| < \lambda$, and for every $\delta \in E \cap S$ there is an unbounded subset $C$ of $\delta$ of order type $< \delta$ such that $\alpha \in C \Rightarrow C \cap \alpha \in {\mathscr  P}_\alpha$.
    \end{enumerate}
 
    2) Equivalently there is a pair $(E,\bar a), E$ a club of $\lambda, \bar
    a = \langle a_\alpha \colon \alpha < \lambda \rangle, a_\alpha \subseteq
    \alpha,\beta \in a_\alpha \Rightarrow a_\beta = a_\alpha \cap \beta$
    and $\delta \in E \cap S \Rightarrow \delta = \sup(a_\delta) > \otp(a_\delta)$ (or even $\delta= \sup(a_\delta)$, $\otp(a_\delta) =
    \cf(\delta) < \delta$.
\end{definition} 

\begin{notation}\label{0.6}  
    1) Let $\lambda,$ $\mu,$ $\kappa, $ $\chi, \theta, $ $\sigma$ denote cardinals, $\Theta$ a set of cardinals, let $\alpha, \beta,\gamma, \delta, \varepsilon, \zeta, I, j, \Upsilon$ denote ordinals but $\Upsilon \ge 2$ and $\theta \ge 2$ and $\delta$ is a
    limit ordinal if not said otherwise. 
\end{notation} 

Clearly Claim \ref{d.6} and some results later continue \emph{Enge\l king Kar\l owicz} (see \cite{EK}) which continue works on the density of some topological spaces, see more history in \cite{Sh:420}.

\newpage

\section{Super black box: Sufficient conditions for {\rm cf}$(2^\mu)$}\label{1} 

As background recall: 

\begin{definition}\label{d.1} 
    1) For a regular uncountable $\lambda$, an ordinal $\Upsilon \ge 2$ and cardinal $\theta \ge 2$  we say $S \subseteq \lambda$ is \emph{$(\Upsilon,\theta)$-small} 
    \underline{if} it is {\bf F}-$(\Upsilon,\theta)$-small for some
    $(\lambda,\Upsilon,\theta)$-colouring $\mathbf F$, also called a 
    $\theta$-colouring of ${}^{\lambda >} \Upsilon$, that is: a function 
    $\mathbf F$ from ${}^{\lambda >}\Upsilon$ to $\theta$ where:
    
    2)  We say $S \subseteq 
    \lambda$ is \emph{$\mathbf F$-$(\Upsilon,\theta)$-small} \underline{if}
    ($\mathbf F$ is a $\theta$-colouring of ${}^{\lambda >} \Upsilon$) and 

    \begin{enumerate}
        \item[$(*)^\theta_{\mathbf F,S}$]  for every  $\bar c \in {}^S \theta$ for some $\eta \in {}^\lambda \Upsilon$ the set $\{\delta \in S \colon \mathbf F(\eta \, {\rest} \, \delta) = c_\delta\}$ is not stationary. 
    \end{enumerate}
    
    3)   Let $D^{\text{M}}_{\lambda, \Upsilon,\theta} 
     \coloneqq \{A \subseteq \lambda \colon \lambda \setminus
    A$ is $(\Upsilon, \theta)$-small$\}$, we call it the \emph{$(\Upsilon,\theta)$-middle diamond filter on $\lambda$}.
    If $\theta = 2$ we call this \emph{weak diamond filter} and may write $D^{\mathrm{wd}}_{\lambda,\Upsilon}$. 
    
    4) In part (3), 

    \begin{enumerate}
        \item[(a)]  We omit $\Upsilon$ if $\Upsilon = \theta$, writing  $\theta$ instead of $(\Upsilon, \theta),$

        \item[(b)] We may omit $\theta$ from ``$S$ is $\mathbf F$-$\theta$-small", and $(*)^\theta_{\mathbf F, S}$ when clear from the context. 
    \end{enumerate}
\end{definition} 

\begin{claim}\label{d.2}  
    1) If $\lambda = {\text {\rm cf\/}}(\lambda) > \aleph_0,\lambda \ge \Upsilon \ge 2$ and
    $\lambda \ge \theta$ \underline{then} 
    $D^{\text{\rm md\/}}_{\lambda,\Upsilon,\theta}$ is 
    a normal filter on $\lambda$. 
    
    2) Assume $S \subseteq \lambda,2 \le \Upsilon_\ell \le \lambda,2 \le
    \theta \le \lambda$ for $\ell = 1,2$.  \underline{then} $S$ is
    $(\Upsilon_1,\theta)$-small iff $S$ is $(\Upsilon_2,\theta)$-small
    \footnote{so $\Upsilon$ seems immaterial, but below when using $\bar C$
    this is not so (see Definition \ref{d.5}).}.
    
    3) In fact $\Upsilon_\ell \in [2,2^{< \lambda}]$ is enough. 
\end{claim} 

\begin{PROOF}{\ref{d.2}}
    Straightforward.  
\end{PROOF}

Recall (very relevant here).

\begin{definition}\label{d.3}  
    1) Let 
    $\chi^{<\sigma>}  = \min \{|{\mathscr  P}| \colon {\mathscr  P} \subseteq [\chi]^\sigma$
    and every $A \in [\chi]^\sigma$ is included in the union of  $< \sigma$ members of ${\mathscr  P}\}$. 
    
    2) Let $\chi^{[\sigma]} = \min\{|{\mathscr  P}| \colon {\mathscr  P} \subseteq
    [\chi]^\sigma$ and every $A \in [\chi]^\sigma$ is equal to the union of $<
    \sigma$ members of ${\mathscr  P}\}$.
    
    3) $\rm{trp}_{\sigma}(\chi) = \chi^{<\sigma>_{\text{tr}}} = \sup\{|\text{lim}_\sigma(t)| \colon t$ is a tree 
    with $\le \chi$ nodes and $\sigma$ levels$\}$, where $\rm{trp}$ stand for tree power, $\chi^{\langle \sigma \rangle_{\tr}}$ is an older notation. 
    
    4) For $J$ an ideal on $\theta$ let
    $\mathbf U_J(\mu) = \min\{|{\mathscr  P}| \colon {\mathscr  P} \subseteq
    [\mu]^\theta$ and for every $f \colon \dom(J) \rightarrow \mu$ for some $u \in {\mathscr  P}$ we have $\{i \in \dom(J) \colon f(i) \in u\} \ne
    \emptyset \mod J\}$; we say an ideal $J$ is an ideal on a set $ \dom(J)$ of cardinality $\theta$.  If $J = [\theta]^{< \theta}$ we may write $\theta$ instead of $J$. 
\end{definition} 

\begin{remark}\label{d.4} 
    1) On $ \trp_{\sigma}(\chi) = \chi^{<\sigma>_{\text{tr}}}$ 
    see \cite{Sh:589}, on $\chi^{<\sigma>}$ see \cite{Sh:430}, on $\chi^{<\sigma>},\chi^{[\sigma]}$ see there and in \cite{Sh:460}, \cite{Sh:829} but no real knowledge of these references is assumed here. 
    
    2) Note that when $\Upsilon = \lambda$ only 
    $\mathbf F \, {\rest} \, \{{}^\delta \delta \colon \delta \in S\}$ is relevant
    in Definition \ref{d.1} because for any $\eta \in {}^\lambda
    \lambda$ for some club $E$ of $\lambda$ we have $\delta \in E \Rightarrow \eta \, {\rest} \, \delta \in {}^\delta \delta$.  
    So we may restrict the colourings to this set.
    
    3) In Definition \ref{d.1}, if $\mu < \lambda = \cf(\lambda),2^\mu = 2^\lambda$ even for $\theta = \Upsilon = 2$, the set $\lambda$ is small (as noted by Uri Abraham, see \cite{Sh:65}). 
\end{remark} 

Below we may usually restrict ourselves to $\Upsilon = \lambda$.  

\begin{definition}\label{d.5}  
    1)  We say $\bar C$ is a $\lambda$-\Bb-parameter \underline{if}:
    
    \begin{enumerate} 
        \item[(a)] $\lambda$ is regular uncountable,
        
        \item[(b)] $S$ is a stationary subset of $\lambda$,
        
        \item[(c)] $\bar C = \langle C_\delta \colon \delta \in S
        \rangle$ and\footnote{we can omit the requirement $\delta =
        \sup(C_\delta)$ but then we trivially have \midia \ if $\lambda = \lambda^{< \kappa}$ when $\kappa =
        \cup\{|C_\delta|^+ \colon \delta \in S\}.$}  $C_\delta \subseteq \delta = \sup(C_\delta)$.
    \end{enumerate} 
    
    1A)  We say $\bar C$ is a \emph{$(\lambda,\kappa,\chi)$-\Bb-parameter} \underline{if} in addition $(\forall \delta \in S)(\forall \alpha \in C_\delta)
    [\cf(\delta) = \kappa \wedge |C_\delta \cap \alpha| < \chi]$.
    
    2)  We say that {\bf F} is a \emph{$(\bar C, \Upsilon,\theta)$-colouring} if:
    $\theta \ge 2,\Upsilon \ge 2, \bar C = 
    \langle C_\delta \colon \delta \in S \rangle$ is a $\lambda$-\Bb-parameter and {\bf F} is a function  from ${}^{\lambda>} \Upsilon$ to $\theta$ such that:

    \begin{enumerate}
        \item[$\bullet$]  if $\delta \in S$ and $\eta_0,\eta_1 \in {}^\delta \Upsilon$ and $\eta_0 \, {\rest} \, C_\delta = \eta_1 \, {\rest} \, C_\delta$ then {\bf F}$(\eta_0) =$ {\bf F}$(\eta_1)$.
    \end{enumerate}

    2A)  If $\Upsilon = \theta$ we may omit it and write  $(\bar C,\theta)$-colouring.
    
    2B)  In part (2) we can replace {\bf F} by $\bar F = \langle F_\delta:\delta
    \in S \rangle$ where $F_\delta:{}^{(C_\delta)} \Upsilon \to \theta$
    such that $\eta \in {}^\delta \Upsilon \wedge \delta \in S 
    \rightarrow \mathbf F(\eta) =
    F_\delta(\eta \, {\rest} \, C_\delta)$.  If $\Upsilon = \lambda$ we may
    use $F_\delta \colon {}^{(C_\delta)} \delta \rightarrow \theta$ (as we can
    ignore a non-stationary $S' \subseteq S$); so abusing notation we may write
    $\mathbf F(\eta \, {\rest} \, C_\delta)$.  Similarly if cf$(\Upsilon) =
    \lambda,\langle \Upsilon_\alpha \colon \alpha < \lambda \rangle$ is
    increasing continuous with limit $\Upsilon$ we may restrict ourselves 
    to ${}^{(C_\delta)}(\Upsilon_\delta)$.

    3)  Assume $\mathbf F$ is a 
    $(\bar C,\Upsilon,\theta)$-colouring, where $\bar C$ is a  $\lambda$-\Bb-parameter. 
    
    We say $\bar c \in {}^S \theta$ (or $\bar c \in {}^\lambda \theta$)
    is an $\mathbf F$-\midia sequence, or an $\mathbf F$-\Bb-sequence if:
    
    \begin{enumerate}    
        \item[$(\ast)$]  for every $\eta \in {}^\lambda \Upsilon$, the 
        set $\{\delta \in S \colon \mathbf F(\eta \, {\rest} \, \delta) = c_\delta\}$  is a stationary subset of $\lambda$.
    \end{enumerate} 
    
    We also may say that $\bar c$ is an $(\mathbf F,S)$-\Bb-sequence or an
    $(\mathbf F,\bar C)$-\Bb-sequence; note that $\Upsilon$ can be
    reconstructed from $\mathbf F$.

    3A)  We say $\bar c \in {}^S \theta$ is a $D-{\mathbf F}$-\Bb-sequence or
    $D-(\mathbf F, \bar C)$-\Bb-sequence
    \underline{if} $D$ is a filter on $\lambda$ (no harm to assume that $S$
    belongs to $D$, or at least $S \in D^+$) and
    
    \begin{enumerate} 
        \item[$(**)$]  for every $\eta \in {}^\lambda \Upsilon$ we have $$ \{\delta \in S \colon \mathbf F(\eta \, {\rest} \, \delta) = c_\delta\} \ne \emptyset \mod D.$$
    \end{enumerate} 
    
    4)  We say that $\bar C$ is a good $(\lambda,\lambda_1)$-\Bb-parameter,
    \underline{if} $\bar C$ is a $\lambda$-\Bb-parameter and for  every $\alpha < \lambda$ we have $\lambda_1 > |\{C_\delta \cap
    \alpha \colon \delta \in S \, \wedge \,\alpha \in C_\delta\}|$.  If $\lambda_1 =
    \lambda$ we may write ``a good $\lambda$-\Bb-parameter".  Similarly for a
    $(\lambda, \lambda_1, \kappa, \chi)$-good-\Bb-parameter and a $(\lambda,\kappa,\chi)$-good-\Bb-parameter (i.e. adding that $\bar C$ is  a $(\lambda,\kappa,\chi)$-\Bb-parameter.  
    
    4A) We say that $\bar C$ is a 
    weakly good $\lambda$-\Bb-parameter \underline{if} it is a $\lambda$-\Bb-parameter and for every $\alpha < \lambda$ we have $\lambda > |\{C_\delta \cap \alpha \colon \delta \in S \, \wedge \,\alpha \in \mathrm{nacc}(C_\delta)\}|$.  Similarly, we define a weakly good $(\lambda, \kappa, \chi)$-\Bb-parameter and a weakly good $(\lambda, \lambda_1, \kappa, \chi)$-\Bb parameter.  

    4B) We say $\bar C$ is a good$^+-\lambda$-\Bb-parameter \underline{if}
    $\bar C = \langle C_\delta \colon \delta \in S \rangle$ is a $\lambda$-\Bb-parameter and $\alpha \in C_{\delta_1} \cap C_{\delta_2}
    \Rightarrow C_{\delta_1} \cap \alpha = C_{\delta_2} \cap \alpha$. Similarly we define a good$^+-(\lambda,\kappa,\chi)$-\Bb-parameter. 
    
    4C) We define a weakly good$^+-\lambda$-\Bb-parameter or $(\lambda,\kappa,\chi)$-\Bb-parameter analogously, i.e.,  $\alpha \in \mathrm{nacc}(C_{\delta_1}) \cap \mathrm{nacc}(C_{\delta_2})
    \Rightarrow C_{\delta_1} \cap \alpha = C_{\delta_2} \cap \alpha$. 

    4D) For a set $\Theta$ of regular cardinals we say that $\bar C =
    \langle C_\delta \colon \delta \in S \rangle$ is $\Theta$-closed \underline{if} each $C_\delta$ is; where a set $C$ of ordinals is $\Theta$-closed if for
    every ordinal $\alpha$, cf$(\alpha) \in \Theta \, \wedge \,\alpha = \sup(C \cap \alpha) < \sup(C) \Rightarrow \alpha \in C$. 
    
    4E) We say $\langle (C_\delta,C'_\delta) \colon \delta \in S \rangle$ is a
    good $(\lambda,\kappa,\chi)$-\Bb$_*$-parameter if $\langle    C_\delta \colon \delta \in S \rangle$ is a
    $(\lambda,\kappa,\chi)$-\Bb-parameter, $C'_\delta \subseteq \delta, \ \sup(C'_\delta) = \delta$ and for every $\alpha < \lambda$ the
    set $|\{C_\delta \cap \alpha \colon \alpha \in C'_\delta,\delta \in S\}|$ has
    cardinality $< \lambda$.
\end{definition} 

\begin{definition} \label{d.5Y}  
    1) We say that $\bar C$ (or $(\lambda,\bar C)$)  has the \emph{$(D,\Upsilon,\theta)$-\Bb-property} \underline{if}:
    
    \begin{enumerate}  
        \item[(a)]   $\bar C$ is a $\lambda$-\Bb-parameter,
        
        \item[(b)] $D$ is a filter on $\lambda,$
        
        \item[(c)]  if $\mathbf F$ is a $(\bar C,\Upsilon,\theta)$-colouring, 
    \end{enumerate}  

    \underline{then} there is a $D$-$\mathbf F$-\Bb-sequence, see Definition
    \ref{d.5}, part (3A).
    
    We may omit $D$ \underline{if} it is the club filter on $\lambda$. 
    
    2)  We say that $\bar C$ has the \emph{$(D,\theta)$-\Bb-property} as in part
    (1) only that $\mathbf F$ is a $(\bar C,\theta,\theta)$-colouring.  We may
    omit $D$ if $D$ is the club filter.  We may omit $\theta$ if
    $\theta=2$. 
    
    3) We say that $S \subseteq \lambda$ has the
    \emph{$(\Upsilon,\theta)$-\Bb-property} when some 
    $\bar C = \langle C_\delta \colon \delta \in S \rangle$ has, and then we say the
    $\theta$-\Bb-property if $\Upsilon = \theta$.  In Definition \ref{d.5}(2),(2A) we replace ``colouring" by
    ``colouring$_2$" if we have all $\eta \colon {}^{\omega >} \lambda \rightarrow \Upsilon$ (or $\eta \colon {}^{\omega >}(C_\delta) \rightarrow \Upsilon$); and let ``colouring$_1$" be the original notion (so in Definition \ref{d.5}(3),(3A) we have new meaning when we use such
    $\mathbf F$'s).

    4) We write \emph{\Bb$_2$-property} \underline{if} we use colouring$_2$ \, $\mathbf F$
    and we write \Bb$_1$-property for the \Bb-property.
\end{definition} 

\begin{remark}\label{d.5Z}  
    We may consider replacing $C_\delta$ by
    $(C_\delta,D_\delta), D_\delta$ a filter on $C_\delta$,  that is $F_\delta(\eta)$ would depend on $(\eta \, {\rest} \, C_\delta)/D_\delta$ only.  For the time being it does not make a real difference.
\end{remark} 

\begin{claim}\label{d.5A} 
    1) The restrictions in 
    Definition \ref{d.5}{\text{\rm (2B)\/}} are O.K.  That is in Definition \ref{d.5}(2B), the first sentence is obvious; for the second sentence (if $\Upsilon = \lambda$) recall that for every $\eta \in {}^\lambda \Upsilon$, for
    some club $E$ of $\lambda,\delta \in E \Rightarrow \eta \, {\rest} \,
    \delta \in {}^\delta \delta$; and similarly in the third sentence.
     
    2) In Definition \ref{d.5}{\text{\rm (4)\/}}, if $\lambda_1 > \lambda$ \underline{then} goodness holds trivially. 
    
    3) If $\bar C = \langle C_\delta \colon \delta \in S \rangle$ is a weakly
    good $(\lambda, \kappa, \chi)$-parameter, \underline{then} for some club $E$ of
    $\lambda$ there is $\bar C' = \langle C'_\delta \colon \delta \in S \cap E
    \rangle$ which is a good$^+ \, (\lambda, \kappa, \chi)$-parameter. 
    
    4) If $\bar C = \langle C_\delta \colon \delta \in S \rangle$ has the $(D,\Upsilon,\theta)$-\Bb-property, and $\bar C' = \langle C'_\delta \colon \delta \in S \rangle,C'_\delta \subseteq C_\delta$,
    \underline{then} also $\bar C'$ has the $(D, \Upsilon, \theta)$-\Bb-property.
\end{claim} 

\begin{PROOF}{\ref{d.5A}}
    1), 2) Easy.
    
    3) Note that the $C'_{\delta}$'s are not required to be a club of $\delta$,  see \cite[\S1]{Sh:420} but we elaborate.
    
    Let $\langle A_\beta \colon \beta < \lambda \rangle$ list without repetition
    $\{C_\delta \cap \alpha \colon \delta \in S$ and $\alpha \in \mathrm{nacc}(C_\delta)\}$ and be such that $A_\beta = A_\gamma \cap \alpha \Rightarrow A_\beta \in \{A_\varepsilon \colon \varepsilon \le \gamma\}$.  Let $E = \{\delta < \lambda$: if $\beta < \delta$ and $\alpha \in S$ and
    $\beta \in \mathrm{nacc}(C_\alpha)$ then $C_\alpha \cap \beta \in
    \{A_\gamma \colon \gamma < \delta\}\}$, clearly $E$ is a club of $\lambda$.
    
    Lastly, for $\delta \in S \cap E$ let 
    $C'_\delta = \{\gamma < \delta$: for some $\alpha \in \mathrm{nacc}(C_\delta)$ we have $C_\delta \cap \alpha = A_\gamma\}$.
    Now $\langle C'_\delta \colon \delta \in S \cap E \rangle$ is as required. 
    
    4) Should be clear.
\end{PROOF}    

From the following definition the case $\Sep(\mu,\theta)$ is used in the main claim below.

\begin{definition}\label{d.8}
    1) $\Sep(\chi,\mu,\theta_0,\theta,\Upsilon)$ means that for some $\bar f$:
    
    \begin{enumerate} 
        \item[(a)] $\bar f = \langle f_\varepsilon \colon \varepsilon < \chi
        \rangle,$
        
        \item[(b)] $f_\varepsilon$ is a function from  ${}^\mu(\theta_0)$ to $\theta,$
        
        \item[(c)]  for every $\varrho \in {}^\chi \theta$ the set $\{\nu \in {}^\mu(\theta_0)$: for every $\varepsilon < \chi$ we have $f_\varepsilon(\nu) \ne \varrho(\varepsilon)\}$ has cardinality $< \Upsilon$.
    \end{enumerate}  
    
    2)  We may omit $\theta_0$ if $\theta_0 = \theta$.  We write $\Sep(\mu,\theta)$ if for some $\Upsilon = \cf(\Upsilon) \le 2^\mu$ we have $\Sep(\mu, \mu,\theta, \theta, \Upsilon)$ and $\Sep(<\mu,\theta)$ if for some $\Upsilon = \cf(\Upsilon) \le 2^\mu$ and $\sigma < \mu$ we have $\Sep(\sigma,\mu,\theta,\theta,\Upsilon)$.
\end{definition}

\begin{claim}[Main claim]\label{d.6}  
    Assume, 
    
    \begin{enumerate}   
        \item[(a)]   $\lambda = {\text {\rm cf\/}}(2^\mu),D$ is a $\mu^+$-complete filter on $\lambda$ extending the club filter, 
        
        \item[(b)]   $\kappa = {\text {\rm cf\/}}(\kappa) < \chi \le \lambda$, 
        
        \item[(c)]      $\bar C = \langle C_\delta \colon \delta \in S \rangle$ is a $(\lambda,\kappa,\chi)$-good-{\rm \Bb}-parameter, or just weakly good or  just $\langle (C_\delta,C'_\delta) \colon \delta \in S \rangle$ is a good $(\lambda,\kappa,\chi)$-{\rm \Bb}$_*$-parameter and $S \in D,$ 
        
        \item[(d)]   $2^{< \chi} \le 2^\mu$ and $\theta \le \mu,$
        
        \item[(e)]   $\alpha < 2^\mu \Rightarrow 
        \rm{trp}_{\kappa}( \vert \alpha \vert) < 2^\mu,$
        
        \item[(f)] $\Sep(\mu,\theta)$,
    \end{enumerate} 
    
    \underline{then} $\bar C$ has the $(D,2^\mu,\theta)$-{\rm \Bb}-property
    {\text{\rm (\/}}recall that this
    means that the number of colours is $\theta$ not just 2{\text{\rm )\/}}. 
\end{claim} 

We may wonder if clause (f) of the assumption is reasonable; the following Claim gives some sufficient conditions for clause (f) of Claim \ref{d.6} to hold.

\begin{claim}\label{d.7}  
    Clause {\text{\rm (f)\/}} of 
    Claim \ref{d.6} holds, i.e., $\Sep(\mu,\theta)$ holds, \underline{if} at least one of the following holds:
    
    \begin{enumerate} 
        \item[(a)] $\mu = \mu^\theta,$
        
        \item[(b)] $\mathbf U_\theta(\mu) = \mu + 2^\theta \le \mu$,
        
        \item[(c)]  $\mathbf U_J(\mu) = \mu$ where for some $\sigma$ we have $J = [\sigma]^{< \theta},\theta \le \sigma,\sigma^\theta \le \mu,$
        
        \item[(d)]  $\mu$ is a strong limit of cofinality $\ne \theta = {\text {\rm cf\/}}(\theta) < \mu,$
        
        \item[(e)]  $\mu \ge \beth_\omega(\theta)$.
    \end{enumerate} 
\end{claim} 

We shall prove this later.

\begin{remark} \label{d.7B} 
    1) Note that $\lambda = \lambda^{< \lambda}
    = \chi,\kappa = \cf(\kappa) < \lambda,S = S^\lambda_\kappa, C_\delta = \delta$ is a particular case (of interest) of the assumptions (a)-(b) of Claim \ref{d.6} and that $\langle C_\delta \colon \delta
    \in S^\lambda_\kappa \rangle$ is a $(\lambda,\kappa,\lambda)$-good$^+$- 
    \Bb-parameter. 
    
    2)  We can replace the requirement ``$S \subseteq S^\lambda_\kappa$"
    by (implicit in (c) of Claim \ref{d.6})
    
    \begin{enumerate} 
        \item[$(\ast \ast \ast)$]  $[\delta \in S \Rightarrow \delta$ is a regular  cardinal].
    \end{enumerate} 
    
    Then clause (e) of Claim \ref{d.6} is replaced by ``$\delta \in S \Rightarrow \trp_{\cf(\delta)}(\vert \delta \vert) < 2^\mu$". 
    
    3)  We may consider 
    finding a guessing sequence $\bar c$ simultaneously for many colourings.  In the notation of 
    \cite[AP,\S1]{Sh:f} this is the same as having large $\mu(0)$.  We may reconsider the problem of $\lambda^+$-entangled linear
    ordering where $\lambda = \lambda^{\aleph_0}$, see \cite{Sh:462}. 
    
    4) For $\theta = 2$, (i.e. weak diamond) the proof below can be simplified: we can
    forget ${\mathscr  F}$ and choose $\varrho_\delta \in {}^\mu \theta \backslash \{\rho_\eta \colon \eta \in {\mathscr  P}_\delta\}$ such that for some $i < \mu$ the sequence $\langle \varrho_\delta(i) \colon \delta \in S \rangle$ is a weak diamond sequence for $\mathbf F$. 
    
    5) Note that the ``minimal" natural choices for getting $\theta$-MD property are $\kappa = \aleph_0,\theta = 2^{\aleph_0},\mu =
    2^\theta,\lambda = \cf(2^\mu)$.
    
    6) Note that in the main claim (Claim \ref{d.6}) if $2^\mu$ is regular, $\Upsilon = 2^\mu$ is allowed.
\end{remark} 

Now we prove Claim \ref{d.6}: 

\begin{PROOF}{\ref{d.6}}
    In clause (c) we can assume $\langle
    (C_\delta,C'_\delta) \colon \delta \in S \rangle$ is a good $(\lambda,\kappa,\chi)$-\Bb$_*$-parameter (in the other cases we choose $C'_\delta = \nacc(C_\delta)$ and use \ref{d.5A}(3)). 
    
    Let $\mathbf F$ be a given $(\bar C,2^\mu,\theta)$-colouring.
    
    By assumption (f) we have $\Sep(\mu,\theta)$ which means (see Definition
    \ref{d.8}(2)) that for some $\Upsilon = \cf(\Upsilon) \le
    2^\mu$ we have $\Sep(\mu,\mu,\theta, \theta,\Upsilon)$.

    Let ${\mathscr  F} = \{f_\xi \colon \xi < \mu\}$ where ${\mathscr  F}$
    exemplify $\Sep(\mu,\mu,\theta,\theta,\Upsilon)$, see Definition \ref{d.8}(1)
    and for $\rho \in {}^\mu \theta$ let 
    $\Sol_\rho \coloneqq \{\nu \in {}^\mu \theta \colon \text{ for every } \varepsilon 
    < \mu$ we have $\rho(\varepsilon) \ne f_\varepsilon(\nu)\}$ where 
    $\Sol$ stands for solutions, so by clause (c) of the Definition
    \ref{d.8}(1) of $\Sep$ it follows that: 
    
    \begin{enumerate} 
        \item[$(\ast)$]   $\rho \in {}^\mu \theta 
        \Rightarrow |\Sol_\rho| < \Upsilon$.
    \end{enumerate} 
    
    Let $h$ be an increasing continuous function from $\lambda$ into $2^\mu$ with unbounded range; if $2^{\mu}$ is regular, we can choose $h = \id_\lambda$   Recall that $\lambda = \cf(2^\mu) \ge 2^{< \chi}$.
    Let $\cd$ be a one-to-one function from ${}^\mu (2^\mu)$ onto $2^\mu$ such that: $$
    \alpha = \cd(\langle \alpha_\varepsilon \colon \varepsilon < \mu
    \rangle) \Rightarrow \alpha \ge \sup\{\alpha_\varepsilon \colon \varepsilon < \mu\}, $$ and let the function $ \cd_i \colon 2^\mu \rightarrow 2^\mu$ for $i < \mu$ be 
    such that $ \cd_i(\cd(\langle \alpha_\varepsilon \colon  \varepsilon < \mu \rangle)) = \alpha_i$.
    
    Let ${\mathscr  C} = \{C_\delta \cap \alpha \colon \delta \in S    {\, \wedge \, } \alpha \in C'_\delta\}$, (recall that trivially
    $|{\mathscr  C}| \le \lambda \le 2^\mu$)
    and let $T= \{\eta \colon $ for some $C \in {\mathscr  C}$ (hence $C$ is of cardinality $< \chi$), we have $\eta \in
    {}^C (2^\mu)\}$.  By assumptions (c) + (d) clearly $|T|  = \dsize \sum_{C \in {\mathscr  C}} 2^{\mu \cdot |C|} \le |{\mathscr  C}| \times (2^\mu)^{< \chi} \le 2^\mu$ but $|T| \ge 2^\mu$ hence $|T| =  
    2^\mu$.  So let us list $T$ possibly with repetitions as $\{\eta_\alpha \colon \alpha < 2^\mu\}$ such that $[\alpha < h(\beta) \Rightarrow \ran(\eta_\alpha) \subseteq h(\beta)]$ 
    and let $T_{< \alpha} = \{ \eta_\beta \colon \beta < h(\alpha)\}$ for $\alpha
    < \lambda$.
    For $\delta \in S$ let ${\mathscr  P}_\delta \coloneqq  \{\eta \colon \eta$ is a function from
    $C_\delta$ into $h(\delta)$ such that for 
    every $\alpha \in C'_\delta$ we have $\eta \, {\rest} \,
    (C_\delta \cap \alpha) \in T_{< \delta}\}$, clearly $\eta \in 
    {\mathscr  P}_\delta \Rightarrow \eta \in {}^{(C_\delta)} h(\delta)$.  Clearly
    $\{\eta \, {\rest} \, (C_\delta \cap \alpha) \colon \eta \in {\mathscr  P}_\delta,
    \alpha \in C'_\delta\}$ is
    naturally a tree with ${\mathscr  P}_\delta$ the set of $\otp(C_\delta)$-branches and cf$(\otp(C_\delta)) = \kappa$
    and set of nodes $\subseteq T_{< \delta}$ hence with $< 2^\mu$ nodes
    hence by assumption (e)
    
    \begin{enumerate}
        \item[$(\ast \ast)$] $|{\mathscr  P}_\delta| < 2^\mu$.
    \end{enumerate} 
    
    For each $\eta \in {\mathscr  P}_\delta$
    and $\varepsilon < \mu$ we define $\nu_{\eta,\varepsilon}
    \in {}^{C_\delta}(h(\delta))$ by $\nu_{\eta,\varepsilon} (\alpha) = 
    \cd_\varepsilon (\eta(\alpha))$ for $\alpha \in C_\delta$. 
    
    Now for $\eta \in {\mathscr  P}_\delta$, clearly $\rho_\eta \coloneqq 
    \langle \mathbf F (\nu_{\eta,\varepsilon}) \colon \varepsilon < \mu \rangle$ 
    belongs to ${}^\mu \theta$. Clearly
    $\{\rho_\eta \colon \eta \in {\mathscr  P}_\delta\}$ is a subset of ${}^\mu \theta$ 
    of cardinality $\le |{\mathscr  P}_\delta|$ which as said above is  $< 2^\mu$.  Recall the definition of $\Sol_\rho$  from the beginning of the proof and remember $(*)$ above, so $\eta \in {\mathscr  P}_\delta \Rightarrow \Sol_{\rho_\eta}$ has cardinality $< \Upsilon = \cf(\Upsilon)
    \le 2^\mu$.  Hence we can find $\varrho^*_\delta \in {}^\mu \theta \setminus
    \cup \{\Sol_{\rho_\eta} \colon \eta \in {\mathscr  P}_\delta\}$.  
    
    [Why?  if $\Upsilon < 2^\mu$ then $|\cup \{\Sol_{\rho_\eta} \colon \eta
    \in {\mathscr  P}_\delta\}| \le \Sigma\{|\Sol_{\rho_\eta}| \colon \eta \in
    {\mathscr  P}_\delta\} \le \Upsilon \times |{\mathscr  P}_\delta| < 2^\mu =
    \theta^\mu = |{}^\mu \theta|$ and if $\Upsilon = 2^\mu$ then $2^\mu$
    is regular, $|{\mathscr  P}_\delta| < 2^\mu,\eta \in {\mathscr  P}_\delta
    \Rightarrow |\Sol_{\rho_\eta}| < 2^\mu$ and again
    $|\cup\{\Sol_{\rho_\eta}:\eta \in {\mathscr  P}_\delta\}| < 2^\mu =
    |{}^\mu \theta|$.] 
    
    Let $\varepsilon < \mu$.  Recall
    that $\varrho^*_\delta \in
    {}^\mu \theta$ for $\delta \in S$ and $f_\varepsilon$ from  the list of ${\mathscr  F}$ is a function from ${}^\mu \theta$
    to $\theta$ so $f_\varepsilon(\varrho^*_\delta) < \theta$.  Hence we can consider the sequence $\bar c^\varepsilon = \langle
    f_\varepsilon(\varrho^*_\delta) \colon \delta 
    \in S \rangle \in {}^S \theta$ as a candidate for being an {\bf
    F}-\Bb-sequence.  If one of them is, we are done. So assume towards a
    contradiction that for each $\varepsilon< \mu$ there is a sequence
    $\eta_\varepsilon \in {}^\lambda (2^\mu)$ that exemplifies
    the failure of $\bar c^\varepsilon$ to be an $\mathbf F$-\Bb-sequence,  so there is a $E_\varepsilon \in D$ such that:

    \begin{enumerate}
        \item[$\boxtimes_1$] $\delta \in S \cap E_\varepsilon \Rightarrow \mathbf F (\eta_\varepsilon \, {\rest} \, C_\delta) \ne f_\varepsilon(\varrho^*_\delta)$.
    \end{enumerate}

    Define $\eta^* \in {}^\lambda(2^\mu)$ by $\eta^*(\alpha) =  \cd(\langle \eta_\varepsilon(\alpha) \colon \varepsilon 
    < \mu \rangle)$.  Now as $\lambda$ is regular uncountable it follows
    by choice of $h$ that $E \coloneqq \{\delta < \lambda \colon$ for every
    $\alpha < \delta$ we have $\eta^* (\alpha) < h(\delta)$ and if $\delta' \in S,\alpha \in C'_{\delta'} \cap \delta$ and $C'' = C_{\delta'}  \cap \alpha$ then $\eta^* \, {\rest} \, C'' \in T_{< \delta}\}$ is 
    a club of $\lambda$ (see the choice of $T$ and $T_{<\delta}$, and recall as well 
    that by assumption (c) the sequence 
    $\langle (C_\delta,C'_\delta) \colon \delta \in S \rangle$ is good, see Definition \ref{d.5}(4) + (4D)).
    
    By clause (a) in the assumption of our Claim \ref{d.6}, the filter $D$ includes the clubs of $\lambda$ so clearly $E \in D$; also $D$ is $\mu^+$-complete hence
    $E^* \coloneqq \cap \{E_\varepsilon \coloneqq \varepsilon < \mu\} \cap E$ belongs to
    $D$.  Now choose $\delta \in E^* \cap S$, fixed for the rest of the
    proof;  clearly $\eta^* 
    \, {\rest} \, C_\delta \in
    {\mathscr  P}_\delta$; just check the definitions of ${\mathscr  P}_\delta$ and
    $E,E^*$.  Let $\varepsilon < \mu$.
    Now recall that $\nu_{\eta^* \, {\rest} \, C_\delta,\varepsilon}$
    is the function from $C_\delta$ to $h(\delta)$ defined by: $$
    \nu_{\eta^* \, {\rest} \, C_\delta,\varepsilon}(\alpha) =
    \cd_\varepsilon (\eta^*(\alpha)).$$
    
    But  by our choice of 
    $\eta^*$ clearly $\alpha \in C_\delta \Rightarrow \cd_\varepsilon
    (\eta^*(\alpha)) = \eta_\varepsilon (\alpha)$, so $$ \alpha \in 
    C_\delta \Rightarrow \nu_{\eta^* \, {\rest} \, C_\delta,\varepsilon}
    (\alpha) = \eta_\varepsilon (\alpha) \quad \text{ so } \quad
    \nu_{\eta^* \, {\rest} \, C_\delta,\varepsilon} = \eta_\varepsilon \, {\rest} \, C_\delta. $$
    
    Hence $\mathbf F (\nu_{\eta^* \, {\rest} \, C_\delta,\varepsilon})
    = \mathbf F (\eta_\varepsilon \, {\rest} \, C_\delta)$.  As $\delta \in E^* \subseteq E_\varepsilon$ clearly $\mathbf F
    (\eta_\varepsilon \, {\rest} \, C_\delta) \ne f_\varepsilon (\varrho^*_\delta)$ and as $\eta^* \, {\rest} \, C_\delta \in {\mathscr 
    P}_\delta$ clearly $\rho_{\eta^* \, {\rest} \, C_\delta} \in {}^\mu
    \theta$ is well defined.  Now easily $\rho_{\eta^* \, {\rest} \, C_\delta}
    (\varepsilon) = \mathbf F(\nu_{\eta^* \, {\rest} \, C_\delta,\varepsilon})$ by 
    the definition of $\rho_{\eta^* \, {\rest} \, C_\delta}$, so we have
    $\rho_{\eta^* \, {\rest} \, C_\delta}(\varepsilon) \ne
    f_\varepsilon(\varrho^*_\delta)$.
    
    As this holds for every $\varepsilon < \mu$ it follows  by the definition of $\Sol_{\rho_{\eta^* \, {\rest} \, C_\delta}}$ that $\varrho^*_\delta \in \Sol_{\rho_{\eta^* \, {\rest} \, C_\delta}}$. But $\eta^* \, {\rest} \, C_\delta \in {\mathscr  P}_\delta$ was proved
    above so $\rho_{\eta^* \, {\rest} \, C_\delta} \in \{\rho_\eta \colon \eta \in {\mathscr 
    P}_\delta\}$ hence $\varrho^*_\delta \in \Sol_{\rho_{\eta^*
    \, {\rest} \, C_\delta}} \subseteq \cup \{ \Sol_{\rho_\nu} \colon\nu \in
    {\mathscr  P}_\delta\}$
    whereas $\varrho^*_\delta$ has been chosen
    outside this set, a contradiction.  
\end{PROOF}    

We return to looking at $\Sep$ (see Definition \ref{d.8}) which appears in assumption (f) of Claim \ref{d.6}.  There are obvious monotonicity properties:

\begin{claim}\label{d.8Z} 
    If $\Sep(\chi,\mu,\theta_0,\theta_1,\Upsilon)$  holds and $\chi' \ge \chi,(\theta'_0)^{\mu'} \le (\theta_0)^\mu, \theta'_1 \le \theta_1$ and $\Upsilon' \ge \Upsilon$, \underline{then} $\Sep(\chi',\mu',\theta'_0,\theta'_1,\Upsilon')$ holds. 
\end{claim} 

\begin{PROOF}{\ref{d.8Z}}
    Let $h$ be a one-to-one function from
    ${}^{\mu'}(\theta'_0)$ into ${}^\mu(\theta_0)$.  Let
    $\bar f = \langle f_\varepsilon \colon \varepsilon < \chi \rangle$ exemplify $\Sep(\chi, \mu, \theta_0, \theta_1, \Upsilon)$.  For
    $\varepsilon < \chi'$ we define $f'_\varepsilon \colon {}^{(\mu')}(\theta'_0)
    \rightarrow \theta'_1$ by: $f'_\varepsilon(\nu) = f_\varepsilon(h(\nu))$
    if $\varepsilon < \chi \, \wedge \,  f_\varepsilon(h(\nu)) < \theta'_1$ and
    $f'_\varepsilon(\nu) = 0$ otherwise.
    Let $\varrho' \in {}^{(\chi')}(\theta'_1)$ and  let $A'_{\varrho'} = \{\nu \in {}^{(\mu')}(\theta'_0)$: 
    for every $\varepsilon < \chi'$ we have
    $f'_\varepsilon(\nu) \ne \varrho'(\varepsilon)\}$.  Let $\varrho \coloneqq
    \varrho' \, {\rest} \, \chi$ so $\varrho \in {}^\chi(\theta'_1)
    \subseteq {}^\chi(\theta_1)$, hence
    by the choice of $\langle f_\varepsilon \colon \varepsilon < \chi \rangle$
    the set $A_\varrho = \{\nu \in {}^\mu(\theta_0)$: for every $\varepsilon
    < \chi$ we have $f_\varepsilon(\nu) \ne \varrho(\varepsilon)\}$ has
    cardinality $< \Upsilon$.  If $\nu
    \in A'_{\varrho'}$ then $\nu \in {}^{\mu'}(\theta'_0)$ and $h(\nu) \in {}^\mu
    (\theta_0)$, and for $\varepsilon < \chi'$ we have $f'_\varepsilon(\nu)
    \ne \varrho'(\varepsilon)$ hence for $\varepsilon < \chi(\le \chi')$
    we have either $f_\varepsilon(h(\nu)) < \theta'_1 \, \wedge \,\varrho(\varepsilon) =
    \varrho'(\varepsilon) \ne f'_\varepsilon(\nu) = f_\varepsilon(h(\nu))$ or $f_\varepsilon(h(\nu)) \ge \theta'_1 \, \wedge \,\varrho(\varepsilon) =
    \varrho'(\varepsilon) < \theta'_1$ so $f_\varepsilon(h(\nu)) \ne
    \varrho(\varepsilon)$.  Hence $\nu \in A'_{\varrho'} \Rightarrow h(\nu) \in
    A_\varrho$ and as $h$ is a one-to-one function we have  $|A'_{\varrho'}| \le |A_\varrho| < \Upsilon \le \Upsilon'$ 
    and we are done.  
\end{PROOF}

To show that Claim \ref{d.6} has reasonable assumptions we should still
prove Claim \ref{d.7}: 

\begin{PROOF}{\ref{d.7}}
    Our claim gives sufficient conditions
    for $\Sep(\mu,\theta)$, i.e. $\Sep(\mu,\mu,\theta,\theta,\Upsilon)$ for some $\Upsilon = \cf(\Upsilon) \le 2^\mu$.  
    Clearly, if $\Upsilon < 2^\mu$ we can waive ``$\Upsilon$ regular" as $\Upsilon^+$ can serve as well.  Let $\theta_1 = \theta$.
    
    \underline{Case 1}:  $\mu = \mu^\theta,\Upsilon = \theta,\theta_0 \in
    [\theta,\mu]$ and we shall prove
    $\Sep(\mu, \mu, \theta_0, \theta_1,\Upsilon)$.  Let $ {\mathscr  F} $ be the set of $ f $ such that for some $ u, \bar{ \rho } $ we have:  
             
    \begin{enumerate} 
        \item[(a)] $ f$   is a function with domain  $ {}^\mu(\theta_0),0,$ 
    
        \item[(b)] $u \in [\mu]^\theta,$
        
        \item[(c)] $ \bar \rho = \langle \rho_i \colon i < \theta \rangle $ is  with no repetition,
        
        \item[(d)] $\rho_i \in {}^u(\theta_0)$,
    
        \item[(e)] $(\forall \nu \in {}^\mu(\theta_0)[\rho_i \subseteq \nu \Rightarrow f(\nu) = i],$  
    
        \item[(f)]  $(\forall \nu \in {}^\mu(\theta_0)[(\bigwedge_{i < \theta}(\rho_i \nsubseteq \nu)) \Rightarrow f(\nu) = 0].$  
    \end{enumerate}  

    We write $f = f^*_{u,\bar \rho}$, if $u,\bar \rho$ witness that $f \in {\mathscr  F}$ as above. Recalling $\mu = \mu^\theta$, clearly $|{\mathscr  F}| = \mu$.  Let ${\mathscr  F} = \{f_\varepsilon \colon \varepsilon < \mu\}$ and we let $\bar f = \langle
    f_\varepsilon \colon \varepsilon < \mu \rangle$.  Clearly clauses (a),(b) of
    Definition \ref{d.8} (with $\mu, \mu, \theta_0, \theta_1 = \theta$ here
    standing for $\chi, \mu, \theta, \theta_1$ there)  holds and let us check clause (c).  So suppose
    $\varrho \in {}^\mu \theta$ and let $R = R_\varrho \coloneqq  \{\nu \in {}^\mu(\theta_0)$: for every $\varepsilon < \mu$ we have $f_\varepsilon(\nu) \ne
    \varrho(\varepsilon)\}$.  We have to prove $|R| < \theta$  (as we have chosen $\Upsilon = \theta$). 
    
    Towards contradiction, assume that 
    $R \subseteq {}^\mu(\theta_0)$ has cardinality $\ge \theta$ and
    choose $R' \subseteq R$ of cardinality $\theta$.  Hence we can find $u \in
    [\mu]^\theta$ such that $\langle \nu \, {\rest} \, u \colon \nu \in R'
    \rangle$ is without repetitions.
    
    Let $\{\nu_i \colon i < \theta\}$ list $R'$ without repetitions  and let $\rho_i \coloneqq \nu_i \, {\rest} \, u$.
    Now let $\bar \rho = \langle \rho_i \colon i < \theta \rangle$, so  $f^*_{u,\bar \rho}$ is well defined and belong 
    to ${\mathscr  F}$.  Hence for some $\zeta < \mu$ we have $f^*_{u,\bar \rho} = f_\zeta$.  Now for each $i < \theta,\nu_i \in
    R' \subseteq R$ hence by the definition of $R,(\forall \varepsilon < \mu)(f_\varepsilon(\nu_i) \ne \varrho(\varepsilon))$ in particular for $\varepsilon = \zeta$ we get $f_\zeta(\nu_i) \ne \varrho(\zeta)$.  But by the choice of $\zeta,f_\zeta(\nu_i) = f^*_{u,\bar \rho}(\nu_i)$ and by the definition of $f^*_{u,\bar \rho}$ recalling $\nu_i \, {\rest} \, u = \rho_i$ we 
    have $f^*_{u,\bar \rho}(\nu_i) = i$, so $i = f_\zeta(\nu_i)  \ne \varrho(\zeta)$.  This holds for every $i < \theta$ whereas $\varrho \in {}^\mu \theta$, contradiction.

    \underline{Case 2}: $2^\theta \le \mu$ and $\mathbf U_\theta[\mu] = \mu$.
    
    Let $\Upsilon = (2^\theta)^+$ or just $(2^{< \theta})^+$ hence $\Upsilon = \cf(\Upsilon)
    \le \mu^+ \le 2^\mu$. Let $\{u_i \colon i < \mu\} \subseteq [\mu]^\theta$ exemplify $\mathbf U_\theta[\mu] = \mu$. Define ${\mathscr  F}$ as in case 1 except that for
    notational simplicity $\theta_0 = \theta$ and we restrict ourselves in the definition of ${\mathscr  F}$
    to $u \in \bigcup\{{\mathscr  P}(u_i) \colon i < \mu\}$.  Clearly $|{\mathscr  F}| 
    \le \mu$ and choose $\bar f = \langle f_\varepsilon \colon \varepsilon < \mu\rangle$ as there.
    
    Assume that $\varrho \in {}^\mu \theta$ and $R = R_\varrho  \subseteq {}^\mu \theta$ is defined as in Case 1, and toward
    contradiction assume that $|R| \ge (2^{<\theta})^+$, hence we can find $\rho^*,\langle (\alpha_\zeta,\rho_\zeta) \colon \zeta < \theta
    \rangle$ such that:

    \begin{enumerate}
        \item[$(\ast)$] 

        \begin{enumerate}
            \item[(a)] $\rho^*,\rho_\zeta \in R,$

            \item[(b)] $\rho_\zeta \, {\rest} \, \{\alpha_\xi \colon \xi < \zeta\}  = \rho^* \, {\rest} \, \{\alpha_\xi \colon \xi < \zeta\},$

            \item[(c)] $\rho_\zeta(\alpha_\zeta) \ne \rho^*(\alpha_\zeta).$
        \end{enumerate}
    \end{enumerate}

    [Why?  I had thought it should be clear but we elaborate. Let $\chi_* = (2^\mu)^+$ hence $\theta^M,R,\bar f \in {\mathscr  H}(\chi)$.  
    Now by induction on $\zeta < \theta$ we choose $M_\zeta$ such that:

    \begin{enumerate}
        \item[$(\ast)'$] 

        \begin{enumerate}
            \item[(a)]  $M_\zeta \prec ({\mathscr  H}(\chi),\in),$

            \item[(b)] $M_\zeta$ has cardinality $2^{|\zeta|} + \aleph_0,$

            \item[(c)] $\mu,\theta,\bar f,\varrho$ belongs to $M_\zeta,$

            \item[(d)] $[M_\zeta]^{|\zeta|} \subseteq M_\zeta,$

            \item[(e)] if $\varepsilon < \zeta$ then $M_\varepsilon  \prec M_\zeta$ and $M_\varepsilon \in M_\zeta.$
        \end{enumerate}
    \end{enumerate}

    Clearly possible and the cardinality of $M \coloneqq \cup\{M_\alpha \colon  \alpha < \theta\}$ is $\Sigma(2' + \aleph_0) = 2^{< \theta}$.
    
    Now choose $\rho^* \in R \backslash M$, possible because $\|M\| = 2^{< \theta} < (2^{< \theta})^+ \le |R|$.  Next, we try to choose a  pair $(\rho_\zeta, \alpha_\zeta)$ by induction on $\zeta < \theta$ such that:
    
    \begin{enumerate}  
        \item[$(\ast)$ (a)] $[M_\zeta]^{|\zeta|} \subseteq M_\zeta,$
        
        \item[(b)] $\rho_\zeta,\alpha_\zeta \in M_{\zeta +1},$
        
        \item[(c)] if $\varepsilon < \zeta$ then $\rho_\zeta(\alpha_\varepsilon) \ne \rho^*(\alpha_\varepsilon),$
        
        \item[(d)]   $\rho^*(\alpha_\zeta) \ne \rho_\zeta(\alpha_\zeta)$.
    \end{enumerate} 
    
    Why can we carry the induction on $\zeta$?  
    Arriving to $\zeta$ the set $X_\zeta 
    = \{(\rho_\varepsilon,\alpha_\varepsilon) \colon \varepsilon  < \zeta\}$ is a subset of $M_\zeta$ of cardinality $\le |\zeta|$ hence it belongs to $M_\zeta$.
    Consider the set ${\mathscr  X}_\zeta = 
    \{\rho \in R$: for every pair $(\rho,\alpha)$ from $\lambda,\rho(\alpha_\varepsilon) \ne \rho_\varepsilon(\alpha_\varepsilon)\}$.
    
    Now obviously this set belongs to $M_\zeta$, also $\rho^* \in {\mathscr  X}_\zeta$ by the 
    induction hypothesis on $\zeta$.  But $\rho^*$ being not in $M$ hence not  in $M_\zeta$; hence $|{\mathscr  X}_\zeta| > 1$  (and much more hence we can choose $\rho' \ne \rho'' \in {\mathscr  X}$ in $M_\zeta$.   As $\rho' \ne \rho''$ are members of ${\mathscr  X}_\zeta$ and 
    ${\mathscr  X}_\zeta \subseteq R \subseteq {}^M \theta$, there is $\alpha < \mu$ such that $\rho'(\alpha) \ne \rho''(\alpha)$ and as $\rho',\rho'' \in M_\zeta$, \wilog \,
    $\alpha \in M_\zeta$.  As $\rho'(\alpha) \ne \rho''(\alpha)$ necessarily for
    some $\rho \in \{\rho',\rho''\}$ we have $\rho^*(\alpha) \ne \rho(\alpha)$.  So
    we can choose $(\rho_\zeta,\alpha_\zeta) = (\rho,\alpha)$.
    
    Having carried out the induction clearly 
    $\rho^*, \langle (\rho_\zeta, \alpha_\zeta) \colon \zeta < \theta\rangle$ are as required.
    
    Clearly $\langle \alpha_\zeta \colon \zeta < \theta \rangle$ is with no repetitions.
    
    So by the choice of $\{u_i \colon i < \mu\}$ as exemplifying $\mathbf U_\theta[\mu] = \mu$, i.e., the definition of $\mathbf U_\theta[\mu]$,  for some $i < \mu$ the set $u_i \cap \{\alpha_\zeta \colon \zeta <
    \theta\}$ has cardinality $\theta$, call it $u$. So $\{\rho \, {\rest} \, u \colon \rho \in R\}$ has cardinality $\ge \theta$ and we can continue as in Case 1.  If we wish to use $\theta_0 \in [\theta,\mu]$ instead of $\theta_0 = \theta$ above let pr be a pairing function on $\mu$.  Without loss of generality each $u_i$ is closed under pr and its inverses, and $f^*_{u,\bar \rho}$ is defined iff for some $\varepsilon < \mu$ we have  $u \subseteq u_\varepsilon$ and $i < \theta \Rightarrow \text{ Rang}(\rho_i) \subseteq u_\varepsilon$.  In the end, choose $i$ such that $u_i \cap \{{\text{\rm pr\/}} (\alpha_\zeta,{\text{\rm pr\/}}
    (\rho_\zeta(\alpha_\zeta),\rho^*(\alpha_\zeta)) \colon \zeta < \theta\}$  has cardinality $\theta$.  
    
    \underline{Case 3}: $\mathbf U_J[\mu] = \mu$ where $J = [\sigma]^{< \theta}$
    for $\sigma$ such that $\sigma^\theta \le \mu, \sigma \ge \theta$.  
    
    Let $\Upsilon = [2^{< \sigma}]^+$.   The proof is similar to Case 2, only find now $\rho^*$ and $\langle (\alpha_\zeta, \rho_\zeta) \colon \zeta < \sigma \rangle$ and choose $i < \mu$
    such that $u_i \cap \{\text{pr}(\alpha_\zeta,\rho_\zeta(\alpha_\zeta)) \colon \zeta < \sigma\}$ has cardinality $\theta$.
    
    \underline{Case 4}: $\mu > \theta \ne \cf(\mu)$ and $\mu$ is strong
    limit.
    
    Follows by case 2. 
    
    \underline{Case 5}:  $\mu \ge \beth_\omega(\theta)$.
    
    By \cite{Sh:460} we can find a regular $\sigma < \beth_\omega(\theta)$ which is $> \theta$ and such that $\mathbf U_\sigma[\mu] = \mu$, so Case 2 applies.
    
    Note that in case 2, we can  replace $2^\theta \le \mu$ by: there is ${\mathscr  F} \subseteq {}^\theta \theta$ of cardinality $\le \mu$ such that for every $A \in [\theta]^\theta$ for some $f \in {\mathscr  F}$ we have $\theta = \text{ Rang}(f \, {\rest} \, A)$.  The proof is obvious.
    
    Similarly in case 3.   
\end{PROOF}

There are some obvious monotonicity properties regarding the \Bb-property.

\begin{claim} \label{d.11} 
    1) If $(\lambda,\bar C)$ has the
    $\theta$-{\rm \Bb}-property (see Definition \ref{d.5Y}(1))  and $\bar C = \langle C_\delta \colon \delta \in S
    \rangle$ and $C^1_\delta$ is an unbounded subset of $C_\delta$ for
    $\delta \in S$ \underline{then} $(\lambda,\langle C^1_\delta \colon \delta \in S \rangle)$ has a $\theta$-{\rm \Bb} property.  We can also increase $S$ and 
    decrease $\theta$; we can deal with the $(\Upsilon, \theta)$-\Bb-property. 
    
    2) Assume $\lambda = {\text {\rm cf\/}}(\lambda) > \aleph_0,S \subseteq
    \{\delta \colon \delta < \lambda,{\text{\rm cf\/}}(\delta) = \kappa\}$ is
    stationary and $\kappa < \chi \le \lambda$. 
    \underline{Then} some $\langle C_\delta \colon \delta \in S \rangle$ is a
    $(\lambda,\kappa,\chi)$-{\rm \Bb}-parameter. 
\end{claim} 

\demo{Proof}  1) Obvious. \newline
2) Easy (note: goodness is not required). 
\hfill$\square_{\ref{d.11}}$

We may consider the following generalization.

\begin{definition}\label{d.11A}
    1) We say $\bar C$ is a
    \emph{$\lambda$-\Bb$_3$-parameter} \underline{if}:
    
    \begin{enumerate} 
        \item[(a)]  $\lambda$ is regular uncountable, 
        
        \item[(b)]  $S$ is a stationary subset of $\lambda$ (consisting of limit ordinals), 
        
        \item[(c)]  $\bar C = \langle C_\delta, C'_\delta, C_{\delta, \varepsilon} \colon \delta \in S,\varepsilon < \cf(\delta) \rangle,$
        
        \item[(d)]  $C_\delta \subseteq \delta = \sup(C_\delta)$ and $C'_\delta \subseteq \delta = \sup(C'_\delta)$, $\otp(C'_\delta) = \kappa$ and $C_\delta = \cup\{C_{\delta,\varepsilon} \colon \varepsilon < \cf(\delta)\},$
        
        \item[(e)]  $\zeta < \cf(\delta) \Rightarrow \sup(\dbcu_{\varepsilon < \zeta} C_{\delta,\varepsilon}) < \delta$ and $C_{\delta,\varepsilon}$ increase with $\varepsilon.$ 
    \end{enumerate} 
    
    2)  We say that a $\lambda$-\Bb$_3$-parameter  $\bar C$ is\footnote{if $\langle C_\delta \colon \delta \in S \rangle$
    is a good $(\lambda,\kappa,\chi)$-\Bb-parameter $C'_\delta \subseteq
    \delta = \sup(C'_\delta)$, $\otp(C'_\delta) = \kappa$ and we let
    $C_{\delta, \varepsilon} = \{\alpha \in C_\delta \colon \otp(\alpha \cap
    C'_\delta) \le \varepsilon +1\}$ for $\varepsilon < \kappa$ then
    $\langle C_\delta, C_{\delta,\varepsilon} \colon \delta \in S, \varepsilon
    <\kappa \rangle$ is a $\lambda$-{\rm \Bb}$_1$-parameter.} \emph{good} if $\alpha <
    \lambda \Rightarrow \lambda > |\{C_{\delta,\zeta} \colon 
    \delta \in S,\zeta < \cf(\delta)$ and the $\zeta$-th member of  $C'_\delta$ is $\le \alpha\}|$. 
    
    3) The $(D,\Upsilon,\theta)$-\Bb$_3$-property for  $\bar C$ means the $(D,\Upsilon,\theta)$-\Bb-property for $\langle C_\delta \colon \delta \in S \rangle$, similarly for the other variants of
    this property (so the $(D,\Upsilon,\theta)$-\Bb$_3$-property for $\lambda$ for a stationary $S \subseteq S^\lambda_\kappa$ is the same as the $(D, \Upsilon, \theta)$-\Bb-property). 
    
    4) If $\{C_{\delta,\varepsilon} \colon \varepsilon < \cf(\delta)\} 
    = \{C_\delta \cap \alpha \colon \alpha \in \mathrm{nacc}(C_\delta)\}$ we may 
    write just $\langle C_\delta \colon \delta \in S \rangle$ (so \Bb$_3$ is  like \Bb in this case).
\end{definition}

We may consider the following variant.

\begin{claim}\label{d.11B} 
    If (A) then  (B), where:

    \begin{enumerate}
        \item[(A)]

        \begin{enumerate}
            \item[(a)]  $\lambda$ is regular uncountable, $S \subseteq \{\delta < \lambda \colon {\text{\rm cf\/}}(\delta) = \kappa\}$ is stationary and $\chi < \lambda,$

            \item[(b)] ${\mathfrak B}$ is a model with universe $\lambda$ such that $A \in [\lambda]^{< \chi} \Rightarrow c \ell_{\mathfrak B}(A) \in [\lambda]^{< \chi}$; we may allow infinitary functions provided that:

            \begin{enumerate}
                \item[$(*)_{{\mathscr  B},\kappa}$]  if $\langle A_i \colon i < \kappa \rangle$ is increasing, $A_i \in [\lambda]^{< \chi}$ then $c \ell_{\mathfrak B}(\dbcu_{i < \kappa} A_i) = \dbcu_{i < \kappa} c \ell_{\mathfrak B}(A_i).$
            \end{enumerate}

            \item[(c)] $\bar C = \langle C_\delta \colon \delta \in S \rangle$ is a good $(\lambda,\kappa,\kappa^+)$-{\rm \Bb}-parameter, 

            \item[(d)] $S^* = \{\delta \in S$: for every $\alpha < \delta, c \ell_{\mathfrak B}(\alpha) \subseteq \delta\},$

            \item[(e)]    $C^*_{\delta,\varepsilon} = \alpha_{\delta,\varepsilon} \cap c \ell_{\mathfrak B}(C_\delta \cap \alpha_{\delta,\varepsilon})$ for $\delta \in S^*$ and  $\varepsilon < \kappa$ where $\alpha_{\delta,\varepsilon} \in C_\delta$ is increasing with $\varepsilon$ with limit $\delta$ and $C^*_\delta = \cup\{C^*_{\delta,\varepsilon}:\varepsilon < { \text{\rm cf\/}} (\delta)\} $ and $C'_\delta = \{\alpha_{\delta,\varepsilon} \colon \varepsilon < { \text{\rm cf\/}}(\delta)\}$.
        \end{enumerate}

        \item[(B)] $\langle C^*_\delta, C_\delta, C^*_{\delta,\varepsilon} \colon \delta \in S^*,\varepsilon < { \text{\rm cf\/}}(\delta) \rangle$ is a good $(S,\kappa,\chi)$-{\rm \Bb}$_3$-parameter (and $S \backslash S^*$ is non-stationary); if $\chi \le \kappa$ then $\dbcu_{\varepsilon < \kappa} C^*_{\delta,\varepsilon} = c \ell_{\mathfrak B}(C_\delta)$.
    \end{enumerate}
\end{claim} 

\begin{PROOF}{\ref{d.11B}}
    Easy.
\end{PROOF}

\begin{claim}\label{d.11C}  
    1) Assume $\bar C = \langle C_\delta,C'_\delta, C_{\delta,\varepsilon} \colon \delta \in S, \varepsilon < 
    { \text{\rm cf\/}}(\delta) \rangle$ is a good $(\lambda, \kappa,\chi)$-{\rm \Bb}$_3$-parameter.
    \underline{Then} Claim \ref{d.6} holds (i.e. we get there the $(D,\theta)$-{\rm \Bb}$_3$-property).

    2) In part (1) hence in Claim \ref{d.6}, we can weaken clause (e) to:

    \begin{enumerate}
        \item[(e)$'$]  (if $\lambda = 2^\mu$ we can use $h =  { \text{\rm id\/}}_\lambda$): there is $h \colon \lambda \rightarrow 2^\mu$ increasing continuous with limit $2^\mu$ and there is a set  $T = \{\eta_\alpha \colon \alpha < 2^\mu\}$ and function $g \colon \lambda \rightarrow \lambda$ such that if $\eta \in {}^\lambda(2^\mu)$ then for stationary many $\delta \in S$ we have:

        \begin{enumerate}
            \item[(a)] $\eta \, {\rest} \, C_\delta \in \{\eta_\alpha \colon \alpha < h(g(\delta))\},$ 

            \item[(b)] $\zeta < { \text{\rm cf\/}} (\delta) \Rightarrow \eta \, {\rest} \, (\dbcu_{\varepsilon < \zeta} C_{\delta,\varepsilon}) \in \{\eta_\beta \colon \beta < \delta\}$.
        \end{enumerate}
    \end{enumerate}
\end{claim} 

\begin{PROOF}{\ref{d.11C}}
    Same proof as Claim \ref{d.6}.
    
    We give details, we are assuming:
    
    \begin{enumerate} 
        \item[$\circledast$]  
        
        \begin{enumerate} 
            \item[(a)]  $\lambda =  {\text {\rm cf\/}}(2^\mu),D$ is a $\mu^+$-complete filter on $\lambda$ extending the club filter,
            
            \item[(b)]  $\kappa  = {\text {\rm cf\/}}(\kappa) < \chi \le \lambda,$
            
            \item[(c)]  $\bar C = \langle C_\delta, C'_\delta, C_{\delta,\varepsilon} \colon \delta \in S \rangle$ is a good $(\lambda,\kappa,\chi)-\text{\rm \Bb}_3$-parameter, 
            
            \item[(d)]  $ 2^{< \chi} \le 2^\mu$ and $\theta \le \mu,$
            
            \item[(e)]  $ \alpha < 2^\mu \Rightarrow 
            |\alpha|^{<\kappa>_{\text{\rm tr\/}}} < 2^\mu$ or,
            
            \item[(f)] clause $(e)'$ above,
            
            \item[(g)] $ \Sep(\mu,\theta)$.
        \end{enumerate} 
    \end{enumerate} 
    
    So we have to deduce
    
    \begin{enumerate} 
        \item[$\boxdot$]   $\bar C^* = \langle C_\delta \colon \delta \in S \rangle$
        has the $(D,2^\mu,\theta)$-{\rm \Bb}-property recall that this means that the number of colours is  $\theta$ not just 2).
    \end{enumerate} 

    Why $\boxdot$ holds?  Let $\mathbf F$ be a  given $(\bar C^*, 2^\mu, \theta)$-colouring.
    
    By assumption (f) we have $\Sep(\mu,\theta)$ which means (see Definition \ref{d.8}(2)) that for some $\Upsilon = \cf(\Upsilon) \le
    2^\mu$ we have $\Sep(\mu,\mu,\theta,\theta,\Upsilon)$.

    Let ${\mathscr  F} = \{f_\xi \colon \xi < \mu\}$ where ${\mathscr  F}$
    exemplify $\Sep(\mu,\mu,\theta,\theta,\Upsilon)$, see Definition \ref{d.8}(1) and for $\rho \in {}^\mu \theta$ let  $\Sol_\rho \coloneqq \{\nu \in {}^\mu \theta \colon \text{ for every } \varepsilon 
    < \mu$ we have $\rho(\varepsilon) \ne f_\varepsilon(\nu)\}$ where  $\Sol$ stands for solutions, so by clause (c) of the Definition \ref{d.8}(1) of $\Sep$ it follows that:
    
    \begin{enumerate} 
        \item[$(*)$]  $\rho \in {}^\mu \theta 
        \Rightarrow |\Sol_\rho| < \Upsilon$.
    \end{enumerate} 
    
    We now assume that sub-clause $(\alpha)$ of clause (e) of $\circledast$ holds and we shall find $h,g,\langle \eta_\alpha \colon \alpha < 2^\mu\rangle$ as required in sub-clause $(e)(\beta)$ of $\circledast$. Let $h$ be an increasing continuous function from $\lambda$ into $2^\mu$ with unbounded range; if $2^\mu$ is regular, $h = \id_\lambda$ is O.K.  Recall that $\lambda = \cf(2^\mu) \ge 2^{< \chi}$. Let $\cd$ be a one-to-one function from ${}^\mu (2^\mu)$ onto $2^\mu$ such that: $$ \alpha = \cd(\langle \alpha_\varepsilon \colon \varepsilon < \mu \rangle) \Rightarrow \alpha \ge \sup\{\alpha_\varepsilon \colon \varepsilon < \mu\}, $$ and let the function $ \cd_i \colon2^\mu \rightarrow 2^\mu$ for $i < \mu$ be  such that $\cd_i(\cd(\langle \alpha_\varepsilon \colon  \varepsilon < \mu \rangle)) = \alpha_i$.
    
    Let ${\mathscr  C} = \{C_{\delta,\zeta} \cap \alpha \colon\delta \in S,\zeta <
    \kappa \text{ and } \alpha$ is the $\zeta$-th member of  $C'_\delta\}$, recall that trivially $|{\mathscr  C}| \le \lambda \le 2^\mu$ and let $T= \{\eta \colon$ for some $C \in {\mathscr  C}$ (hence $C$ is of cardinality $< \chi$), we have $\eta \in {}^C (2^\mu)\}$.  By assumptions (c) + (d) clearly $|T|  = \dsize \sum_{C \in {\mathscr  C}} 2^{\mu \cdot |C|} \le |{\mathscr  C}| \times (2^\mu)^{< \chi} \le 2^\mu$ but $|T| \ge 2^\mu$ hence $|T| =   2^\mu$.  So let us list $T$ possibly with repetitions as
    $\{\eta_\alpha \colon\alpha < 2^\mu\}$, \wilog \, such that  $[\alpha < h(\beta) \Rightarrow
    \text{ Rang}(\eta_\alpha) \subseteq \beta]$ 
    and let $T_{< \alpha} = \{ \eta_\beta:\beta < h(\alpha)\}$ for $\alpha
    < \lambda$. For $\delta \in S$ let ${\mathscr  P}_\delta \coloneqq \{\eta \colon\eta$ is a function from $C_\delta$ into $h(\delta)$ such that for 
    every $\zeta < \kappa$ we have $\eta \, {\rest} \, C_{\delta,zeta}  \in T_{< \delta}\}$, clearly $\eta \in {\mathscr  P}_\delta  \Rightarrow \eta \in {}^{(C_\delta)} h(\delta)$.  Clearly
    $\{\eta \, {\rest} \, C_{\delta,\zeta} \colon \eta \in {\mathscr  P}_\delta,\zeta <  \kappa\}$ is naturally a tree with ${\mathscr  P}_\delta$ the set of $\otp(C_\delta)$-branches and cf$(\otp(C_\delta)) = \kappa$ and set of nodes $\subseteq T_{< \delta}$ hence with $< 2^\mu$ nodes hence by assumption (e),

    \begin{enumerate}
        \item[$(\ast \ast)$] $|{\mathscr  P}_\delta| < 2^\mu$.
    \end{enumerate}

    If $\lambda = 2^\mu$ then there is a function $g \colon \lambda \rightarrow \lambda$ as required in sub-clause $(e)(\beta)$ of $\circledast$, just
    let $g(\alpha) = \text{\rm Min}\{\gamma < \lambda \colon \gamma > \alpha$ and for every $\delta \in S \cap (\gamma +1)$ we have ${\mathscr  P}_\delta
    \subseteq \{\eta_\beta \colon \beta < h(\gamma)\}$.  If $\lambda = \text{\rm cf}(2^\mu) < 2^\mu$ we should choose $\langle \eta_\alpha \colon \alpha <
    2^\mu \rangle$, by letting $\eta'_{2\alpha} = \eta_\alpha$ and making, e.g. $\delta \in S \Rightarrow {\mathscr  P}_\delta \subseteq \{\eta'_{2
    \alpha+1} \colon \alpha < |{\mathscr  P}_\delta|^+ (< 2^\mu)\}$. In any case, 

    \begin{enumerate}
        \item[$(\ast \ast \ast)$] $h, g, \langle \eta_\alpha \colon \alpha < 2^\lambda\rangle$ are as in $(c)(\beta)$ of $\circledast$ and let ${\mathscr  P}_\delta = \{\eta_\alpha \colon \alpha < h(g(\alpha))\}$.
    \end{enumerate}

    For each $\eta \in {\mathscr  P}_\delta$
    and $\varepsilon < \mu$ we define $\nu_{\eta,\varepsilon}
    \in {}^{C_\delta}(h(\delta))$ by $\nu_{\eta,\varepsilon} (\alpha) = 
    \cd_\varepsilon (\eta(\alpha))$ for $\alpha \in C_\delta$.  Now for $\eta \in {\mathscr  P}_\delta$, clearly $\rho_\eta \coloneqq   \langle \mathbf F (\nu_{\eta,\varepsilon}) \colon \varepsilon < \mu \rangle$  belongs to ${}^\mu \theta$.
    Clearly $\{\rho_\eta \colon \eta \in {\mathscr  P}_\delta\}$ is a subset of ${}^\mu \theta$ 
    of cardinality $\le |{\mathscr  P}_\delta|$ which as said above is  $< 2^\mu$.  Recall the definition of $\Sol_\rho$  from the beginning of the proof and remember $(*)$ above, so $\eta \in {\mathscr  P}_\delta \Rightarrow \Sol_{\rho_\eta}$ has cardinality $< \Upsilon = \cf(\Upsilon) \le 2^\mu$.  Hence we can find $\varrho^*_\delta \in {}^\mu \theta \setminus \cup \{\Sol_{\rho_\eta} \colon \eta \in {\mathscr  P}_\delta\}$.  

    [Why?  if $\Upsilon < 2^\mu$ then $|\cup\{\Sol_{\rho_\eta} \colon \eta
    \in {\mathscr  P}_\delta\}| \le \Sigma\{|\Sol_{\rho_\eta}| \colon \eta \in
    {\mathscr  P}_\delta\} \le \Upsilon \times |{\mathscr  P}_\delta| < 2^\mu =
    \theta^\mu = |{}^\mu \theta|$ and if $\Upsilon = 2^\mu$ then $2^\mu$
    is regular, $|{\mathscr  P}_\delta| < 2^\mu,\eta \in {\mathscr  P}_\delta \Rightarrow |\Sol_{\rho_\eta}| < 2^\mu$ and again $|\cup\{\Sol_{\rho_\eta} \colon \eta \in {\mathscr  P}_\delta\}| < 2^\mu =
    |{}^\mu \theta|$.] 
    
    Let $\varepsilon < \mu$.  Recall
    that $\varrho^*_\delta \in {}^\mu \theta$ for $\delta \in S$ and $f_\varepsilon$ from 
    the list of ${\mathscr  F}$ is a function from ${}^\mu \theta$ to $\theta$ so $f_\varepsilon(\varrho^*_\delta) < \theta$.  Hence
    we can consider the sequence $\bar c^\varepsilon = \langle f_\varepsilon(\varrho^*_\delta) \colon \delta 
    \in S \rangle \in {}^S \theta$ as a candidate for being an {\bf F}-\Bb-sequence.  If one of them is, we are done. So assume towards a contradiction that for each $\varepsilon< \mu$ there is a sequence
    $\eta_\varepsilon \in {}^\lambda (2^\mu)$ that exemplifies the failure of $\bar c^\varepsilon$ to be an $\mathbf F$-\Bb-sequence,  so there is a $E_\varepsilon \in D$ such that:

    \begin{enumerate}
        \item[$\boxtimes_1$] $\delta \in S \cap E_\varepsilon \Rightarrow \mathbf F (\eta_\varepsilon \, {\rest} \, C_\delta) \ne f_\varepsilon(\varrho^*_\delta)$.
    \end{enumerate}

    Define $\eta^* \in {}^\lambda(2^\mu)$ by $\eta^*(\alpha) =  \cd(\langle \eta_\varepsilon(\alpha) \colon \varepsilon  < \mu \rangle)$.  Now as $\lambda$ is regular uncountable it follows
    by choice of $h$ that $E \coloneqq \{\delta < \lambda:$ for every $\alpha < \delta$ we have $\eta^* (\alpha) < h(\delta)$ and if $\delta'
    \in S,\zeta < \kappa,\alpha$ the $\zeta$-th member of 
    $C'_{\delta'} \cap \delta$ then $\eta^* \, {\rest} \, C_{\delta',\zeta}\}$ is  a club of $\lambda$ (see the choice of $T$ and recall as well 
    that by assumption (c) the sequence 
    $\langle (C_\delta,C'_\delta,C_{\delta,\zeta}):\delta \in S,\zeta < \kappa\rangle$ is good, see Definition \ref{d.11A}(2).
    
    By clause (a) in the assumption of our  Claim \ref{d.6}, the filter $D$ includes the clubs of $\lambda$ so clearly $E \in D$; also $D$ is $\mu^+$-complete hence $E^* \coloneqq \cap \{E_\varepsilon \colon \varepsilon < \mu\} \cap E$ belongs to $D$.  Now choose $\delta \in E^* \cap S$, fixed for the rest of the proof;  clearly $\eta^*  \, {\rest} \, C_\delta \in {\mathscr  P}_\delta$; just check the definitions of ${\mathscr  P}_\delta$ and
    $E,E^*$.  Let $\varepsilon < \mu$.
    Now recall that $\nu_{\eta^* \, {\rest} \, C_\delta,\varepsilon}$ is the function from $C_\delta$ to $h(\delta)$ defined by $$
    \nu_{\eta^* \, {\rest} \, C_\delta,\varepsilon}(\alpha) = \cd_\varepsilon (\eta^*(\alpha)).$$
    
    But  by our choice of 
    $\eta^*$ clearly $\alpha \in C_\delta \Rightarrow \cd_\varepsilon (\eta^*(\alpha)) = \eta_\varepsilon (\alpha)$, so $$
    \alpha \in  C_\delta \Rightarrow \nu_{\eta^* \, {\rest} \, C_\delta,\varepsilon} (\alpha) = \eta_\varepsilon (\alpha) \quad \text{ so } \quad
    \nu_{\eta^* \, {\rest} \, C_\delta,\varepsilon} = \eta_\varepsilon \, {\rest} \, C_\delta. $$
    
    Hence $\mathbf F (\nu_{\eta^* \, {\rest} \, C_\delta,\varepsilon}) = \mathbf F (\eta_\varepsilon \, {\rest} \, C_\delta)$.  As
    $\delta \in E^* \subseteq E_\varepsilon$ clearly $\mathbf F (\eta_\varepsilon \, {\rest} \, C_\delta) \ne f_\varepsilon (\varrho^*_\delta)$ and as $\eta^* \, {\rest} \, C_\delta \in {\mathscr  P}_\delta$ clearly $\rho_{\eta^* \, {\rest} \, C_\delta} \in {}^\mu
    \theta$ is well defined.  Now easily $\rho_{\eta^* \, {\rest} \, C_\delta}
    (\varepsilon) = \mathbf F(\nu_{\eta^* \, {\rest} \, C_\delta,\varepsilon})$ by 
    the definition of $\rho_{\eta^* \, {\rest} \, C_\delta}$, so we have $\rho_{\eta^* \, {\rest} \, C_\delta}(\varepsilon) \ne
    f_\varepsilon(\varrho^*_\delta)$.
    
    As this holds for every $\varepsilon < \mu$ it follows  by the definition of $\Sol_{\rho_{\eta^* \, {\rest} \, C_\delta}}$ that $\varrho^*_\delta \in \Sol_{\rho_{\eta^* \, {\rest} \, C_\delta}}$. But $\eta^* \, {\rest} \, C_\delta \in {\mathscr  P}_\delta$ was proved
    above so $\rho_{\eta^* \, {\rest} \, C_\delta} \in \{\rho_\eta \colon \eta \in {\mathscr 
    P}_\delta\}$ hence $\varrho^*_\delta \in \Sol_{\rho_{\eta^* \, {\rest} \, C_\delta}} \subseteq \cup \Sol_{\rho_\nu} \colon \nu \in
    {\mathscr  P}_\delta\}$
    whereas $\varrho^*_\delta$ has been chosen
    outside this set, a contradiction.  \hfill$\square_{\ref{d.11C}}$
\end{PROOF}

\begin{claim}\label{d.20} 
    In the proof of Claim \ref{d.6} instead of
    $f \colon \lambda \rightarrow \Upsilon$ and 
    $f_\delta \colon C_\delta \rightarrow \theta$ we can use $f \colon [\lambda]^{< \aleph_0} \rightarrow \Upsilon$ and
    $f_\delta \colon [C_\delta]^{< \aleph_0} \rightarrow \theta$ (so get the
    $\theta$-{\rm MD}$_\ell$-property for $\ell = 0,1,2$).
\end{claim}

\begin{PROOF}{\ref{d.20}}
    Let $\cd_*$  be a one-to-one function from ${}^{\omega >}\lambda$
    into $\lambda$ and without loss of generality $\alpha_\ell \le
    \text{\rm cd}_*(\langle \alpha_0,
    \alpha_1,\dotsc,\alpha_{n-1} \rangle) < { \text{\rm Max\/}}(\{\omega\}
    \cup \{\omega \alpha_{\pi(0)} + \ldots \omega \alpha_{\pi(n-1)} \colon \pi$ a
    permutation of $n\})$.  We are given $\bar C
    = \langle C_\delta \colon \delta \in S \rangle$, \wilog \, every $\delta \in
    S$ is infinite and divisible by $(\delta)$.
    
    Let $C_{1,\delta} = \{{\text{\rm cd\/}}_*(\eta) \colon \eta \in {}^{\omega >}
    (C_\delta)\}$ and $C'_{1,\delta} = { \text{\rm nacc\/}}(C_\delta)$. \hfill$\square_{\ref{d.20}}$ 
\end{PROOF}

\begin{definition}\label{d.12G}  
    1) For a set $A$ of ordinals and an ordinal $\alpha$ let $$A * \alpha \coloneqq \{\beta \colon \beta < \alpha \text{ or for some } \gamma \in A \text{ we have } \beta \in [\gamma,\gamma + \alpha)\}.$$
    
    2) The pair $(\lambda,\bar C)$ or just $\bar C$  has the \emph{$(\chi,\theta)$-MD$_0$-property} (or use MD instead of MD$_0$) \underline{if} we can find $\langle f_\delta \colon \delta \in S \rangle$ such that:
    
    \begin{enumerate}
        \item[(a)]  $f_\delta$ is a function from $C'_\delta = C_\delta * \alpha^*_\delta$ to $\theta$ for some $\alpha^*_\delta < \chi,$
    
        \item[(b)]  if $f$ is a function from $\lambda$ to $\theta$ and $\alpha < \chi$ is non-zero then  for stationary many $\delta \in S$ we have  $f_\delta \subseteq f$ and $\alpha^*_\delta = \alpha$.
    \end{enumerate}
    
    If we omit $\chi$ we mean $\chi = 2$ so $\alpha^*_\delta = 1$ for every $\delta \in S$.
    
    3) \emph{$(D,\chi,\theta)$-MD$_2$-property} means: if  $\tau$ is a relational vocabulary of cardinality $< \theta$,  \underline{then} we can find $\langle M_\delta \colon \delta \in S \rangle$ such that:
    
    \begin{enumerate}
        \item[(a)]  $M_\delta$ is a $\tau$-model with universe $C'_\delta = C_\delta * \alpha^*_\delta$ for some $\alpha^*_\delta < \chi,$
    
        \item[(b)]  if $M$ is a $\tau$-model with universe $\lambda$ and $\alpha < \chi$ \underline{then} the set  $A_{M,\alpha} = \{ \delta \in S \colon M \, {\rest} \, C_\delta = 
        M_\delta$ and $\alpha^*_\delta = \alpha\}$ is stationary and even belongs to $D^+$.
    \end{enumerate}
    
    If we omit $\chi$ we mean $\chi = 2$ so $\alpha^*_\delta =1$ for every $\delta \in S$.  Similarly with $D$. 
    
    4) We say $(\lambda,S)$ has the \emph{$(\chi,\theta)$-MD$_\ell$-property} \underline{if} for some $\bar C = \langle C_\delta \colon \delta \in S \rangle$ the pair $(\lambda,\bar C)$ has the $(\chi,\theta)$-MD$_\ell$-property.  We say that $\lambda$ has the $(\kappa, \chi, \theta)$-property if $(\lambda, S^\lambda_\kappa)$ has the $(\chi,\theta)$-MD$_\ell$-property.
\end{definition} 

The following claim shows the close affinity of the properties stated
in Theorem \ref{0.1} in the introduction and the
$\theta$-\Bb-property with which we deal in this section here. This claim is needed to deduce Theorem \ref{0.1} from Claim \ref{d.6}.

\begin{claim}\label{d.12} 
    1) Assume $(\lambda,\bar C)$ has the
    $\theta$-{\rm \Bb}-property where $\bar C = 
    \langle C_\delta \colon \delta \in S \rangle$ is a $(\lambda,\kappa,\chi)$-{\rm \Bb}-parameter.  If $\kappa < \chi,\theta = 
    \theta^{< \chi}$ \underline{then} $(\lambda,\bar C)$ has the
    $(\chi,\theta)$-{\rm MD}$_0$-property. 
    
    2) Assume in addition that $\bar C$ is a good$^+ \, (\lambda,\kappa,\chi)$-{\rm \Bb}-parameter \underline{then} $(\lambda,\bar C)$ has the $(D,\chi,\theta)$-{\rm MD}$_2$-property. 
    
    3) In part (2) we can allow $\tau$ to have function symbols but have:

    \begin{enumerate}
        \item[(a)]   $M_\delta$ is a $\tau$-model with universe $\supseteq C'_\delta = C_\delta * \alpha^*_\delta$ for some 
        $\alpha^*_\delta < \chi,$

        \item[(b)] if $M$ is a $\tau$-model with universe $\lambda$ and $\alpha < \lambda$ \underline{then} the following set belongs to $D^+$: $$ A'_{M,\alpha} = \{\delta \in S \colon \alpha^*_\delta = \alpha \text{ and } (M \, {\rest} \, c \ell_M(C'_\delta),c)_{c \in C'_\delta} \cong (M_\delta,c)_{c \in C'_\delta}\}.$$
    \end{enumerate}

    4) In part (2) we can replace good$^+$ by good \underline{if} we add: 
    
    \begin{enumerate} 
        \item[$\boxdot$]  each $C_\delta$ is closed under pr, where {\rm pr} is a
        fixed pairing function on $\lambda$.
    \end{enumerate} 
\end{claim} 

\begin{PROOF}{\ref{d.12}}
    The proof consists of a translation of one version of \midia to another (note: ``$\tau$ is a relational" is a strong restriction). 
    
    1) First fix $\alpha^* < \chi$ and deal with $\langle C''_\delta \colon \delta \in S \rangle,C''_\delta = C'_\delta \backslash
    [0,\alpha^*)$ recalling $C'_\delta = C_\delta * \alpha^*$.  
    Let $\cd$ be a one-to-one function from ${}^{\chi >}\theta$ onto 
    $\theta$.  Define $\langle F_\delta:\delta \in S \rangle$ as follows: if $f \in {}^{(C_\delta)}\theta$ then
    $F_\delta(f) = \cd(\langle
    f(\beta_{\delta,\varepsilon}) \colon \varepsilon < \otp(C_\delta)\rangle)$ 
    where $\langle \beta_{\delta,\varepsilon} \colon \varepsilon <
    \otp(C_\delta) \rangle$ list $C_\delta$ in
    increasing order so $F_\delta$ is a function 
    from ${}^{(C_\delta)} \theta$ into $\theta$.  As we assume
    ``$(\lambda,\bar C)$ has the $\theta$-\Bb-property", we can apply its
    definition, hence for 
    $\bar F = \langle F_\delta:\delta \in S \rangle$ there is a
    diamond sequence $\langle c_\delta:\delta \in S \rangle$.  
    
    We now define $f^*_\delta \colon C''_\delta \rightarrow \theta$ as follows:
    if $\beta_0$ is a member of $C_\delta$ and
    $\alpha \in [\beta_0,\beta_0 + \alpha^*)$ and $\alpha < \min(C''_\delta \backslash (\beta_0 +1))$ and $c_\delta = 
    { \text{\rm cd\/}}(\langle \xi_{\delta,\varepsilon} \colon \varepsilon 
    < { \text{\rm otp\/}}(C_\delta) \rangle$ 
    and $\xi_{\delta,\varepsilon} = { \text{\rm cd\/}} (\langle \zeta_{\delta,\varepsilon,i} \colon i < \alpha^* \rangle$ then
    $f^*_\delta(\alpha) = \zeta_{\delta,\otp(C_\delta \cap \beta_0),
    \alpha - \beta_0}$; otherwise $f^*_\delta(\alpha) = 0$.  So
    $f^*_\delta \colon C''_\delta \rightarrow \theta$. Why is $\langle f^*_\delta \colon \delta \in S \rangle$ as required for $\langle C''_\delta \colon \delta \in S \rangle$ (and $\alpha^*$)?
    Assume that $f \colon\lambda \rightarrow \theta$, so we can define
    $f' \colon \lambda \rightarrow \theta$ by $f'(\alpha) = { \text{\rm cd\/}}
    (\langle f(\alpha + \varepsilon) \colon \varepsilon < \alpha^* \rangle)$, hence for
    stationarily many $\delta \in S$ we have $c_\delta = F_\delta(f'
    \, {\rest} \, C_\delta)$.  Now we check that $f \, {\rest} \, C''_\delta$
    is equal to $f^*_\delta$ defined above.
    
    Let $\xi_{\delta,\varepsilon,i} = f'(\beta_{\delta,\varepsilon})$ and
    $\zeta_{\delta,\varepsilon +i} = f(\beta_{\delta,\varepsilon} +i)$ for
    $i < \alpha^*, \varepsilon < \otp(C_\delta)$ so $\xi_{\delta, \varepsilon} = f'(\beta_{\delta, \varepsilon}) = \cd(\langle \zeta_{\delta,\varepsilon,i} \colon i < \alpha^* \rangle), c_\delta = F_\delta(f' \, {\rest} \, C_\delta) 
    = \cd(\langle \xi_{\delta,\varepsilon} \colon \varepsilon 
    < \otp(C_\delta) \rangle$.  Hence by the definition of $f^*_\delta$ we
    have: if $\beta = \beta_{\delta,\varepsilon} \in C_\delta,\alpha \in
    [\beta,\beta + \varepsilon)$ and $\alpha < \min(C_\delta
    \backslash (\beta +1))$ then $f^*_\delta(\alpha) =
    \zeta_{\delta,\varepsilon,\alpha - \beta} =
    f(\beta_{\delta,\varepsilon} + (\alpha - \beta)) = f(\beta + (\alpha -
    \beta)) = \beta$ as required.
    
    To have ``every $\alpha^* < \chi$" (and guess also $f \, {\rest} \, (C'_\delta \backslash C''_\delta) \subseteq f \, {\rest} \, [0,\alpha^*))$
    we just have to partition $S$ to $\langle S_\alpha \colon \alpha < \chi \rangle$ such that each $\bar C \, {\rest} \, S_\alpha$ has the $\theta$-\Bb-property which we can do by \ref{d.12B}(3) below (or work a little more above, i.e., $f'(\alpha) = \cd(\langle f(\alpha + \varepsilon) \colon \varepsilon < \alpha_\delta \rangle \char 94 \langle f(\varepsilon) \colon \varepsilon < \alpha_\delta \rangle)$). 
    
    2), 3), 4)  Left to the reader.   
\end{PROOF}

\begin{claim}\label{d11Z}
    1) Assume  $\lambda = { \text{\rm cf\/}}(\lambda) > \kappa =  { \text{\rm cf\/}}(\kappa),S \subseteq S^\lambda_\kappa$ and

    \begin{enumerate}
        \item[(a)]  $\bar C^\ell$ is a $\lambda$-{\rm \Bb}-parameter for $\ell = 1,2,$

        \item[(b)] $\bar C^1 = \langle C^1_\delta \colon \delta \in S \rangle$ has the $(\Upsilon,\theta)$-{\rm \Bb} property, 

        \item[(c)] $A_\alpha \in [\alpha]^{< \chi} \, \wedge \,  \Upsilon^{|A_\alpha|} = \Upsilon$ for $\alpha < \lambda$ and $\delta \in S \Rightarrow C^2_\delta = \cup\{A_\alpha \colon \alpha \in C^1_\delta\}$,
    \end{enumerate}

    \underline{then} $\bar C^2$ has the $(\Upsilon,\theta)$-{\rm \Bb}-property. 

    2) If $\bar C$ is a $(\lambda,\kappa,\chi^+)$-{\rm \Bb}-parameter with the 
    $(\Upsilon,\theta)$-{\rm \Bb}-property, $\Upsilon^\chi =  \Upsilon$ and $C'_\delta =
    \{\alpha + I \colon \alpha \in C_\delta \vee \alpha = 0$ and $i < \chi\}$,
    \underline{then} $\bar C' = \langle C'_\delta \colon \delta \in {\text{\rm Dom\/}} (\bar C)\rangle$ is a $(\lambda,\kappa,\chi^+)$-{\rm \Bb}-parameter with  at he $\theta$-{\rm \Bb}-property.
\end{claim} 

\begin{PROOF}{\ref{d11Z}}
    Easy by coding (as in the proof of Claim \ref{d.12}).
\end{PROOF}

\begin{definition}\label{d.12A} 
    Assume  $\lambda > \kappa$ are regular and $\bar C = \langle C_\delta \colon \delta \in S \rangle$ is a $(\lambda, \kappa, \chi)$-\Bb-parameter with the $(\Upsilon,\theta)$-\Bb-property.  Let ID$_{\lambda, \kappa,\chi, \Upsilon, \theta}(\bar C) = \{A \subseteq  \lambda \colon \bar C \, {\rest} \, (S \cap A)$ does not have the $(\Upsilon, \theta)$-\Bb-property$\}$.  If $\Upsilon = \theta$  we may omit it.
\end{definition}

\begin{claim}\label{d.12B} 
    Let $\bar C,\lambda,\kappa,\chi,\Upsilon,\theta$
    be as above in Definition \ref{d.12A} and $J =  { \text{\rm ID\/}}_{\lambda,\kappa,\chi,\Upsilon,\theta}(\bar C)$. 
    
    1) $J$ is an ideal on $\lambda$ such that $\lambda \backslash S \in J,\lambda \notin J$ if $\Upsilon \ge \aleph_0$. 
    
    2) $J$ is $\lambda$-complete and even normal if $\Upsilon \ge \lambda$ (if $\aleph_0 \le \Upsilon < \lambda$ it is still $\Upsilon^+$-complete). 
    
    3) If $S \subseteq \lambda$ and $S \ne \emptyset$ {\rm mod} $J$, \underline{then}
    for any $\bar F^i = \langle F^i_\delta \colon \delta \in S \rangle,F^i_\delta \colon {}^{(C_\delta)} \theta \rightarrow \theta$,  $\lambda$ can be partitioned to $\theta$ subsets each having a 
    $\bar F^i$-{\rm \Bb}-sequence.
\end{claim} 

\begin{PROOF}{\ref{d.12B}}
    1) Let $\cd \colon \Upsilon \times \Upsilon \rightarrow
    \Upsilon$ be one-to-one onto and $\cd_\ell \colon \Upsilon \rightarrow
    \Upsilon$ be such that $ \cd_\ell(\cd(\langle \alpha_0,\alpha_\ell \rangle)) = \alpha_\ell$.  Clearly $J$ is a subfamily of ${\mathscr 
    P}(\lambda)$ and it is closed under subsets, so we need just to prove that it is closed under union of two.  So suppose $A = A_1 \cup A_2,A_\ell \in J$ for $\ell =1,2$.  As $J$ is closed under subsets, \wilog \, $A \cap A_2 = \emptyset$.  Let the sequence $\bar F^\ell = \langle F^\ell_\delta \colon \delta \in S \cap A_\ell \rangle$ exemplify $A_\ell \in J$.  Now we define $\bar F = \langle F_\delta \colon \delta \in S \cap A \rangle$ as follows:
    
    \begin{enumerate}
        \item[$\circledast$] if $\delta \in S \cap A$ and $f \in {}^{(C_\delta)}\Upsilon$ \underline{then}:  $\delta \in A_\ell \Rightarrow F_\delta(f) = F^\ell_\delta(\cd_\ell \circ f)$, i.e., $F_\delta(f) = F^\ell_\delta(f^\ell)$ where $\alpha \in C_\delta \Rightarrow f^\ell(\alpha) = \cd_\ell(f(\alpha))$.
    \end{enumerate}
    
    2) Straight as the sequence $\langle F_\delta \colon \delta \in S \rangle$ is not required to be nicely definable. 
    
    3) Guessing $\theta \times \theta$ is the same as guessing $\theta$.
\end{PROOF}

\centerline{$* \qquad * \qquad *$}
 
Well, are there good $(\lambda,\kappa,\kappa^+)$-parameters?  Some version of this is crucial for having interesting consequences of
Claim \ref{d.6}, but we shall have to quote
(on $\hat{I}[\lambda]$ see Definition \ref{0.5}; see on it \cite[\S1]{Sh:420}).

\begin{claim}\label{d.13} 
    1) If $S$ is  a stationary subset of a regular cardinal $\lambda$ and $S \in I [\lambda]$, see Definition \ref{0.5}, and 
    $(\forall \delta \in S)[{\text{\rm cf\/}}(\delta) = \kappa]$,  \underline{then}:

    \begin{enumerate}
        \item[$(\alpha)$]  for some club $E$ of $\lambda$, there is a good$^+ \,(S \cap E,\kappa,\kappa^+)$-{\rm \Bb}-parameter. (Of course, we can replace the $\kappa^+$ by any $\chi \le \lambda$) (recall that we do not required the $C_\delta$'s to be a closed subset of $\delta$),

        \item[$(\beta)$]  for some club $E$ of $\lambda$ there is a weakly good$^+ \,(S \cap E,\kappa,\kappa^+)$-{\rm \Bb}-parameter consisting\footnote{this means: the parameter $\bar C = \langle C_\delta \colon \delta \in S \cap E \rangle, C_\delta$ a club of $\delta$.} of clubs (this is of interest when $\kappa > \aleph_0$),

        \item[$(\gamma)$]   for some club $E$ of $\lambda$ there is a good $(S \cap E,\kappa,\kappa^+)$-{\rm \Bb}-parameter $\bar C = \langle  C_\delta \colon \delta \in S \cap E \rangle$ such that:

        \begin{enumerate}
            \item[(i)]   $\bar C$ is weakly good$^+$, i.e., $\alpha \in \nacc(C_{\delta_1}) \cap \text{\rm nacc}(C_{\delta_2})  \Rightarrow C_{\delta_1} \cap \alpha = C_{\delta_2} \cap \alpha,$

            \item[(ii)] $\bar C$ is good, i.e., for each $\alpha < \lambda$ we have $\lambda > |\{\alpha \cap C_\delta \colon \delta \in S, \alpha \in C_\delta\}|,$

            \item[(iii)] $\bar C$ is $\Theta_\lambda$-closed that is each $C_\delta$ is $\Theta_\lambda$-closed which means $\alpha = \sup(C_\delta \cap \alpha)  \, \wedge \,{ \text{\rm cf\/}}(\alpha) \in \Theta_\lambda \Rightarrow \alpha \in C_\delta$ where $\Theta_\lambda = \{\sigma \colon \sigma = { \text{\rm cf\/}}(\sigma)$ and $\alpha < \lambda \Rightarrow \trp_{\sigma}(\vert \alpha \vert) < \lambda\}.$
        \end{enumerate}
    \end{enumerate}
    
    2)  If $\kappa = {\text {\rm cf\/}}(\kappa),
    \kappa^+ < \lambda = {\text {\rm cf\/}}(\lambda)$, \underline{then} there is a  stationary $S \in \hat{I}[\lambda]$ with $(\forall \delta \in S) [{\text{\rm cf\/}}(\delta) = \kappa]$. 
    
    3) In part (1), moreover, we can find a good $(S \cap E,\kappa,\kappa^+)$-{\rm \Bb}-parameter and a $\lambda$-complete filter $D$ on $\lambda$ such that: for every club $E'$ of $\lambda$, the set $\{\delta \in S \cap E \colon C_\delta \subseteq E' \cap E\}$ belongs to $D$.
\end{claim}

\begin{PROOF}{\ref{d.13}}
    1)  By the definition of $I [\lambda]$ see Definition \ref{x.20.0109} and Definition \ref{d.5} above (see more in \cite[\S1]{Sh:420}) and for clause $(\gamma)$ the definition of $\trp_{\sigma}(\vert \alpha \vert)$.
      
    2)  By \cite[\S1]{Sh:420}. 
    
    3) By \cite[III]{Sh:g}.  
\end{PROOF}

\begin{definition}\label{d.14}   
    Let $\bar C$ be a $(\lambda,\kappa,\chi)$-\Bb-parameter.
    
    $\text{ }$  For a family ${\mathscr  F}$ of $\bar C$-colourings and ${\mathscr  P}
    \subseteq {}^\lambda 2$,  let $\id_{\bar C,{\mathscr  F},{\mathscr  P}}$ be 
    the set of $ W \subseteq \lambda $ such that:
    $ \{W \subseteq \lambda \colon \& \text{ for some } \mathbf F  \in {\mathscr  F} \text{ for every } \bar c \in {\mathscr  P}$   for some $  \eta \in {}^\lambda \lambda $ the set $\{\delta \in W \cap S \colon \mathbf F(\eta \, {\rest} \, C_\delta) = c_\delta\}$  is a non-stationary subset of $\lambda.$  
    
    If ${\mathscr  F}$  is the family of all $\bar C$-colouring we may omit it.
    If we write Def instead of ${\mathscr  F}$ this mean as in \cite[\S1]{Sh:576} (not used so the reader can ignore this).
\end{definition} 

We can strengthen Definition \ref{d.5}(2) as follows:

\begin{definition}\label{d.15} 
    We say the  $\lambda$-colouring $\mathbf F$ is \emph{$(S,\kappa,\chi)$-\Bb-good} \underline{if}:

    \begin{enumerate}
        \item[(a)]   $S \subseteq \{\delta < \lambda \colon \cf(\delta)  = \kappa\}$ is stationary, 

        \item[(b)]   we can find $E$ and $\langle C_\delta \colon \delta \in S \cap E \rangle$ such that:

        \begin{enumerate}
            \item[$(\alpha)$]  $E$ is a club of $\lambda,$

            \item[($\beta$)] $C_\delta$ is an unbounded subset of $\delta,|C_\delta| < \chi,$

            \item[$(\gamma)$] if  $\rho,\rho' \in {}^\delta \delta,\delta \in S \cap E$, and $\rho' \, {\rest} \, C_\delta = \rho \, {\rest} \, C_\delta$ then $\mathbf F(\rho') = \mathbf F(\rho),$

            \item[$(\delta)$]  for every $\alpha < \lambda$ we have $\lambda > |\{C_\delta \cap \alpha \colon \delta \in S \cap E \text{ and } \alpha \in C_\delta\}|,$

            \item[$(\varp)$]  $\delta \in S \Rightarrow \trp_{\cf(\delta)}(\vert \delta \vert) < \lambda$ or just $\delta \in S \Rightarrow \lambda > |\{C \colon C \subseteq \delta$ is unbounded  in $\delta$ and for every $\alpha \in C$ for some $\gamma \in S$ we have $C \cap \alpha = C_\gamma \cap \alpha$ and $\alpha \in C_\gamma\}|$.
        \end{enumerate}
    \end{enumerate}
\end{definition} 

We can sum up (relying again on some quotations).

\begin{conclusion}\label{d.16} 
    If  $\lambda = \text{cf}(2^\mu)$ and we let $\Theta \coloneqq \{\theta \colon \theta =
    \text{cf}(\theta)$ and $(\forall \alpha < 2^\mu) (\trp_{\theta}(\vert \alpha \vert) < 2^\mu)\}$ \underline{then}:

    \begin{enumerate}
        \item[(a)]  $\Theta$ is ``large'' e.g. contains every large enough $\theta \in \text{ Reg } \cap \beth_\omega$ if $\beth_\omega < \lambda,$  

        \item[(b)] if $\theta \in \Theta_\lambda \wedge \theta^+ < \lambda$ then there is a stationary $S \in \hat{I}[\lambda]$ such that $\delta \in S \Rightarrow \cf(\delta) = \theta,$

        \item[(c)]  if $\theta,S$ are as in clause (b) above  \underline{then} for some club $E$ of $\lambda,$ we have that: 

        \begin{enumerate}
            \item[(i)] there is a weakly good $(\lambda,\kappa,\kappa^+)$-\Bb-parameter $\bar C = \langle C_\delta \colon \delta \in S \cap E \rangle,$

            \item[(ii)] there is a good \, $(\lambda,\kappa,\kappa^+)$-\Bb-parameter  $\langle C_\delta \colon \delta \in E \cap S \rangle,$

            \item[(iii)] there is a good$^+ \,(\lambda, \kappa,\kappa^+)$-\Bb-parameter  $\bar C = \langle C_\delta \colon \delta \in E \cap S \rangle.$
        \end{enumerate}

        \item[(d)] for $\theta,S,\bar C$ as in clause (b), (c) above and $\Upsilon \le 2^\mu$, if $\mathbf F$ is a $(\bar C,\Upsilon,\theta)$-colouring and $D = D_\lambda$ (or $D$ is a $\mu^+$-complete filter on $\lambda$ extending $D_\lambda$), and lastly $S \in D^+$ \underline{then} there is a $\mathbf F$-\midia sequence for $\bar C$ and there is an $D-\mathbf F$-\Bb-sequence; see Definition \ref{d.5}(3), (3A),

        \item[(e)]  if $\theta,S,\bar C$ as in clauses (b), (c): $\Upsilon \le 2^\mu$,  \underline{then} $\bar C$ has the $(\Upsilon,\theta)$-property,

        \item[(f)]   for $\theta,S$ as in clause (b) above $\Upsilon \le 2^\mu$, if $\theta > \aleph_0$ $\bar C' = \langle C'_\delta \colon \delta \in S \rangle,C_\delta$ a stationary subset of $\delta$ for each $\delta$ then we can find $\bar C'' = \langle C''_\delta \colon \delta \in E \cap S \rangle$ such that $C''_\delta$ is a club of $\delta$ for each $\delta$ and $\langle C'_\delta \cap C''_\delta \colon \delta \in S \rangle$ has the $(\Upsilon,\theta)$-\Bb-property (really $\bar C''$ does not depend on $\bar C'$),

        \item[(g)]   in clause (f) if $\theta = \theta^\sigma,\sigma \ge \kappa,2^\chi \le \theta$ we can conclude also that $\langle C'_\delta \cap C''_\delta \colon \delta \in S \rangle$ has the $(\chi,\sigma)$-MD$_2$-property.
    \end{enumerate}
\end{conclusion}

\begin{PROOF}{\ref{d.16}}
     Clause (a) holds by \cite{Sh:460}, clause (b), (c) holds by Claim \ref{d.13}, clause (d) by the major  Claim \ref{d.6} and clause (e) is a rephrasing of clause (d), and lastly for clause (f) use above case (iii) of clause (c) and use the obvious monotonicity (and $C'_\delta \cap C''_\delta$ is unbounded in $\delta$ being stationary).
\end{PROOF}

Note some open questions.

\begin{question} \label{d.18} 
    1) Assume  $\lambda$ is a weakly inaccessible which is not strong limit (i.e. a limit regular uncountable cardinal but not strong limit) and $\langle
    2^\theta \colon \theta < \lambda \rangle$ is not eventually constant. 
    
    Does $\lambda$ have the weak diamond?  Other consequences?   In this case cf$(2^{< \lambda}) = \lambda$ as $\langle 2^\theta \colon \theta < \lambda \rangle$ witness it hence $2^\lambda > 2^{< \lambda} \coloneqq 
    \Sigma\{2^\theta \colon \theta < \lambda\}$.  If $\lambda = 2^{< \lambda}$ then 
    $\lambda = \lambda^{< \lambda}$ and such $\lambda$'s are investigated in
    \cite{Sh:755} (for weak diamond on many stationary sets).  
    
    2) Consider generalizing 
    \cite[AP,\S3]{Sh:f} to our present context.  Also, we may consider $\langle c_\eta \colon \eta \in \lambda^{+n} \times \lambda^{+(n-3)} \times \ldots \times \lambda
    \rangle$ where $c_\eta \in \theta$.
\end{question} 

\begin{claim}\label{d.n.21}
    Assume, 
    
    \begin{enumerate}
        \item[$\circledast_{\mu,\theta}$] 
        
            \begin{enumerate} 
            
            \item[(i)] $\theta \le \mu,$
            
            \item[(ii)] if ${\mathscr  F} \subseteq {}^\mu \theta$ has cardinality $< 2^\mu$ then there is $\rho \in {}^\mu \theta$ such that $(\forall f \in {\mathscr  F})(\exists \varepsilon < \mu)(\rho(\varepsilon) = f(\varepsilon))$.
        \end{enumerate} 
    \end{enumerate} 
    
    1) $\Sep(\mu,\mu,\theta,2^\mu)$ 
    hence if $2^\mu$ is a regular cardinal then $\Sep(\mu,\theta)$ holds so clause (f) of Claim \ref{d.6} is satisfied. 
    
    2) Even if $2^\mu$ is singular still in Claim \ref{d.6} we can replace clause (f) in the assumption by $\circledast_{\mu,\theta}$.
\end{claim} 

\begin{remark} 
    On relatives of $\circledast_{\mu,\theta}$ see ${\mathfrak e}_\kappa$, see \cite{Sh:713}.
\end{remark}

\begin{PROOF}{\ref{d.n.21}}
    1) Clearly ${}^\mu \theta$ has cardinality $2^\mu$, so let $\langle \nu_\alpha \colon \alpha < 2^\mu \rangle$ be a list without
    repetitions of ${}^\mu \theta$.  For $\alpha \le 2^\mu$, let $\Gamma_\alpha = 
    \{\nu_\beta \colon \beta < \alpha\}$ and we choose $\bar f^\alpha$ by induction on $\alpha \le 2^\mu$ such that: 
    
    \begin{enumerate} 
        \item[(a)]  $\bar f^\alpha = \langle f^\alpha_\varepsilon \colon \varepsilon < \mu \rangle,$
        
        \item[(b)]  $f^\alpha_\varepsilon$ is a function from $\Gamma_\alpha$ to $\theta,$
        
        \item[(c)]  $f^\alpha_\varepsilon$ is increasing continuous in $\alpha$, i.e., $\beta < \alpha \Rightarrow f^\beta_\varepsilon = f^\alpha_\varepsilon \, {\rest} \, \text{ Dom}(f^\beta_\varepsilon),$
        
        \item[(d)]  for every $\gamma < \beta < \alpha$ the set $\{\varepsilon < \mu \colon f^\alpha_\varepsilon(\nu_\beta) = \nu_\gamma(\varepsilon)\}$ is not empty.
    \end{enumerate}
        
    If we succeed, then we let $f_\varepsilon =
    f^{2^\mu}_\varepsilon$ that is $f_\varepsilon = \bigcup\{f^\alpha_\varepsilon \colon \alpha < 2^\mu\}$ for $\varepsilon < \mu$ and $\langle f_\varepsilon \colon\varepsilon < \mu \rangle$ is clearly as required.
    
    For clause (d) it is enough to consider the case $\alpha = \beta +1$, as the non-trivial case in the induction is $\alpha = \beta +1$.
    We can choose $\langle f^{\beta +1}_\varepsilon(\nu_\beta) \colon \varepsilon
    < \mu \rangle \in {}^\mu \theta$ such that it holds by the assumption $(*)$. 
    
    2) In the proof of Claim \ref{d.6}, now $\rho \in {}^\mu \theta \Rightarrow |\Sol_\rho| < 2^\mu$ does not suffice.  However,
    before choosing $\langle f_\varepsilon \colon \varepsilon < \mu \rangle$ we
    are given $\mathbf F$ and we can choose $\cd$, $\langle \cd_\varepsilon \colon \varepsilon < \mu \rangle$ and define $T,\langle \eta_\alpha \colon \alpha < 2^\mu \rangle,T_{<\alpha}$ and $\langle {\mathscr 
    P}_\delta \colon \delta \in S \rangle$.  Note that $\alpha < \lambda \Rightarrow \trp_{\kappa}(\vert h(\alpha) \vert) < 2^\mu$ hence $\alpha < \lambda \Rightarrow |\cup\{{\mathscr  P}_\delta \colon \delta \in S \cap \alpha\}| < 2^\mu$.  We can also choose $\langle \nu_{\eta,\varepsilon} \colon \eta \in {\mathscr  P}_\delta, \varepsilon < \mu
    \rangle$ and define  $\rho_\eta \in {}^\mu \theta$ for $\eta \in  {\mathscr  P}_\delta,\delta \in S$.
    
    Now we can phrase our demands on $\langle f_\varepsilon \colon \varepsilon < \mu \rangle$, which is in fact weaker than the full generality of $\Sep(\mu,\theta)$;  it is
    
    \begin{enumerate} 
        \item[$(\boxdot )$] 
        ${}^\mu \theta \ne \{\nu \in {}^\mu \theta$: for some $\eta \in {\mathscr  P}_\delta$ we have $(\forall \varepsilon < \mu)(f_\varepsilon(\nu) \ne \rho_\eta(\varepsilon))\}$.
    \end{enumerate} 
    
    For this, it is enough that for the list $\langle \nu_\alpha \colon \alpha <
    2^\mu \rangle$ chosen in the proof of part (1),  for  each $\delta \in S$,
    for some $\gamma < 2^\mu$ we have $\{\rho_\eta \colon \eta \in 
    {\mathscr  P}_\delta\} \subseteq \{\nu_\alpha \colon\alpha < \gamma\}$ 
    which is easy to obtain.
\end{PROOF}

\begin{observation}\label{d.n.22} 
    If $\mu \ge \aleph_0,\theta=2$, then
    $\circledast_{\mu,\theta}$ holds.
\end{observation}  

\begin{PROOF}{\ref{d.n.22}}
    If ${\mathscr  F} \subseteq {}^\mu \theta,|{\mathscr  F}| < 2^\mu$, let $\nu \in {}^\mu \theta \backslash {\mathscr  F}$ and let $\rho = \langle 1 - \nu(\varepsilon) \colon \varepsilon < \mu \rangle$.
\end{PROOF}

\newpage

\section{Upgrading to larger cardinals}\label{2} 

In \S1 we get the $\theta$-\midia \ property for  regular cardinals $\lambda$ of the form cf$(2^\mu)$ on ``most" $S \subseteq S^\lambda_\kappa$ for ``many" $\kappa = \cf(\kappa)$ (and $\theta$).   We would like to prove Theorem \ref{0.1}. For this, we first give an upgrading, which is a stepping-up claim
(for MD$_\ell$) given as Claim \ref{x.20} below.  Then we prove that it applies to $\lambda_2$ (for quite a few $\lambda_1$), when  $\lambda_2$ is a successor of regular, then when $\lambda_2$ is a successor of singular, and lastly when $\lambda_2$ is weakly inaccessible (not strong limit).   We can upgrade each of the variants from \S1, but we shall concentrate on the one 
from Definition \ref{d.5}(5A) because of Claim \ref{d.12}.  The main result of this section is claimed in Claim \ref{x.24}.

\begin{question}\label{x.17}
    What occurs to $\lambda_2$ strongly inaccessible, does it have to have some \Bb-property?
\end{question} 

Possibly combining \S1 and \S2 we may eliminate the $\lambda_2 > 2^{\lambda_1 + \Upsilon + \theta}$ in Claim \ref{x.20}, but does not seem urgent now.

\begin{claim}\label{x.20} 
    1) Assume
    
    \begin{enumerate}  
        \item[(a)]  $\lambda_1 < \lambda_2$ are regular cardinals  such that $\lambda_2 > 2^{\lambda_1+\Upsilon^{< \chi} +\theta}$ (this last condition can be waived for ``nice" colourings, see Definition \ref{x.20.gr} below),
        
        \item[(b)] for some stationary 
        $S_1 \subseteq S^{\lambda_1}_\kappa$, some  sequence $\bar C = \langle C^1_\delta \colon \delta \in S_1 
        \rangle,$ 
        
        \item[(c)]  $\lambda_2 > \lambda^+_1,\lambda_2$ is a successor of
        regular. 
    \end{enumerate}  
    
    \underline{then} for some stationary $S_2 \subseteq S^{\lambda_2}_\kappa$, there is $\bar C = \langle C_\delta \colon \delta \in S_2 \rangle$, a good \footnote{the goodness is a bonus.} $(\lambda_2,\kappa,\chi)$-{\rm \Bb}-parameters which has the $(D_{\lambda_2}, \Upsilon, \theta)$-{\rm \Bb}-property.  

    2) We can in part (1) replace clause (c) with any of the following clauses:

    \begin{enumerate}
        \item[(c)$'$]   there is a $\le \lambda_1$ square on some stationary $S \subseteq S^{\lambda_2}_{\lambda_1}$, (see Definition \ref{x.20.0109}(1A); this holds if $\lambda_2 = \mu^+, \mu^{\lambda_1} = \mu),$

        \item[(c)$''$] $\lambda_2$ has a $(\lambda_1,S_1)$-square (where $S_1$ is from clause (b)) which means:

        \begin{itemize}
            \item[$\boxplus^{\lambda_2}_{\lambda_1,S_1}$]   $\lambda_1 < \lambda_2$ are regular, $S_1 \subseteq \lambda_1$ is stationary and there are $S_2 \subseteq S^{\lambda_2}_{\lambda_1}$ stationary and $\bar C^2 = \langle C^2_\delta \colon \delta \in S_2 \rangle$ satisfying $C^2_\delta \subseteq \delta  = \sup(C^2_\delta),C^2_\delta$ is closed, {\text{\rm otp\/}}$(C^2_\delta) = \lambda_1$ and if $\alpha \in  C^2_{\delta_1} \cap C^2_{\delta_2},\{\delta_1,\delta_2\} \subseteq S_2$ and $\alpha = \sup(C^2_{\delta_1} \cap \alpha)$  and {\text{\rm otp\/}}$(\alpha \cap C^2_{\delta_2}) \in S_1$, \underline{then} $\alpha \cap C^2_{\delta_2}  = \alpha \cap C^2_{\delta_1}$.
        \end{itemize}

        \item[(c)$'''$]   $\lambda_2$ has a $(\lambda_1,\bar C_1)$-square
        $\bar C_2$ (where $\bar C_1$ is from clause $(b)$) which means:

    \begin{enumerate}
        \item[$\boxplus^{\lambda_2}_{\lambda_1,\bar C_1,\bar C_2}$] 

        \begin{enumerate}
            \item[$(\alpha)$] $\bar C_2 = \langle C^2_\delta \colon \delta \in S_2 \rangle$, 

            \item[$(\beta)$] $S_2 \subseteq
            S^{\lambda_2}_{\lambda_1}$ is stationary,

            \item[$(\gamma)$] $C^2_\delta$ is a closed unbounded subset of $\delta$ of order type $\lambda_1$,

            \item[$(\delta)$]  if $\alpha \in C^2_{\delta_1} \cap C^2_{\delta_2}$ (so $\{\delta_1,\delta_2\} \subseteq S_2)$, {\text{\rm otp\/}} $(\alpha \cap C^2_{\delta_1}) \in S_1$,   and $\alpha = \sup(C^2_{\delta_1} \cap \alpha),$ \underline{then} {\text{\rm otp\/}}$(\alpha \cap C^2_{\delta_1}) = {\text {\rm otp\/}}(\alpha  \cap C^2_{\delta_2})$, call it  $\gamma$ and $\{\beta \in C^2_{\delta_1} \colon {\text{\rm otp\/}}(\beta \cap  C^2_{\delta_1}) \in C^1_\gamma\} = \{\beta \in C^2_{\delta_2}:{\text{\rm otp\/}} (\beta \cap C^2_{\delta_2}) \in C^1_\gamma\}$.
        \end{enumerate}
    \end{enumerate}
    \end{enumerate}

    3) Assume $\ell \in \{0,2\}$ and,

    \begin{enumerate}
        \item[(a)] $\lambda_1 < \lambda_2$ are regular uncountable,

        \item[(b)] some $(\lambda_1, \kappa, \chi)$-{\rm \Bb}-parameter has the $\theta$-{\rm MD}$_\ell$-property,

        \item[(c)]  at least one of the variants of clause (c) above holds.
    \end{enumerate}

    \underline{Then}, some $(\lambda_2,\kappa,\chi)$-{\rm \Bb}$_\ell$-parameter has the
    $\theta$-{\rm MD}$_\ell$-property.
\end{claim}  

\begin{definition}\label{x.20.gr} 
    We say that a 
    $(\Upsilon,\theta)$-colouring
    $\mathbf F = \langle F_\delta:\delta \in S \rangle$ for a
    $\lambda$-\Bb-parameter $\bar C = \langle C_\delta \colon \delta \in S
    \rangle$ is \emph{nice} when:

    \begin{enumerate}
        \item[$\circledast$] if $\delta_1,\delta_2 \in S$ and
        $\otp(C_{\delta_1}) = \otp(C_{\delta_2})$ and $\eta_1 \in
        {}^{(C_{\delta_1})}\Upsilon,\eta_2 \in {}^{(C_{\delta_2})}\Upsilon$
        and $\alpha_1 \in C_{\delta_1} \wedge \alpha_2 \in
        C_{\delta_2} \wedge \otp(C_{\delta_1} \cap \alpha_1) = \otp(C_{\delta_2} \cap \alpha_2) \Rightarrow \eta_1(\alpha_1) =
        \eta_2(\alpha_2)$, \underline{then} $F_{\delta_1}(\eta_1) =
        F_{\delta_2}(\eta_2)$.
    \end{enumerate}
\end{definition}  

\begin{PROOF}{\ref{x.20.gr}}
    1) If clause $(c)$ holds then clause $(c)'$ of part (2) holds by \cite[\S4]{Sh:351}, so it is enough to prove part (2). 
    
    2) Trivially clause $(c)'$ implies clause $(c)''$ and clause $(c)''$
    implies clause $(c)'''$; so assume $(c)'''$; and let $\bar C^2 = \langle
    C^2_\delta \colon \delta \in S_2 \rangle$ be as there.

    Let $\bar h = \langle h_\delta:\delta \in S_2
    \rangle$ be such that $h_\delta$ is an increasing continuous function from
    $\lambda_1$ onto $C^2_\delta$ recalling clause $(\gamma)$, i.e.,
    $C^2_\delta$ is a closed unbounded subset of $\delta$ of order type $\lambda_1$.  Let $S = \{h_\delta(i) \colon \delta \in S_2$
    and $i \in S_1\}$ and if $\alpha \in S,\alpha = h_\delta(i),\delta \in
    S_2$ and $i \in S_1$ then we let $C^*_\alpha = 
    \{h_\delta(j) \colon j \in C^1_i\}$; and let
    $\bar C^* = \langle C^*_\alpha:\alpha \in S \rangle$.  Why is
    $C^*_\alpha$ (hence $\bar C^*$) well defined?  We shall check that,  if
    $\alpha \in S,\alpha = h_{\delta_\ell}(i_\ell), \delta_\ell \in
    S_2,i_\ell \in S$ for $\ell=1,2$ then $\{h_{\delta_1}(j) \colon j \in C^1_{i_1}\} =
    \{h_{\delta_2}(j) \colon j \in C^1_{i_2}\}$; this holds by clause $(\delta)$
    of the assumption $(c)'''$.   We should check the conditions 
    for ``$\bar C^*$ is a good $(\lambda_2,\kappa,\chi)$-\Bb-parameter with
    the $(D_{\lambda_2},\Upsilon,\theta)$-\Bb-property.
    
    The main point is: 
    assume that $\mathbf F$ is a $(\bar C^*,\Upsilon,\theta)$-colouring, find
    a $\bar c$ as required. 
    
    For each $\delta \in S_2$ we define a 
    $(\bar C^1,\Upsilon,\theta)$-colouring $\mathbf F^\delta$ by $\mathbf
    F^\delta(\eta) = \mathbf F(\eta \circ h^{-1}_\delta)$ that is for $i \in
    S_1$ and $\eta \in {}^{C^1_i}\Upsilon$ we define $\mathbf F^\delta(\eta)$
    as $\mathbf F(\nu)$ where $\nu \in {}^{C^*_{h_\delta(i)}}\Upsilon$, where
    $\nu$ is defined by $\forall \beta \in C^*_{h_\delta(i)}$, we let
    $\nu(\beta) = \eta(h^{-1}_\delta(\beta))$.  As\footnote{only
    place this is used, if $\mathbf F$ is nicely defined this is not necessary.}
    $\lambda_2 > 2^{\lambda_1+|\Upsilon|^{< \chi}+\theta}$
    because $\theta^{\Sigma\{|\Upsilon|^{|C'_\delta|} \colon \delta \in S_1\}}
    \le 2^{\theta\lambda_1\Sigma\{|\Upsilon|^{|\alpha|} \colon \alpha \le
    \chi\}} = 2^{\lambda_1 + |\Upsilon|^{< \chi}} + \theta$, for some 
    stationary $S^*_2 \subseteq S_2$ we have
    $\delta \in S^*_2 \Rightarrow \mathbf F^\delta = \mathbf F^*$, so $\mathbf
    F^*$ is a $(\bar C^1,\Upsilon,\theta)$-colouring, hence it has a good
    $(D_{\lambda_1},\Upsilon,\theta)$-\Bb-sequence $\bar c^1 = \langle c^1_i:i \in S_1 \rangle$. 
    
    We define $\bar c^* = \langle c^*_\alpha:\alpha \in S \rangle$ by
    $\delta \in S_2 \wedge i \in S_1 \Rightarrow c^*_{h_\delta(i)} =
    c^1_i$.  By clause $(\delta)$ of $(c)'''$ it is well defined.  Now for
    any $\eta \in {}^{(\lambda_2)} \Upsilon$ and club $E_2$ of $\lambda_2$,
    we can find $\delta \in S_2 \cap \text{ acc}(E_2)$, hence $E_1 =
    \{\alpha < \lambda_1:h_\delta(\alpha) \in E_2\}$ is a club of
    $\lambda_1$ and let $\eta_1 \in {}^{(\lambda_1)} \Upsilon$ be
    $\eta_1(\alpha) = \eta(h_\delta(\alpha))$, so by the choice of $\bar
    c^1$ for some $\delta_1 \in S_1 \cap E_1,c^1_{\delta_1} = \mathbf
    F^*(\eta_1 \, {\rest} \, C^1_{\delta_1})$ and let $\delta_2 =
    h_\delta(\delta_1)$, so $\delta_2 \in E_2$ and $c^*_{\delta_2} = c^1_i =
    \mathbf F^*(\eta_1 \, {\rest} \, C^1_{\delta_1}) = \mathbf F^\delta(\eta_1
    \, {\rest} \, C^1_{\delta_1}) = \mathbf F(\eta \, {\rest} \,
    C^*_{\delta_2})$, so we are done. 
    
    3) Similar only easier (or see the proof of Claim \ref{x.20.9}).
\end{PROOF}

\begin{claim}\label{x.21A} 
    We can strengthen the conclusion of
    Claim \ref{x.20}(1) and Claim \ref{x.20}(2) to:  
    we can find a sequence $\langle S_{2,i} \colon i < \lambda_2 \rangle$ of pairwise disjoint stationary subsets of $\lambda_2$ and for each $i < \lambda_2$ a sequence  $\langle C^i_\delta \colon \delta \in S_{2,i} \rangle$ which is
    a $(\lambda_2, \kappa, \chi)$-{\rm \Bb}-parameter having the $(D_{\lambda_2}, \Upsilon, \theta)$-{\rm \Bb}-property.
\end{claim} 

\begin{PROOF}{\ref{x.21A}}
    By the proof of Claim \ref{x.20}, \wilog \, clause $(c)''$ of Claim \ref{x.20}(2) and clauses (a), (b) of Claim \ref{x.20}(1) hold.
    For some $\zeta < \lambda_1$ for un-boundedly many $i < \lambda_2$ the set $S_{2,i} \coloneqq \{\delta \in S_2 \colon i$ is the $\zeta$-th member of $C^2_\delta\}$ is stationary.
\end{PROOF} 

\begin{claim}\label{x.21} 
    1) If $\mu$ is a strong limit, 
    $\lambda$ is a  successor cardinal, $\lambda = \cf(\lambda) > 2^\mu$, moreover $\lambda > 2^{2^\mu}$ (superfluous for nice colourings, e.g., for MD$_\ell,\ell \in \{0,2\}$), \underline{then} for every large enough regular $\kappa < \mu$, 
    some good$^+ \,(\lambda,\kappa,\kappa^+)$-\Bb-parameter  $\bar C$ has the $(2^\mu,\theta)$-\Bb-property for every 
    $\theta < \mu$. 
    
    2) We can add: $\bar C$ has the $\theta$-MD$_\ell$-property for $\ell
    = \{0,2\}$. 
    
    3) We can find a sequence $\langle S_i \colon i < \lambda \rangle$ of
    pairwise disjoint stationary subsets of $S^\lambda_\kappa$ such that
    for each $i < \lambda$ there is $\bar C^i = \langle C^i_\delta:\delta
    \in S \rangle$ as in parts (1) and (2). 
\end{claim} 

\begin{PROOF}{\ref{x.21}}
    1) Let $\kappa_0 < \mu$ be large enough such  that $[\theta \in \{\lambda, 2^\mu\} \, \wedge \,\alpha < \theta \, \wedge \,\kappa_0
    \le \kappa = \cf(\kappa) < \mu] \Rightarrow |\alpha|^{[\kappa]} < \theta$, (see Definition \ref{d.3}, exists by \cite{Sh:460}).
    
    By Claim \ref{d.6}, for any regular $\kappa \in [\kappa_0,\mu)$ and $\theta < \mu$ we have the result for cf$(2^\mu)$ in place of
    $\lambda$.  If $\lambda$ is a successor of regular we can  apply Claim \ref{x.20}(1), i.e., for clause $(c)$ there is a square as required by \cite[\S4]{Sh:351}.   In the case that $\lambda$ is a successor of a singular cardinal, by  \ref {x.2.9} below, clause $(c)''$ of Claim \ref{x.20}(2) holds for $\lambda_2 = \lambda,\lambda_1 = \text{cf}(2^\mu)$ so we can apply Claim \ref{x.20}. 
    
    2) Follows by Claim \ref{d.12}. 
    
    3) In the proof of parts (1) + (2) above, instead of Claim \ref{x.20}  we use Claim x.21A.  
\end{PROOF}

Below we generalize \emph{Enge\l king Kar\l owicz  theorem} (see \cite{EK}) (which says that if $\mu = \mu^\kappa, \lambda = 2^\mu$ then there are $f_\varepsilon \colon \lambda
\rightarrow \mu$ for $\varepsilon < \mu$ such that every partial function from $\lambda$ to $\mu$ with domain of cardinality $\le \kappa$, is extended by $f_\varepsilon$ for some $\varepsilon < \mu$); we do more than is necessary for the \midia.

\begin{claim}\label{x.23} 
    1) Assume
    
    \begin{enumerate} 
        \item[(a)]  $\mu = {\text {\rm cov\/}}(\mu,\theta^+,\theta^+,\kappa^+),$
        
        \item[(b)]  $2^\theta \le \mu,$
        
        \item[(c)]   $\mu < \lambda \le 2^\mu,$
    \end{enumerate}
    
    \underline{then} there is $\bar F$ such that:
    
    \begin{enumerate} 
        \item[$(\alpha)$] $\bar F = \langle F_\varepsilon \colon \varepsilon < \mu \rangle,$
        
        \item[($ \beta $)]    each $F_\varepsilon$ is a function from $\lambda$ to $\mu,$ 
        
        \item[($ \gamma $)]  if $u \in [\lambda]^{\le \theta}$ \underline{then} we can find $\bar E$ such that: 
        
        \begin{enumerate} 
            \item[(i)]  $\bar E = \langle E_i \colon i < \kappa \rangle$ is a sequence of equivalence relations on $u,$
            
            \item[(ii)]  $\alpha \in u \Rightarrow \{\alpha\} = \cap \{\alpha/E_i \colon i < \kappa\},$
        
            \item[(iii)]  if $f \colon u \rightarrow \mu$ respects some $E_i$ (i.e. $\alpha E_i \beta \Rightarrow f(\alpha) = f(\beta))$  \underline{then} we can find a partition $\langle u_\xi \colon \xi < \kappa \rangle$ of $u$ such that $\xi < \kappa \Rightarrow f \, {\rest} \, u_\xi  \subseteq F_\varepsilon$ for some $\varepsilon < \mu$. 
        \end{enumerate} 
    \end{enumerate} 
    
    2) In part (1) we  can strengthen clause $(\gamma)^-$ to $(\gamma)$ below if
    $(d)_\ell$ for some $\ell \in \{1,\dotsc,5\}$ below holds, where:

    \begin{enumerate}
        \item[$(\gamma)$] if $u \in [\lambda]^{\le \theta}$ and $f \colon u \rightarrow \mu$, \underline{then} we can find $\bar u = \langle u_i \colon I < \kappa \rangle$ such that:

        \begin{enumerate}
            \item[(i)]   $u = \cup\{u_i \colon i < \kappa\}$ and,

            \item[(ii)] for every $i < \kappa$ for some $\varepsilon < \mu$ we have $f \, {\rest} \, u_i \subseteq F_\varepsilon.$  
        \end{enumerate}
    \end{enumerate}

    \begin{enumerate}
        \item[(d)$_{1}$] there is a sequence $\bar f = \langle f_\alpha \colon \alpha < \lambda \rangle$ of pairwise distinct functions from $\mu$ to $\mu$ such that: if $u \in [\lambda]^{\le \theta}$ \underline{then}  we can find sequences $\bar u = \langle u_i \colon i < \kappa \rangle$ and  $\langle B_i \colon i < \kappa \rangle$ such that  $u = \cup \{u_i \colon i < \kappa\},B_i \in [\mu]^{< \aleph_0}$ for $i < \kappa$ and for each $i < \kappa$ the sequence $\langle f_\alpha \, {\rest} \, B_i \colon \alpha \in u_i \rangle$ is with no repetitions, 

        \item[(d)$_{2}$] there is a sequence $\bar f = \langle f_\alpha \colon \alpha < \lambda \rangle$ of functions from $\kappa$ to $\mu$ and an ideal $J$ on $\kappa$ such that $\bar f$ is $\theta^+$-free modulo $J$, i.e., $(\theta^+,J)$-free (which means that 
        for $u \in [\lambda]^{\le \theta}$, we can find a sequence $\bar a = \langle a_\alpha \colon \alpha\in u \rangle$ of subsets of $\kappa$  such that $a_\alpha = \kappa$ {\text{\rm mod\/}} $J$  and $i \in a_\alpha \cap a_\beta \, \wedge \,\alpha \ne \beta \Rightarrow f_\alpha(i) \ne f_\beta(I)),$

        \item[(d)$_{3}$]  like $(d)_2$ but $a_\alpha \in J^+$ that is $a_\alpha \notin J,$ 

        \item[(d)$_{4}$]  $(\forall \alpha < \mu)(|\alpha|^\theta < \mu)$ and $\lambda < {\text{\rm pp\/}}_J(\mu)$  for some ideal $J$ on $\kappa,$

        \item[(d)$_{5}$]   $\mu = \mu^{\aleph_0} < \lambda < \mu^\kappa$.
    \end{enumerate}
    
    3) If we weaken the conclusion of the part (2), replacing $(\gamma)$ by: 

    \begin{enumerate}
        \item[$(\gamma)'''_{\theta,\kappa}$] for every club $E$ of $\lambda$ for some closed $u \subseteq E$ of order type $\theta$  there is $\langle u_i \colon i < \kappa \rangle$ such that $u_i \subseteq u$ and  $\zeta \in u \, \wedge \,{\text {\rm cf\/}}(\zeta) > \kappa \Rightarrow \zeta \in \cup\{u_i \colon i < \kappa\}$ and $(\forall i < \kappa)(\exists \varepsilon < \mu)[f \, {\rest} \, u_i \subseteq F_\varepsilon),$ 
    \end{enumerate}
     
    \underline{then} we can weaken the assumption in a parallel way,
    e.g.

    \begin{enumerate}
        \item[$(d)^-_{1,\theta,\kappa}$]   there is  $\bar f = \langle f_\alpha \colon \alpha < \lambda \rangle, f_\alpha \in {}^\mu \mu$ such that for any club $E$ of $\lambda$ we can find a closed $e \subseteq E$ of order type $\theta +1$, and $\bar u = \langle u_i \colon i < \kappa \rangle,\bar B = \langle B_i \colon i < \kappa \rangle$ such that $B_i \in [\mu_1] < \aleph_0, u_i \subseteq \theta$ and $\zeta \in u \, \wedge \,{\text {\rm cf\/}}(\zeta) \ge \kappa \Rightarrow \zeta \in \cup\{u_i \colon i < \kappa\}$ and for  each $i, \langle f_\alpha \, {\rest} \, B_i \colon \alpha \in u_i \rangle$  is with no repetitions (and in $(d)_3$ it becomes easy).
    \end{enumerate}
\end{claim} 

\begin{remark}
    We can replace $\lambda$ (as a domain) by $[\lambda]^2$, etc.
\end{remark}

\begin{PROOF}{\ref{x.23}}
    1)  Let ${\mathscr  A}$ be a family of $\le \mu$ subsets of $\mu$ each of cardinality $\le \theta$ such that any $A \in
    [\mu]^{\le \theta}$ is included in the union of $\le \kappa$ many of them; it exists by assumption (a).  By assumption (b) \wilog \, $A \in {\mathscr  A} \, \wedge \, B \subseteq A \Rightarrow B \in {\mathscr  A}$ so we replace ``included" by ``equal" and also \wilog \, $A_1 \in {\mathscr  A} \, \wedge \,  A_2
    \in {\mathscr  A} \Rightarrow A_1 \cup A_2 \in {\mathscr  A}$.  Let pr be a
    pairing function on $\mu$.
    
    As $\mu < \lambda \le 2^\mu$ there is a sequence $\langle f^*_\zeta \colon \zeta < \lambda \rangle$ of pairwise distinct functions from $\mu$ to $\{0,1\}$. Let
    $ {\mathscr  Y} \coloneqq \{(A,B,\bar g,\bar \beta) \colon  B \in {\mathscr  A},\bar g = \langle g_i \colon i < \theta \rangle 
    \& g_i \in {}^B 2 \text{ and } i < j < \theta \Rightarrow g_i \ne g_j  \text{and }  \bar \beta = \langle \beta_i \colon i < \theta\rangle  \text{ satisfies } \{\beta_i \colon i < \theta\} \subseteq A \in {\mathscr  A}\}. $
    
    Clearly $|{\mathscr  Y}| \le \mu \times 2^\theta = \mu$.  For each $x = (A,\bar g,\bar \beta) \in {\mathscr  Y}$ 
    we define a function $F_x \colon \lambda \rightarrow \mu$ as follows:
    
    $F_x(\zeta)$ is $\beta_i$ if $g_i \subseteq f^*_\zeta$ and is 0 if
    there is no such $i$ (clearly we cannot have $i_0 \ne i_1 < \theta \, \wedge \, 
    g_{i_0} \subseteq f^*_\zeta \, \wedge \,  g_{i_1} \subseteq f^*_\zeta$).
    
    Let $\bar F = \langle F_\varepsilon \colon \varepsilon < \mu \rangle$ list
    $\{F_x \colon x \in {\mathscr  Y}\}$, clearly $\bar F$ satisfied demand $(\alpha)
    + (\beta)$.  As for demand $(\gamma)^-$ let $u \in [\lambda]^{\le
    \theta}$ be given.  For any $\zeta_1 \ne \zeta_2 \in u$ let $\alpha =
    \alpha(\zeta_1,\zeta_2) < \mu$ be minimal such that
    $f^*_{\zeta_1}(\alpha) \ne f^*_{\zeta_2}(\alpha)$, exists as $f^*_{\zeta_2}
    \ne f^*_{\zeta_1} \in {}^\mu 2$ and let $B =
    \{\alpha(\zeta_1,\zeta_2):\zeta_1 < \zeta_2$ are from $u\}$.  By the
    choice of ${\mathscr  A}$ we can find pairwise disjoint
    $B_\xi \in {\mathscr  A}$ for $\xi <
    \kappa$ such that $B = \bigcup \{B_\xi \colon \xi < \kappa\}$.  Lastly, we
    define $\bar E = \langle E_\xi \colon \xi < \kappa \rangle$ by

    \begin{enumerate}
        \item[$\circledast$] for $\zeta_1,\zeta_2 \in u, \zeta_1 E_\xi \zeta_2$ \underline{iff} $f^*_{\zeta_1} \, {\rest} \, B_\xi = f^*_{\zeta_2} \, {\rest} \, B_\xi$. 
    \end{enumerate}

    Assume that $f \colon u \rightarrow \mu$ respects $E_\xi$, then we can find
    non-empty $A_\xi \in {\mathscr  A}$ for $\xi < \kappa$ such that Rang$(f) =
    \dbcu_{\xi < \kappa} A_\xi$.  Now for any pair $(\xi_1,\xi_2) \in \kappa \times \kappa$, let $u_{\xi_1} = \text{ Min}\{\zeta \in u \colon f(\zeta) \in A_{\xi_1}\}$, let
    $\langle \zeta_i \colon i < \theta \rangle$ list $u_{\xi_2}$, let $x_{\xi_1,\xi_2}$ be $(A_{\xi_2},B_{\xi_1}  \cup B_{\xi_2}, \langle
    f^*_{\zeta_i} \, {\rest} \, B_{\xi_1} \colon i < \theta \rangle,\langle f(\zeta_i) \colon i < \theta \rangle)$.  Now check.
    
    2) Naturally the proof splits into five cases, but first, let ${\mathscr  A}$
    be as in the proof of part (1), pr$(...,...)$ be a pairing function on
    $\mu$, let $\langle f^*_\zeta \colon \zeta < \lambda \rangle$ be a sequence
    of functions from $\mu$ to $\mu$ (chosen separately in each case). Let $\mathscr{X}$ the set of elements $x = (A, \bar{g}, \bar{\beta)}$ such that:

    \begin{enumerate}
        \item[(i)] $A \in \mathscr{A},$

        \item[(ii)] $\bar{g} = \langle g_i \colon i < \theta \rangle$,

        \item[(iii)] $\text{ each } g_i \text{ is a partial function from } A \text{ to } A$ and,

        \item[(iv)] $\bar \beta = \langle \beta_i \colon i \le \theta \rangle \text{ and } i \le \theta \Rightarrow \beta_i \in A.$
    \end{enumerate}
    
    Clearly $|{\mathscr  X}| \le \mu$ and for each $x = (A,\bar g,\bar \beta)
    \in {\mathscr  X}$ we define $F_x \colon \lambda \rightarrow \mu$ by $$ F_x(\zeta) = \beta_i \text{ \underline{iff} } i = \mathbf i_x(\zeta) = \min\{\varepsilon \le \theta \colon \text{ if } \varepsilon < \theta 
    \text{ then } g_\varepsilon \subseteq f^*_\zeta\}. $$
     
    So let $u \in [\lambda]^{\le \theta}$ and $f \colon u \rightarrow \mu$ be
    given and let $\langle v_i \colon i < \kappa \rangle$ be such that $u =
    \bigcup\{v_i \colon i < \kappa\}$ and $\langle A_i \colon i < \kappa \rangle$ be 
    such that $i < \kappa \Rightarrow A_i \in {\mathscr  A}$ and 
    Rang$(f^{v_i}) \subseteq A_i$.
    
    \underline{Case 1}:  Clause $(d)_1$ holds.
    
    We choose $f^*_\zeta = f_\zeta$ for every $\zeta < \lambda$ where
    $f_\zeta$ is from clause $(d)_1$.  Let
    $\langle (u_i,B_i) \colon i < \kappa \rangle$ be as guaranteed to exist for
    $u$ in clause $(d)_1$.  For each $i,j < \kappa$ such that $u_i \cap
    v_j \ne \emptyset$ and let $\langle
    \zeta^{i,j}_\varepsilon:\varepsilon < \theta \rangle$ list $u_i \cap
    v_j$ (possibly with repetition) and let $x_{i,j} = (A^*_{i,j},\langle
    g^{i,j}_\varepsilon \colon \varepsilon < \theta \rangle,\langle
    \beta^{i,j}_\varepsilon \colon \varepsilon < \theta \rangle)$ where: 

    \begin{enumerate}
        \item[(a)] $A^*_{i,j} = B_i \cup A_j,$

        \item[(b)]  $g^{i,j}_\varepsilon = f^*_{\zeta^{i,j}_\varepsilon} \, {\rest} \, B_i,$

        \item[(c)] $\beta^{i,j}_\varepsilon = f(\zeta^{i,j}_\varepsilon)$.
    \end{enumerate}

    We can check that 
    
    \begin{enumerate}
        \item[$(\alpha)$]  $x_{i,j} \in {\mathscr  X}$ and,

        \item[$(\beta)$]  $\bigcup\{u_i \cap v_j \colon i,j < \kappa\} = u,$

        \item[$(\gamma)$] $F_{x_{i,j}} \supseteq f \, {\rest} \, (u_i \cap v_j)$.
    \end{enumerate}  
     
    [Why?  Let $\zeta(*) \in u_i \cap v_j$ so $\zeta(*) = \zeta^{i,j}_{\varepsilon(*)}$ for some $\varepsilon(*) < \theta$ so
    $g^{i,j}_{\varepsilon(*)} = f^*_{\zeta^{i, j}_{\varepsilon(*)}}
    \, {\rest} \, B_i$ by clause (b) above hence $\varepsilon(*) \in
    \{\varepsilon \subseteq \theta$: if $\varepsilon < \theta$ then
    $g_\varepsilon \subseteq f^*_{\zeta(*)}\}$.  But $\langle f^*_\zeta
    \, {\rest} \, B_i \colon \zeta \in u_i \rangle$ is a sequence without
    repetitions, so  $$\varepsilon < \varepsilon(*) \Rightarrow \zeta^{i,j}_\varepsilon \in
    u_i \cap \zeta \Rightarrow g^{i,j}_\varepsilon  =
    f^*_{\zeta^{i,j}_\varepsilon} \, {\rest} \, B_i \\  \ne f^*_{\zeta^{i,j}_{\varepsilon(*)}} = g^{i,j}_{\varepsilon(*)}$$
    
    (the middle inequality holds by the choice of $\langle (u_i,B_i) \colon i <
    \kappa \rangle$ above, i.e., by clause $(d)_1$).
    
    So $\mathbf i_{x_{i,j}}(\zeta(*)) = \varepsilon(*)$ so
    $F_{x_{i,j}}(\zeta(*)) = \beta^{i,j}_{\varepsilon(*)} =
    f(\zeta^{i,j}_{\varepsilon(*)})$, the last equality holds by clause
    (c).  So we are done.
    
    \underline{Case 2}:  Clause $(d)_2$ holds.  
    
    Follows from Case 3.
    
    \underline{Case 3}: $(d)_3$ holds.
    
    Let $f'_\zeta$ be any function from $\mu$ to $\mu$ extending
    $f_\zeta$.  It suffices to show that $\langle f'_\zeta \colon \zeta < \lambda
    \rangle$ satisfies clause $(d)_1$.  So let $u \in  [\lambda]^{\le \theta}$, so by $(d)_3$  we can find a sequence $\bar a = \langle a_\zeta \colon \zeta \in u \rangle$ of subsets of $\kappa$ such that
    $a_\alpha \notin J$ and 

    \begin{enumerate}
        \item[$\circledast$]    $\zeta_1 \ne \zeta_2 \in u \, \wedge \,  i \in a_{\zeta_1} \cap a_{\zeta_2} \Rightarrow f_\alpha(i) \ne f_\beta(i)$.
    \end{enumerate}

    Let $B_i = \{i\} \in [\mu]^{< \aleph_0}$ and $u_i = \{\zeta \in u \colon i
    \in a_\zeta\}$, clearly they are as required.
    
    \underline{Case 4}:  Clause $(d)_4$ holds.
    
    This implies Clause $(d)_3$ by \cite[II,1.3,p.46+1.4]{Sh:g}(3),p.50 (or see \cite{Sh:506}).
    
    \underline{Case 5}:  $(d)_5$.
    
    Let $\mu_1 = \min \{\partial \colon \partial^\kappa > \lambda\}$, so clearly $\mu_1 \le \mu$.  As $\kappa \le \theta,2^\theta \le \mu$ we
    have $2^\kappa < \mu_1$, hence $\kappa_1 = \cf(\mu_1)$ is $\le \kappa$ and $(\forall \alpha < \mu_1)(|\alpha|^\kappa < \mu_1)$.  By $(d)_5$ we have $\kappa_1 > \aleph_0$.  So by \cite[VIII]{Sh:g}, pp$_{J^{\bd}_\kappa}(\mu_1) \ge (\mu_1)^\kappa  = \mu^\kappa > \lambda$ so Clause $(d)_4$ holds.   
\end{PROOF}

\begin{claim}\label{x.23.9} 
    1) Assume
    
    \begin{enumerate}
        \item[(a)] $\mu = \,{\text{\rm cov\/}}(\mu,\theta^+,\theta^+,\kappa^+)$ so $\mu > \theta \ge \kappa,$
    
        \item[(b)] $2^\theta \le \mu$, 
    
        \item[(c)] $\mu < \lambda = \,{\text{\rm cf\/}}(\lambda) \le \mu^\kappa,$
    
        \item[(d)] $\theta = \,{\text{\rm cf\/}}(\theta),S \subseteq S^\lambda_\theta$ is stationary and $S \in \hat{I}[\lambda],$
    
        \item[(e)] $T \subseteq \theta$ and $[\alpha < \lambda \, \wedge \,  \delta \in T \Rightarrow \trp_{\cf(\delta)}(\vert \alpha \vert)  < \lambda]$.
    \end{enumerate}
    
    \underline{then} there is $\bar F$ such that: 
    
    \begin{enumerate}
        \item[$(\alpha)$]  $\bar F = \langle F_\varepsilon \colon \varepsilon < \mu \rangle,$
    
        \item[$(\beta)$] $F_\varepsilon$ is a function from $\lambda$ to $\mu,$
    
        \item[$(\gamma)$]  for a club of $\delta \in S$, there is an increasing continuous sequence $\langle \alpha_\varepsilon \colon \varepsilon < \theta \rangle$ with limit $\delta$ such that:
    
        \begin{enumerate}
            \item[$(\ast)$]   for every function $f$ from $\{\alpha_\varepsilon \colon \varepsilon \in T\}$ to $\mu$ we can find a sequence $\langle u_i \colon i < \kappa \rangle$ of subsets of $T$ with union $T$ such that $(\forall i < \kappa)(\exists \varepsilon < \mu)(f \, {\rest} \, \{\alpha_j \colon j \in u_i\} \subseteq F_\varepsilon)$. 
        \end{enumerate}
    \end{enumerate}
\end{claim} 

\begin{PROOF}{\ref{x.23.9}}
    1) Let $\bar a = \langle a_\alpha \colon \alpha < \lambda \rangle$ and $E_0$ as in Definition \ref{0.5}(2) witness $S \in \hat{I}[\lambda]$. Let $\bar M = \langle M_\alpha \colon \alpha < \lambda \rangle$ be an increasing continuous sequence of elementary sub-models
    of $({\mathscr  H}(\chi), \in, <^*_\chi)$, each of cardinality $<
    \lambda$ satisfying  $\bar M \, {\rest} \, (\alpha +1) \in M_{\alpha +1}, \|M_0\| = \mu, M_{1+\alpha} \cap \lambda \in \lambda$ and $\{\lambda,\theta, T, S, \bar a, E_0\} \cup (\mu+1) \subseteq M_0$.
    Let $E = \{\delta < \lambda \colon M_\delta \cap \lambda = \delta\}$, clearly
    a club of $\lambda$ and $E \subseteq E_0$ as $E_0 \in M_0$.  
    
    Now we choose by induction on $\alpha < \lambda$ a function $f_\alpha$
    from $\kappa$ to $\mu$ such that:

    \begin{enumerate}
        \item[$(\ast)_{1}$]  if $b \in 
        [\alpha]^{\le \theta} \, \wedge \,  b \in M_{\alpha +1}$ then $\ran(f_\alpha) \nsubseteq \cup\{\ran(f_\beta) \colon \beta \in b\}.$ 
    \end{enumerate}

    \begin{enumerate}
        \item[$(\ast)_{2}$] $f_\alpha$ is the $<^*_\chi$-first object satisfying $(*)_1$. 
    \end{enumerate}
    
    This is clearly possible as $\kappa \le 2^\theta \le \mu < \lambda \le
    \mu^\kappa$. \newline
    For every $i < \kappa$ and $u \in M_0 \cap [\mu]^{\le \theta}$ and 
    $g \colon u \rightarrow \mu$ from $M_0$, let $F_{g,i}$ be the unique
    $F_{g,i} \colon \lambda \rightarrow \mu$ and it satisfies: if $f_\alpha(i)
    \in \dom(g)$ then $F_{g,i}(\alpha) = g(f_\alpha(i))$. [Why
    $F_{g,i} \in M_0$?  As $\mu+1 \subseteq M_0$ and $u,g \in M_0$.]
    
    Let $\langle (g_\varepsilon, i_\varepsilon) \colon \varepsilon < \mu \rangle$
    list the set of all such pairs $(g,i)$ and let $F_\varepsilon =
    F_{g_\varepsilon,i_\varepsilon}$.  We shall show that $\langle
    F_\varepsilon \colon \varepsilon< \mu \rangle$ is as required; so clauses $(\alpha)$
    + $(\beta)$ hold trivially.
    
    Now if $\delta \in \nacc(E) \cap S \subseteq S \cap E_0$ so
    $a_\delta \subseteq \delta = \sup(a_\delta)$, $\otp(a_\delta) =
    \cf(\delta) = \theta$ and let $\langle
    \alpha_\varepsilon:\varepsilon < \theta \rangle$ list the closure in
    $\delta$ of $a_\delta$ in the increasing order.  So  for some club $C$ of $\theta$ consisting of limit ordinals 
    we have $\varepsilon \in C
    \Rightarrow \alpha_\varepsilon \in E$.  
    Let $\Theta = \cf(\delta) \colon \delta \in T\}$, now suppose $\varepsilon \in C$ is a limit ordinal with cofinality $\in \Theta$. So for every $\zeta < \varepsilon, \alpha_{\zeta +1} \in
    a_{\alpha_\varepsilon}$ so, as $\bar a \in M_0 \prec M_{\alpha_\varepsilon} \, \wedge \,  \alpha_{\zeta +1} \in
    M_{\alpha_\varepsilon}$ clearly $a_{\alpha_{\zeta +1}} \in
    M_{\alpha_\varepsilon}$ but (see Definition \ref{0.5}(2))
    $c \ell(a_{\alpha_{\zeta+1}}) = \{\alpha_\varepsilon:\varepsilon 
    \le \zeta\}$ hence
    $\{\alpha_\varepsilon \colon \varepsilon \le \zeta\} \in 
    M_{\alpha_\varepsilon}$.  Now the tree
    $(\{\nu \colon \nu \in M_{\alpha_\varepsilon}$ is a sequence of ordinals of length
    $< \delta\},\trianglelefteq)$ is a tree which belongs to $M_{\alpha_\varepsilon +1}$, has $\varepsilon$ levels 
    and $\|M_{\alpha_\varepsilon}\|
    < \lambda$ nodes so by assumption (e) it has $< \lambda$ $\varepsilon$-branches.  As this set (of $\varepsilon$-branches) belongs to $M_{\alpha_\varepsilon +1}$, every $\varepsilon$-branch of this tree belongs to $M_{\alpha_\varepsilon +1}$ but $\langle \alpha_\zeta \colon \zeta <
    \varepsilon \rangle$ is (essentially) such a branch; so together
    $\varepsilon \in C \wedge
    \cf(\varepsilon) \in \Theta \Rightarrow \langle \alpha_\zeta \colon \zeta < \varepsilon \rangle \in M_{\alpha_\varepsilon +1}$.

    Let $$ v_i = \{\varepsilon \in C:\& \text{cf}(\varepsilon) \in \Theta \text{ and} \\
      \& f_{\alpha_\varepsilon}(i) \notin \cup\{\text{Rang}(f_\beta) \colon \beta \in
    \{\alpha_\zeta \colon \zeta < \varepsilon\}\}.$$ 

    So $\langle v_i \colon i < \kappa \rangle$ is a sequence of subsets of $C$
    with union $\supseteq T$.  For each $i < \kappa$ by assumption (a)
    there is $\langle v_{i,j} \colon j < \kappa \rangle$ be a sequence of subsets
    of $v_i$ with union $v_i$ such that for every $j < \kappa
    \{f_{\alpha_\varepsilon}(i) \colon \varepsilon \in v_{i,j}\} \subseteq
    v^+_{i,j} \in M_0$.

    Let $h$ be a function from $\{\alpha_\varepsilon \colon \varepsilon \in T\}$  to $\mu$.  So by assumption (a) we  can find a sequence $\langle w_j \colon j < \kappa \rangle$ of subsets of $\mu$ each of cardinality $\le \theta$ such that
    Rang$(h) \subseteq \dbcu_{j < \kappa} w_j$ and $\{w_j \colon j < \kappa\}
    \subseteq M_0$.  Let $u_{i,j_1,j_2} = 
    \{\varepsilon \colon \varepsilon \in v_i,f_{\alpha_\varepsilon}(i) \in v^+_{i,j_1}$
    and $h(\alpha_\varepsilon) \in w_{j_2}\}$.
    Hence for $i,j_1,j_2 < \kappa$ there is $g = g_{i,j_1,j_2} \colon 
    \{f_{\alpha_\varepsilon}(i) \colon \varepsilon \in u_{i,j_1,j_2}\}
    \rightarrow \mu$ satisfying $\varepsilon \in u_{i,j_1,j_2} \Rightarrow
    g(f_{\alpha_\varepsilon}(i)) = h(\alpha_\varepsilon)$, see the choice
    of $v_i,v_{i,j_1},w_{j_2},u_{i,j_1,j_2}$.  
    Now as $v_i \subseteq \theta,2^\theta \le \mu,\mu +1
    \subseteq M_0$ and $v^+_{i,j_1},w_{j_2} \in M_0$ are subsets of $\mu$
    of cardinality $\le \theta$ clearly the set of partial functions from 
    $v^+_{i,j_1}$ to $w_{j_2}$ belongs to $M_0$ and has cardinality
    $\le^{|u_{i,j}|}|w_j| \le \theta^\theta = 2^\theta \le \mu$
    hence is included in $M_0$.  So $g_{i,j}$ belongs to $M_0$ for every
    $i,j < \kappa$.  Also easily $F_{g_{i,j,i}}$ extends  $h \, {\rest} \, (\{\alpha_\varepsilon \colon \varepsilon \in u_{i,j}\})$.  So
    $\langle \{\alpha_\varepsilon \colon \varepsilon \in u_{i,j}\} \colon 
    i,j < \kappa \rangle$ and $\langle
    \alpha_\varepsilon:\varepsilon \in C \rangle$ are as required 
    (modulo some renaming).   But $\dom(g_{i,j}) = \{f_{\alpha_\varepsilon}(i) \colon \varepsilon \in u_{i,j}\} \subseteq w_j$.
\end{PROOF}

Compare the following with \cite{Sh:237e}.

\begin{claim}\label{x.2.9} 
    1) Assume

    \begin{enumerate}
        \item[(a)]   $\lambda = \mu^+$ and $\kappa = { \text{\rm cf\/}}(\mu) < \mu$ and $\aleph_0 < \theta = \,{\text{\rm cf\/}}(\theta) < \mu,$

        \item[(b)] 

        \begin{enumerate}
            \item[(i)] $\Theta  = \{\sigma \colon \aleph_0 \le \sigma = { \text{\rm cf\/}}(\sigma) < \mu$ and $ \trp_{\sigma}(\mu) = \mu\},$

            \item[(ii)]  $ T_\theta =  \{\delta < \theta \colon {\text{\rm cf\/}}(\delta) \in \Theta\},$

            \item[(iii)] $T^+_\theta =  T_\theta \cup \{0\} \cup \{\zeta +1 \colon \zeta < \theta\}$,

            \item[(iv)] $T^-_\theta = T^+_\theta \backslash T_\theta.$
        \end{enumerate}

        \item[(c)] $\mathbf P \subseteq \{\bar A \colon \bar A$ a sequence  of length $\le \mu$ of subsets of $\theta$ each included in $T^+_\theta$ and including $T^-_\theta\},$

        \item[(d)]  ${\mathscr  P} \subseteq [\mu]^{\le \theta}$ has cardinality $\le \mu,$ 

        \item[(e)]  if $\alpha_\varepsilon < \mu$ for $\varepsilon < \theta$ and $\langle w^*_i \colon i < \kappa \rangle$ a sequence of subsets of $T^+_\theta$ and including $T^-_\theta,$     \underline{then} for some $\bar A \in \mathbf P$ we have: 

        \begin{enumerate}
            \item[(i)] for some $\zeta < \ell g(\bar A)$ and $j < \kappa,A_\zeta \subseteq w^*_j$,
    
            \item[(ii)]  $\zeta < \ell g(\bar A) \, \, \wedge \, \,  j < \kappa \, \wedge \,  A_\zeta \subseteq w^*_j \Rightarrow \{\alpha_\varepsilon \colon \varepsilon \in A_\zeta\} \in {\mathscr  P}$.
        \end{enumerate}
    \end{enumerate}

    \underline{Then} we can find $(\bar C_i, A^*_i,S_i)$ for $i < \mu$ such that:

    \begin{enumerate}
        \item[$(\alpha)$] $\{0\} \cup \{\zeta +1 \colon \zeta < \theta\} \subseteq A^*_i \subseteq \theta,$

        \item[$(\beta)$]  $\bar C_i = \langle C^i_\delta \colon \delta \in S_i \rangle; S_i \subseteq \lambda$ and $$[{\text{\rm otp\/}}(\alpha \cap C^i_\delta)  \in A^*_i \, \wedge \,  C^i_\delta \subseteq \delta \, \, \wedge \,  \,  \alpha \in  C^i_\delta \Rightarrow \alpha \in S_i \, \wedge \,  C^i_\alpha = \alpha \cap C^i_\delta],$$

        \item[$(\gamma)$] $S^*_i = \{\delta \in S_i \colon {\text{\rm cf\/}}(\delta)  = \theta\}$ is stationary,

        \item[$(\delta)$]  for every club $E$ of $\lambda$ there are $\delta^* \in E$ and a sequence $\langle w^*_i \colon i < \kappa \rangle$ of subsets of $T^+_\theta$ including $T^-_\theta$ with union $T^+_\theta$ such that  for each $j < \kappa$ there is $\bar A \in \mathbf P$ satisfying:  $(\forall \zeta < \ell g(\bar A))(A_\zeta \subseteq w^*_j \rightarrow (\exists i < \mu)[\delta^* \in S^*_i \, \wedge \,   A^*_i = A_\zeta])$

        \item[$(\varp)$] if $B \in \hat{I}[\lambda]$ is stationary and $B \subseteq S^\lambda_\theta$ then we can demand $B \backslash \cup \{S^*_i \colon I < \mu\}$  is not stationary and in clause $(\delta)$ any $\delta \in B$ satisfies the requirement on $\delta^*$. 
    \end{enumerate}
    
    2) Assume, in addition $\alpha < \mu \Rightarrow |\alpha|^\theta < \mu,
    T'_\theta = \{\delta < \theta \colon {\text{\rm cf\/}}(\delta) \ne \kappa\}$ and we let $\mathbf P = \{\bar A \colon \bar A = \langle A_i \colon i < \kappa \rangle$ is a
    sequence of sets with union $T^+_\theta\}$.  
    \underline{Then} $T'_\theta = T^+_\theta$ and
     for some ${\mathscr  P}$
    all the assumptions hence the conclusion of part (1) holds.
    
    3) Assume $\mu_*$ is strong limit $< \mu,\kappa = \,{\text{\rm cf\/}}(\mu) <
    \mu,\lambda = \mu^+$, \underline{then} for some $\theta_0 < \mu_*$, for any
    regular $\theta \in (\theta_0,\mu_*)$ the set $T_\theta$ include $\{\delta <
    \theta \colon {\text{\rm cf\/}}(\delta) \ge \theta_1$ or $\delta$ non-limit$\}$ and
    for $\mathbf P = \{\bar A \colon \bar A$ is a sequence of $< \theta_0$ subsets
    of $\theta$ with union $T_\theta$ such that $\{A_\zeta \colon \zeta < \ell
    g(\bar A)\}$ is closed under finite unions$\}$. 
    
    \underline{Then} for some ${\mathscr  P}$ the assumptions hence the conclusion of
    part (1) holds. 
    
    4) In (1) (and (2), (3), if $\bar A = \langle
    A_\varepsilon \colon \varepsilon < \varepsilon(*) \rangle \in \mathbf P$ we can
    replace it by $\bar A' = \langle A'_\varepsilon \colon \varepsilon <
    \varepsilon(*) \rangle,A'_\varepsilon = A_\varepsilon \cup \{\delta <
    \theta \colon {\text{\rm cf\/}}(\delta) > \aleph_0$ and $A_\varepsilon \cap \delta$
    is stationary in $\delta\}$ and still obtain an element of $\mathbf P$.
\end{claim} 

\begin{remark}\label{x.2.9A} 
    1) If  $\lambda$ is (weakly) inaccessible we have weaker results: letting $\Theta = \{\sigma \colon \sigma = \cf(\sigma)$ and $\alpha < \lambda \Rightarrow \trp_{\sigma}(\vert \alpha \vert) < \lambda\}$,  there is $\langle {\mathscr  P}_\alpha \colon \alpha < \lambda \rangle,
    {\mathscr  P}_\alpha \subseteq [\alpha]^{< \theta},|{\mathscr  P}_\alpha| < \lambda$ such that for regular $\theta < \lambda$ for every stationarily many $\delta \in S^\lambda_\theta$ for some increasing continuous sequence $\langle
    \alpha_\varepsilon \colon \varepsilon < \alpha \rangle$ with limit $\delta,[\varepsilon < \theta \, \wedge \,  \cf(\varepsilon) \in \Theta_\lambda \Rightarrow \{\alpha_\zeta \colon \zeta < \varepsilon\} \in {\mathscr  P}_{\alpha_\varepsilon}]$ (see \cite[\S1]{Sh:420}). 
    
    2) Can we restrict ourselves to: $\langle w^*_j \colon j < \kappa \rangle$ is
    increasing continuous in part (1) and similarly in (2).
    
    It seems very reasonable (if $\lambda = \mu^+,\kappa = \cf(\mu)
    < \mu$) that we can find $\langle \lambda_i \colon i < \kappa \rangle$, an
    increasing sequence of regular  cardinals $< \mu$ with limit $\mu$ such that $\lambda^+ = \mathrm{tfc}(\, \prod_{i < \ell}
    \lambda_i,\sigma \le_{J^{\bd}_\kappa})$ and $j < \kappa
    \Rightarrow \lambda_j > \text{ max pcf}\{\lambda_i \colon i < j\}$.  Why? This is trivial if $\kappa = \aleph_0$ and if $\kappa > \aleph_0$, it follows by the weak hypothesis (see \cite{Sh:410}), whose negation is not known to be consistent and which has large consistency strength.
    
    3) If there is $\langle \lambda_i \colon i < \kappa \rangle$ is as above then in part (1) of Claim \ref{x.2.9} we can strengthen clauses $(\delta),(\varepsilon)$ in the conclusion.  By demanding that ``$\langle w^*_j \colon j < \kappa \rangle$" is increasing continuously.  Why? Choosing $\langle f_\alpha \colon \alpha < \lambda \rangle$ we add:
    
    \begin{enumerate}  
        \item[$\boxdot$]  if $f_{\alpha_1}(i_1) = f_{\alpha_2}(i_2)$ then $i_1 = i_2 \, \wedge \,  f_{\alpha_1} \, {\rest} \, i_1 = f_{\alpha_2} \, {\rest} \, t_2$. 
    \end{enumerate}  
    
    So there is a sequence $\langle \mathbf g_i \colon i < \kappa \rangle,\mathbf
    g_i \colon \mu \rightarrow \mu$ and $\mathbf g_{i_1}(f_{\alpha_1}(i_1)) =
    f_{\alpha_1}(i_0)$ when $i_0 < i_1 < \kappa,\alpha < \lambda$ and a
    function $\mathbf I \colon \mu \rightarrow \kappa$ such that $\mathbf i(f_\alpha(i)) = i$.
            
    Now at the end of the proof we replace $\alpha'_\varepsilon = \cd(j,f_{\alpha_\varepsilon}(j),\xi(\alpha_\varepsilon))$ (so fixing
    $j$) by $\alpha'_\varepsilon = \cd(\xi)$.
\end{remark}

\begin{PROOF}{\ref{x.2.9}}
    There is $B$ as in clause $(\varepsilon)$, that is 
    $B \in \hat{I}[\lambda]$ and $B \subseteq S^\lambda_\theta$ is a 
    stationary subset of $\lambda$, see Claim \ref{d.13}(2).  So let $B$ be
    such a set; clearly omitting from $B$ a non-stationary subset of
    $\lambda$ is O.K. as if $B' \subseteq B,B \backslash B'$ is
    non-stationary and the demands in clause $(\varepsilon)$ for $B'$ are 
    exemplified by the sequence $\langle \bar C_i,A^*_i,S^*_i) \colon i < \mu
    \rangle$ then the same sequence exemplifies it for $B$.
    
    We can find $\bar a = \langle a_\alpha \colon \alpha < \lambda \rangle$
    witnessing $B \in \hat{I}[\lambda]$ see Definition \ref{0.5}(2), so 
    (possibly omitting from $B$ a
    non-stationary set) we have: $\alpha \in B \Rightarrow a_\alpha
    \subseteq \alpha =
    \sup(a_\alpha) \, \wedge \,  \theta = \otp(a_\alpha)$ and $[\beta \in a_\alpha
    \Rightarrow a_\beta = \beta \cap a_\alpha]$ and $[\alpha \in \lambda
    \backslash B \Rightarrow a_\alpha \subseteq \alpha \, \wedge \, 
    \otp(a_\alpha) < \theta]$.
    
    Let $\chi = \lambda^+$ and $\bar M = \langle M_\alpha \colon \alpha < \lambda
    \rangle$ be an increasing continuous sequence of elementary sub-models
    of $({\mathscr  H}(\chi),\in,<^*_\chi)$ to each of which $\{\bar a,\mu,\mathbf P,
    {\mathscr  P}\}$ belongs, and satisfying $\bar M \, {\rest} \, (\alpha + 1) \in
    M_{\alpha +1},\mu \subseteq M_\alpha$ and $\|M_\alpha\| < \lambda$. 

    As $\lambda = \mu^+$ clearly $\alpha < \lambda \Rightarrow
    \|M_\alpha\| = \mu$.  
    For any $\alpha < \lambda$ let $\langle b_{\alpha,\zeta} \colon \zeta < \mu
    \rangle$ list $M_{\alpha +1} \cap [\alpha]^{< \theta}$.
    
    Let $\langle \lambda_i \colon i < \kappa \rangle$ be an increasing sequence of regular cardinals with limit $\mu, \lambda_i > \theta$ with $\lambda = \tcf(\,  \prod_{i < \kappa} \lambda_i,<_{J^{\bd}_\kappa})$
    and let $\bar f = \langle f_\alpha \colon \alpha < \lambda \rangle$ be an
    increasing sequence in $( \, \prod_{i < \kappa} \lambda_i,<_{J^{\bd}_\kappa})$ cofinal in it and obeying $\bar a$ (i.e. $i < \kappa \, \wedge \,  \beta \in a_\alpha
    \Rightarrow f_\beta(i) < f_\alpha(i))$ and $\zeta < \lambda_i \wedge i
    < j < \kappa \wedge \beta \in b_{\alpha,\zeta} \Rightarrow f_\beta(j)
    < f_\alpha(j)$ (there is such $\langle \lambda_i \colon I < \kappa \rangle$ by \cite[II]{Sh:g}); there is such $\bar f$ by \cite[I]{Sh:g}. For notational help later \wilog \, each $f_\alpha(i)$ is divisible by $\theta$.
    
    Now for any $x = (h_1,h_2,A,i) = (h^x_1, h^x_2, A^x, i^x)$ where $h_1$ is a partial function from $\mu$ to $\theta,h_2$ a function from $A$ to $\mu$ where $T^-_\theta \subseteq A \subseteq T^+_\theta$ and $i < \kappa$, we define:
    
    \begin{enumerate} 
        \item[$(*)$] $S_x$ is the set of $ \alpha $ such that: 
        
        \begin{enumerate}
            \item[(a)] $ \alpha < \lambda,$
            
            \item[(b)] $ \cf(\alpha) \ne \theta,$ 
            
            \item[(c)] $ f_\alpha(i) \in  \Dom (h_1),$
            
            \item[(c)] $ \zeta \coloneqq h_1(f_\alpha(i)) \in A$  satisfies:  
        
            \begin{enumerate} 
                \item[(i)] $ b_{\alpha,h_2(\zeta)}$  
                is a subset of $ \alpha $  of order type $ \zeta,$
                
                \item[ii)]  $b_{\alpha,h_2(\zeta)}$ 
                is unbounded in $ \alpha$  if $  \zeta $  is a limit ordinal and,
                
                \item[(iii)] if $  \varepsilon \in A \cap \zeta $  $ \beta \in b_{\alpha, h_2(\zeta)} $  $ \varepsilon = \otp(b_{\alpha, h_2(\zeta)} \cap \beta),$  
            \end{enumerate}  
            
            \underline{then}  $\varepsilon = h_1(f_\beta(i))$  and  $ \beta \cap b_{\alpha, h_2(\zeta)} = b_{\beta, h_2(\varepsilon)}$ and for $\alpha \in S_x$ let $C^x_\alpha = b_{\alpha, h^*_2(h^*_1(f_\alpha(i))}$. 
        \end{enumerate} 
    \end{enumerate}  
    
    For $x$ as above let $S^*_x = \{\delta < \lambda \colon \cf(\delta) =
    \theta$ and there is an increasing continuous sequence $\langle
    \alpha_\varepsilon \colon \varepsilon < \theta \rangle$ with limit $\delta$
    such that $\varepsilon \in A \Rightarrow \alpha_\varepsilon \in S_x
    \, \wedge \,  C^x_{\alpha_\varepsilon} = \{\alpha_\zeta \colon \zeta < \theta\}\}$.
    Now for any such $x$, for $\delta \in S^*_x$ choose $\langle
    \alpha_\varepsilon \colon \varepsilon < \theta \rangle$ witnessing it and let
    $C^x_\delta = \bigcup\{C^x_{\alpha_\varepsilon} \colon \varepsilon \in A\}$;
    note that $S_x \cap S^*_x = \emptyset$ because $[\alpha \in S_x
    \Rightarrow \cf(\alpha) \ne \theta],[\alpha \in S^*_x
    \Rightarrow \cf(\alpha) = \theta]$.
    
    Let $\langle x_\varepsilon \colon \varepsilon < \varepsilon^* \le \mu
    \rangle$ list the set of $x \in M_0$ which are as above and $S^*_x$
    is stationary, and for each $\varepsilon < \varepsilon^*$, let $S_\varepsilon \coloneqq  S_{x_\varepsilon} \cup S^*_{x_\varepsilon}$,  $\bar C_\varepsilon \coloneqq  \langle C^{x_\varepsilon}_\alpha \colon \alpha \in S_\varepsilon \rangle,A_\varepsilon \coloneqq A^{x_\varepsilon}$.  We have to check the five clauses of Claim \ref{x.2.9}(1):
    
    \underline{Clause $(\alpha)$}:  See the demand on $x$ above.
    
    \underline{Clause $(\beta)$}:  Holds by clause (iii) in the definition of
    ``$\alpha \in S_x"$ and of $C^x_\alpha$ for $\alpha \in S_x$ and similarly for $\alpha \in S^*_x$.

    \underline{Clause $(\gamma)$}:  Holds by the choice of $\langle x_\varepsilon \colon \varepsilon < \varepsilon^* \rangle$.

    \underline{Clause $(\delta) + (\varepsilon)$}:  Let $\delta \in B$, and we shall show that $\delta \in \cup\{S^*_i \colon i < \mu\}$, this clearly
    suffices for $(\varepsilon)$ hence for $(\delta)$.  Let
    $\langle \alpha_\zeta \colon \zeta \le \theta \rangle$ list the closure of
    $a_\delta$ in the increasing order by the choice of $\bar a, a_\delta$
    is a set of non-limit ordinals hence cf$(\alpha_\zeta) = \theta
    \Leftrightarrow \zeta = \theta$ and $\alpha_\theta = \delta$.
    
    Now for every $\zeta \in T^+_\theta$, for some $\xi(\zeta) < \mu$ we
    have $\{\alpha_\xi \colon\xi < \zeta\} = b_{\alpha,\xi(\zeta)}$ 
    ($\xi(\zeta)$ exists as $\trp_{\cf(\zeta)}(\mu) =  \mu$
    by the definition of $T_\theta$) and let
    $j(\zeta) = \min\{j < \kappa \colon f_{\alpha_\zeta}(j) >
    f_{\alpha_\varepsilon}(i)$ for every $\varepsilon < j\}$.
    
    For each $j < \kappa$ let $w^*_j = \{\zeta \colon \zeta \in T^-_\theta$ or $\zeta \in T_\theta \, \wedge \,  j(\zeta) =j\}$ and apply clause $(e)$ of the
    assumptions to the sequence $\langle \alpha'_\zeta \colon \zeta < \theta
    \rangle$ where $\alpha'_\zeta = \cd (j(\zeta),f_{\alpha_\zeta}(j_\zeta),\xi(\alpha_\zeta))$, $\cd$ a coding function on $\mu$ so (by (e)) we can find $\bar A \in \mathbf P$
    satisfying clauses (i) and (ii) of assumption (e).  So by (i) of (e)
    there are $i(*) < \ell g(\bar A)$ and $j(*) < \kappa$ such that
    $A_{\zeta(*)} \subseteq w^*_{j(*)}$.
    For any pair $(i,j)$ such that $i < \ell g(\bar A), j <\kappa$ and $A_i
    \subseteq w^*_j$ we shall define $x = x_{i,j}$.  By (ii) of (e) we
    have $\{\alpha'_\varepsilon \colon \varepsilon \in A_i\} \in {\mathscr  P}$.  
    
    We now define $x= (h^x_1,h^x_2,A^x,i^x)$ as follows: we let $i^x = j, A^x = A_i$, $\dom(h^x_1) = \{f_{\alpha_\varepsilon}(j) \colon \varepsilon \in A^x\},h^x_1(f_{\alpha_\varepsilon}(j)) = \varepsilon$, $\dom(h_2) = A^x,h_2(\varepsilon) = \xi(\alpha_\zeta)$.  Easily $x = x_{i,j} \in M_0$ (as ${\mathscr  P} \in M_0,|{\mathscr  P}| \le \mu,\mu +1 \subseteq M_0$) and $\alpha_\theta \in S^*_x$, so we are done.
    
    2) Let ${\mathscr  P} = \{A:A$ is a bounded subset of $\mu$ of cardinality
    $\le \theta\}$.  Clearly $|{\mathscr  P}| \le \Sigma\{|\alpha|^\theta \colon \alpha < \mu\} \le \Sigma\{|\alpha| \colon \alpha < \mu\} = \mu$, so clause (d) holds.
    
    Now if $\sigma = \cf(\sigma) < \theta$ and $\sigma \ne \kappa$
    then $\trp_{\sigma}(\mu) = \mu$. [Why?  For completeness if
    $T$ is a tree with set of nodes $\subseteq \mu$ and $\sigma$ levels,
    then every $\sigma$-branch of $T$ includes a $<_T$-unbounded subset
    which is bounded as a subset of $\mu$ hence belongs to ${\mathscr  P}$, and
    clearly the branch is determined by any $<_T$-unbounded subset hence
    the number of $\sigma$-branches of $T$ is $\le |{\mathscr  P}| \le \mu$.]
    
    So $T'_\theta = T_\theta$ 
    and the non-obvious clause to check is (e). So assume $\alpha_\varepsilon < \mu$ for $\varepsilon < \theta$ and $\langle w^*_i \colon i < \kappa \rangle$ is a sequence of subsets of $T^+_\theta$. Let $\mu = \Sigma\{\mu_i \colon i < \kappa\},\mu_i < \mu$ and let $\bar A = \langle A_i \colon i < \kappa \rangle$ list the family $\{\{\alpha_\varepsilon \colon \varepsilon \in w^*_i$ and $\alpha_\varepsilon
    < \mu_j\} \colon I, j < \kappa\}$.
    
    3), 4)  Left to the reader.  
\end{PROOF}

\begin{claim}\label{x.24} 
    1) If $\lambda$ is regular and not
    strongly inaccessible and $\mu < \lambda$ is a strong limit singular of uncountable cofinality, \underline{then} for every large enough regular $\kappa < \mu$,  for every $\Upsilon, \theta < \mu$, some   $(\lambda, \kappa, \kappa^+)$-parameter $\bar C$ has the
    $(\Upsilon,\theta)$-{\rm \Bb}-property (see Definition \ref{d.12G}(4)); in fact, $\bar C$ is a sequence of clubs.

    2) Moreover, we can find such $\bar C^\alpha$ for $\alpha < \lambda$
    satisfying: $\langle {\text{\rm Dom\/}}(\bar C^\alpha) \colon \alpha < \lambda
    \rangle$ are pairwise disjoint.
    
    3) Moreover for every large enough regular $\kappa < \mu$ for every
    $\chi, \sigma < \mu$ some $(\lambda, \kappa, \kappa^+)$-parameter has the $(\chi, \sigma)$-{\rm \Bb}$_\ell$-parameter for $\ell = 0,1,2$.
\end{claim} 

\begin{PROOF}{\ref{x.24}}
    1) Let $\lambda_1 < \lambda$ be such that $\lambda \le 2^{\lambda_1}$ and let $\kappa_0 < \mu$ be large enough such
    that: 

    \begin{enumerate}
        \item[$(*)_1$]

        \begin{enumerate}
            \item[(i)] if $\lambda$ is not a successor of regular, \underline{then} $\alpha < \lambda \Rightarrow \mathrm{cov}(|\alpha|,\mu,\mu,\kappa_0) < \lambda,$

            \item[(ii)]  if $\lambda$ is inaccessible, then $\cov(\lambda_1,\mu,\mu,\kappa_0) = \lambda_1$ (exist by \cite{Sh:460}).
        \end{enumerate}
    \end{enumerate}
    
    Note that $\lambda_1$ is needed only if $\lambda$ is (weakly) inaccessible.  
    
    Let $\mu_* \in (\kappa_0,\mu)$ be strong limit; this is possible as
    $\mu$ is a strong limit 
    of uncountable cofinality.  Let $\kappa_1 = \kappa_1[\mu_*] <
    \mu_*$ be $> \kappa_0$ and such that $\alpha < 2^{\mu_*} \Rightarrow
    \cov(|\alpha|,\mu,\mu,\kappa_1) < 2^{\mu_*}$; (such $\kappa_1$
    exists by \cite{Sh:460}).
    
    Let $\kappa \in \text{ Reg} \cap \mu_* \backslash \kappa_1$
    and we  shall prove the result for $\kappa$ and every $\Upsilon,\theta <
    \mu_*$, so \wilog \, $\Upsilon = \theta$.  By \emph{Fodor's Lemma} on $\mu$ this is enough.   Let $\lambda_* = \cf(2^{\mu_*})$.  By Claim \ref{d.13}(1) for every stationary $S \subseteq \{\delta < \lambda_* \colon \text{cf}(\delta) = \kappa\}$ from $\hat{I}[\lambda_*]$ (and there are such, see Claim \ref{d.13}(2)) there are $\bar C^S = \langle C^S_\delta \colon \delta \in S \rangle, D$ such that $\bar C^S$ is a weakly good$^+ \, (\lambda_*,\kappa,\kappa^+)$-\Bb-parameter, $D$ is a $\lambda_*$-complete filter on $\lambda_*$ containing the clubs of $\lambda_*$ and by Claim \ref{d.6} + Claim \ref{d.7} $S \in D^+$, (e.g., $D
    = {\mathscr  D}_\lambda$), the pair
    $(\bar C^S,D)$ has the $\theta$-\Bb-property. 
    
    For the case $\lambda$ is (weakly) inaccessible we add ``each $C_\delta$ is a club of $\delta$".  This is permissible as we can use clause
    $(\beta)$ of Claim \ref{d.13}(1).  Fix such $S_1 \subseteq S^{\lambda_*}_\kappa$ and $\bar C = \langle C_\delta \colon \delta \in S_1  \rangle = \bar C^{S_1}$.
    
    Note that we actually have:

    \begin{enumerate}
        \item[$\boxtimes$] 

        \begin{enumerate}
            \item[(a)]  if $S' \subseteq S_1$ is stationary, $C'_\delta \subseteq C_\delta \subseteq \delta = \sup(C'_\delta)$  for $\delta \in S'$ then $\bar C' = \langle C'_\delta \colon \delta \in S' \rangle$ is a weakly good$^+ \,(\lambda_*, \kappa, \kappa^+)$-\Bb-parameter and $\bar C'$   has the $(D, \theta)$-\Bb-property, 

            \item[(b)] if $\lambda_* = \cup\{A_i \colon i < i(*) < \kappa\}$ then for some $i, j<i(*), \bar C'$ is a good $(\lambda_*, \kappa, \kappa^+)$-\Bb-parameter and $\bar C'$ has the $(D,\theta)$-\Bb-property where $S' \coloneqq \{\delta \in S_1 \colon \delta \in A_j \text{ and }\delta = \sup(\mathrm{nacc} (C_\delta) \cap A_i)\}$ and for $\delta \in S'$ we let $C'_\delta \coloneqq \mathrm{ nacc}(C_\delta) \cap A_j$.
        \end{enumerate}
    \end{enumerate}
    
    [Why?  Let $S'_{i,j} = \{\delta \in S_1 \colon \delta \in A_j \text{ and }
    \delta = \sup(\mathrm{nacc}(C_\delta) \cap A_i)\}$.   Clearly $\langle S_{i,j} \colon I, j < i(*) \rangle$ is a sequence of subsets of $S_1$ with union $S_1$.  By Claim \ref{d.12B}, for some
    $i, j <i(*)$ the sequence $\bar C \, {\rest} \, S_{i,j}$ has the
    $(D,\theta)$-\Bb-property.  By monotonicity, i.e., by Claim \ref{d.11}(1), the sequence $\langle \mathrm{nacc}(C_\delta) \cap A_i \colon 
    \delta \in S_{i,j} \rangle$ has the $(D,\theta)$-\Bb-property.]

    Now, we consider several cases:  
    
    \underline{Case 1}:  $\lambda$ is the successor of regulars.
    
    We apply Claim \ref{x.20}(1) with $(\lambda,\lambda_*)$ here standing for
    $(\lambda_2,\lambda_1)$ there; there is suitable partial square (in ZFC by \cite[\S4]{Sh:351}).  Note that we do not ``lose" any $\kappa < \mu$
    because of ``the lifting to" $\lambda$.
    
    \underline{Case 2}: $\lambda$ is the successor of a singular.
    
    The proof is similar only now we use Claim \ref{x.2.9}(3) to get the
    appropriate partial square using $\boxtimes$ to derive the $\bar C'$
    used as input for it (recalling that each sequence  $\bar A \in \mathbf P$ is closed under the union of two!)
    
    \underline{Case 3}:  $\lambda$ is (weakly) inaccessible. 
    
    Really the proof, in this case, covers the other cases (but $\lambda_1$
    in this case, contributes to increasing $\kappa_0$).  This time we  use Claim \ref{x.23.9} + Claim \ref{x.20.9} which is proved below.  Let $S_2$ be a stationary subset of $S^\lambda_{\lambda_*}$ from
    $\hat{I}[\lambda]$, it exists, see Claim \ref{d.13}(2). Let $\langle F_\varepsilon \colon \varepsilon <
    \lambda_1 \rangle$ be as in Claim \ref{x.23.9} with  $(\lambda, \lambda_1, \lambda_*, \kappa)$  here standing for $(\lambda, \mu, \theta, \kappa)$ there.  Let $g^*$ be a function from $\lambda_1$ onto ${\mathscr  H}_{\le \kappa}(\lambda_*)$. Possibly replacing $S_2$ by its intersection with some club of
    $\lambda$, there is a good $\lambda_*$-\Bb-parameter $\bar C^2 =
    \langle C^2_\delta \colon \delta \in S_2 \rangle$.  Let $\bar e = \langle
    e_\alpha \colon \alpha \in S^\lambda_\kappa \rangle,e_\alpha$ a club of $\alpha$ of
    order type $\kappa$.
    
    Clearly,
    
    \begin{enumerate} 
        \item[$(*)$]  for some $\varepsilon < \lambda_1$ for any club $E$ of $\lambda$ for some increasing continuous sequence $\langle \alpha_i \colon i \le \lambda_* \rangle$ of member of $E$ the set $B_\varepsilon(\langle \alpha_i \colon i \le \lambda_* \rangle) =  \{i \in S_1 \colon g^{\ast} (F_\varepsilon(\alpha_i)) = \{(i,\zeta, \mathrm{otp}(\alpha_\zeta \cap C^2_{\alpha_{\lambda_*}}), \mathrm{otp}(\alpha_\zeta \cap e_{\alpha_i})) \colon \zeta <i$ and $\alpha_\zeta \in e_{\alpha_i}\}\}$ is a stationary subset of $\lambda_*$.
    \end{enumerate} 
    
    [Why?  For every club $E$ we can choose $\alpha_i \in E \, (i \le \lambda_*)$ increasing continuous so for some $\varepsilon$ the set $B_\varepsilon(\langle \alpha_i \colon i \le \lambda_* \rangle)$ is a stationary subset of $\lambda_*$.  As the number of possible $\varepsilon$ is $\lambda_1 < \lambda = \cf(\lambda)$, there is
    $\varepsilon$ as above.]

    For this $\varepsilon$ let

    \begin{enumerate}
        \item[$(\ast)$] 

        \begin{enumerate}
            \item[(i)]  $S^*_3 = \{\alpha < \lambda \colon \text{cf}(\alpha) = \kappa$ and $g^*(F_\varepsilon(\alpha))$ has the right form$\},$

            \item[(ii)]  $\delta \in S^*_3$ let $C^3_\delta = \{\beta \in e_\delta \colon (i,\zeta,j, \mathrm{otp}(\beta \cap e_\delta)) \in  g_*(F_\varepsilon(\delta))\},$

            \item[(iii)] let $F' \colon \lambda \rightarrow \lambda_*$ be defined by $F'(\alpha) = \zeta$ if: $\zeta = \zeta_1$ for some $(\zeta_1,\zeta_2,\zeta_3,\zeta_4) \in  g^{\ast} (F_\varepsilon(\alpha))$.
        \end{enumerate}
    \end{enumerate}
    
    Finally, we apply Claim \ref{x.20.9}.
    
    2), 3)  Similarly using Claim \ref{x.21A} when using Claim \ref{x.21}.
\end{PROOF}

We now state another version of the upgrading of \midia.

\begin{claim}\label{x.20.9} 
    1) Assume

    \begin{enumerate}
        \item[(a)]  $\lambda_1 < \lambda_2$ are regular cardinals, 

        \item[(b)]  $S_1 \subseteq S^{\lambda_1}_\kappa$ is stationary, $\bar C^1 = \langle C^1_\delta \colon \delta \in S_1 \rangle$ is a $(\lambda_1,\kappa,\chi)$-{\rm \Bb}-parameter which has the $(D_{\lambda_1},\Upsilon,\theta)$-{\rm \Bb}$_\ell$-property, $D_{\lambda_1}$ is a filter on $\lambda_1$ to which $S_1$ belongs,

        \item[(c)]  $S_2 \subseteq S^{\lambda_2}_\kappa$ is stationary and $\bar C^2 = \langle C^2_\delta \colon \delta \in S_2 \rangle$ is a $(\lambda_2,\kappa,\chi)$-{\rm \Bb}-parameter,

        \item[(d)]  $F \colon \lambda_2 \rightarrow \lambda_1$ satisfies: for every club $E$ of $\lambda_2$ there is an  increasing continuous function $h$ from $\lambda_1$ into $\lambda_2$ with $\sup({\text{\rm Rang\/}}(h)) \in {\text{\rm acc\/}}(E)$  such that $\{\delta \in S_1 \colon h(\delta) \in S_2$ and $C^2_{h(\delta)} \subseteq \{h(\alpha) \colon \alpha \in C^1_\delta\}$ and $F(h(\delta)) = \delta\}$ belongs to $D_{\lambda_1},$

        \item[(e)] $2^{2^{\lambda_1}} < \lambda_2$,
    \end{enumerate}
    
    \underline{then} $\bar C^2$ has the $(\Upsilon,\theta)$-{\rm \Bb}$_\ell$-property.

    2) Assume that in part (1) we omit assumption (e) and replace the 
    {\rm \Bb}-property by {\rm MD}$_\ell$-property where $\ell \in
    \{0,2\}$. Then $\bar C^2$ has the $(\Upsilon,\theta)$-{\rm MD}$_\ell$-property. 
    
    3) In part (1) if we omit assumption (e) and strengthen clause (d) by
    demanding $\{\delta \in S_1 \colon h(\delta) \in S_2$ and $C^2_{h(0)} =
    \{h(\alpha) \colon \alpha \in C^1_\delta\}$ and $F(h(\delta)) = \delta\}$ belongs to $D_{\lambda_1}$, \underline{then} we can still conclude that: for every nice $(\bar C^2,\Upsilon,\theta)$-colouring $\mathbf F$ (see below) there is an $\mathbf F$-\midia \ sequence.
\end{claim} 

\begin{remark}\label{x.29}
    1) Is clause (d) of the assumption reasonable?  Later in \cite[4.5]{Sh:F611} we prove that it holds under reasonable conditions. 
    
    2) This + \ref{x.20.10} are formalized elsewhere.
\end{remark} 

\begin{PROOF}{\ref{x.20.9}}
    Similar to Claim \ref{x.20}.
    
    First, assume $\ell=2$ (so we are in part (2)) and  let $\tau$ be a relational vocabulary of
    cardinality $< \Upsilon$ and by clause (b) of the assumption there is
    a sequence $\bar M^1 = \langle M^1_\delta \colon \delta \in S_1
    \rangle, M^1_\delta$ a $\tau$-model with universe $C^1_\delta$ such
    that: 
    
    \begin{enumerate}  
        \item[$\circledast_1$]   if $M$ is a $\tau$-model with universe $\lambda_1$, \underline{then} the set of $\delta \in S_1$ for which $M^1_\delta = M \, {\rest} \, C^1_\delta$ is $\ne \emptyset \mod D_{\lambda_1}$.
    \end{enumerate}    
    
    Now for every $\delta \in S_2$ if $F(\delta) \in S_1$ and $\{F(\alpha) \colon \alpha \in C^2_\delta\} \subseteq C^1_{f(\delta)}$ is an unbounded subset of $f(\delta)$, then let $M^2_\delta$ be the unique $\tau$-model with universe $C^2_\delta$ such that $F \, {\rest} \,
    C^2_\delta$ is an isomorphism from $M^2_\delta$ onto $M^1_{F(\delta)}
    \, {\rest} \, \{F(\alpha) \colon \alpha \in C^2_\delta\}$. If $\delta \in S_1$
    does not satisfy the requirement above let $M^2_\delta$ be any $\tau$-model with universe $C^2_\delta$.
    
    Now it is easy to check that $\bar M^2 = \langle M^2_\delta \colon \delta \in
    S_2 \rangle$ is as required. 
    
    Second, assume $\ell =1$ (and we are in part (1)) and let  $\mathbf F^2$ be a $(\bar C^2,\Upsilon,\theta)$ colouring.  For each $(\bar C^1,\Upsilon,\theta)$-colouring $\mathbf F^1$ we ask:
    
    \underline{Question on $\mathbf F^1$}:  Do we have, for every club $E$ of
    $\lambda_2$ an increasing continuous function $h \colon \lambda_1 \rightarrow \lambda_2$ such that:  
    
    \begin{enumerate}  
        \item[(i)]  sup$(\text{Rang}(h)) \in E,$
    
        \item[(ii)]  the following set is a stationary subset of $\lambda_1$  $\{\delta \in S_1 \colon h(\delta) \in S_2$ and $C^2_{h(\delta)} \subseteq \{h(\delta) \colon \alpha \in C^1_\delta\}$ and $F(h(\delta)) = \delta$ and $(\forall \eta_1 \in {}^{(C^1_\delta)}\Upsilon)(\forall \eta_2 \in {}^{(C^2_{h(\delta)})} \Upsilon)[(\forall \alpha \in C^1_\delta)(\eta_1(\alpha) = \eta_2(h(\alpha))) \Rightarrow \mathbf F^1(\eta_1) = \mathbf F^2(\eta_2)\}$.
    \end{enumerate} 
        
    If for some $\mathbf F^1$ the answer is yes, we proceed as in the proof
    for $\ell =1$. Otherwise, for each such $\mathbf F^1$ there is a club
    $E[\mathbf F^1]$ exemplifying the negative answer.  Let $E =
    \bigcap\{E[\mathbf F^1] \colon \mathbf F^1$ is a $(\bar C^1,\Upsilon,\theta)$-colouring$\}$.  As there are $\le
    2^{2^{\lambda_1}}$ such functions $\mathbf F^1$ clearly $E$ is a 
    club of $\lambda_2$ and using (d) of the assumption, there is $h \colon \lambda_1 \rightarrow \lambda_2$ as there.  We can define $\mathbf F^1$ by $h, \mathbf F^2$ naturally and get a contradiction.  
    
    The case $\ell=0$ and part (3) are left to the reader.  
\end{PROOF}

Trying to apply Claim \ref{x.20.9}, clause (d) of the assumption looks the most
serious.  By the following, if each $C^1_\delta$ is a stationary subset
of $\delta$ (so $\kappa > \aleph_0$), it helps.

\begin{claim}\label{x.20.10}
    In Claim \ref{x.20.9}, if $2^{\lambda_1} <
    \lambda_2$, and we replace $(d)$ by $(d)^- + (f)$ below, \underline{then} we can conclude that  some $\bar C^3 = \langle C^3_\delta \colon \delta \in S
    \rangle$ has the $\theta$-{\rm MD}$_\ell$-property where $C^3_\delta \in
    D_\delta$ \underline{where}: 
    
    \begin{enumerate}
        \item[(f)]  $\bar D = \langle D_\delta \colon \delta \in S_2
        \rangle,D'_\delta$ is a filter on $C^2_\delta$ containing the
        co-bounded subsets of $C^2_\delta$ (or just $D_\delta \subseteq {\mathscr P}(C^2_\delta) \backslash \{A \subseteq C^2_\delta \colon \sup(A) < \delta\},$
    
        \item[(d)$^{-}$] $F \colon \lambda_2 \rightarrow \lambda_1$ 
        satisfies: for every club
        $E$ of $\lambda$ for some increasing continuous function $h$ from
        $\lambda_1$ to $\lambda_2$ with $\sup({\text{\rm Rang\/}}(h)) 
        \in {\text{\rm acc\/}}(E)$ the set
        $\{\delta \in S_1:h(\delta) \in S_2$ and $\{h(\alpha) \colon \alpha \in
        C^1_\alpha\} \cap C^2_\delta \in D'_\delta\}$ belong to
        $D_{\lambda_1}$.
    \end{enumerate}
\end{claim} 

\begin{PROOF}{\ref{x.20.10}}
    Let ${\mathscr  C}  \coloneqq \{\bar e \colon \bar e = \langle e_\delta \colon \delta \in  S_1 \rangle, e_\delta$ is a subset of $C^2_\delta\}$. So $|{\mathscr  C}| \le 2^{\lambda_1} < \lambda$ and for each $\bar e \in
    E$ we define $\bar C^{\bar e} \coloneqq \langle C^{\bar e}_\delta \colon \delta \in
    S^{\bar e}_2 \rangle$ by $S^{\bar e}_\delta = \{\delta \in S_2 \colon C^{\bar e}_\delta \in D_\delta\},$ where $C^{\bar e}_\delta = \{\alpha \in
    C^2_\delta \colon \mathrm{otp}(C^2_\delta \cap \alpha) \in e_{F(\delta)}\}$. If there exists some $\bar e \in {\mathscr  C},$ then $ \bar C^{\bar e}$ is as required, so we are done. If not, continue as in the proof of Claim \ref{x.20.9}
\end{PROOF}

\newpage

\section{Glossary}

\S0 \underline{Introduction}: 

\begin{itemize}
    \item Theorem \ref{0.1}: Many cases of MD-diamond.
 
    \item Theorem \ref{0.2}:  With division to $\lambda$.
    
    \item Conclusion \ref{0.3}: Connection to the principle from
    \cite{Sh:667}.
    
    \item Definition \ref{x.20.0109}: Having (partial) squares.
\end{itemize}

\S1 \underline{Super black box: sufficient conditions for $\cf(2^{\mu})$}:

\begin{itemize}
    \item Definition \ref{d.1}: Definition of $S$ being $\mathbf F-\theta$-small and the normal filter.

    \item Claim \ref{d.2}:  The filter is normal and $\Upsilon$ does not matter.

    \item Definition \ref{d.3}: Recalling pcf definitions.

    \item Remark \ref{d.4}:

    \item Definition \ref{d.5}:  The main notions are defined: 

    \begin{itemize}
        \item $\bar C$ is $\lambda$-\Bb-parameter, $(\lambda,\kappa,\chi)$-\Bb-parameter, 

        \item $\mathbf F$ is a $(\bar C,\mu,\theta)$-coloring, 

        \item $\bar C$ is a $(D,\mathbf F)$-\midia sequence,

        \item $\bar C$ is  good and $\bar C$ has the $(D,\mu,\theta)$-\Bb-property or the $(D,\mu,\theta)$-\Bb-property,  also (good$^+$-\Bb-property and \Bb$_\ell$.
    \end{itemize}

    \item Definition \ref{d.5Y}:

    \item Claim \ref{d.5A}:  Easy implications.

    \item Major Claim \ref{d.6}:  Sufficient conditions for $\bar C$ having the $(D,\theta)$-\Bb or properties ($D$ a filter on $\lambda,\lambda = \cf(2^\mu))$.

    \item Claim \ref{d.7}:  E.g., $\mu = \mu^\theta$ implies $\Sep(\mu,\theta)$, one of the assumptions of Claim \ref{d.6}.

    \item Remark \ref{d.7B}:  On variants.

    \item Definition \ref{d.8}:  $\Sep(\chi,\mu,\theta,\theta_1,\Upsilon)$, $\Sep(\mu,\theta)$.

    \item Claim \ref{d.8Z}:  Monotonicity properties of $\Sep$.

    \item Claim \ref{d.11}:  On some implications.

    \item Definition \ref{d.11A}: $\lambda$-MD-parameter.

    \item Claim \ref{d.11B}: Existence of good MD-parameters.

    \item Claim \ref{d.11C}:  Like Claim \ref{d.6} for MD-parameters (\ref{d.11} - \ref{d.11C} are not central).

    \item Claim \ref{d.12}:  Getting, e.g., $(D,\Upsilon,\theta)$-MD-property, i.e., $\langle M_\delta \colon \delta \in S \rangle,|M_\delta| = C_\delta,$ which guess (as in Theorem \ref{0.1}) for having $\bar C$ with the   $(D,\theta)$-property; also md-properties MD-properties are defined.

    \item Definition \ref{d.12A}: $\theta$-MD$_\ell$-property and the ideal ID$_{\lambda,\kappa,\chi,\Upsilon,\theta}(\bar C).$

    \item Claim \ref{d.12B}: Non saturation ideal above (needed in Claim \ref{d.12})

    \item Claim \ref{d.13}:  Quotations on the existence of good \Bb-parameters (by
    quoting).
    
    \item Definition \ref{d.14}:  A related ideal.
    
    \item Definition \ref{d.15}: $\mathbf F$ is $(S,\kappa,\chi)$-good
    (\ref{d.14}, \ref{d.15} not central).
    
    \item Conclusion \ref{d.16}:  Phrasing some conclusions.
    
    \item Question \ref{d.18}:  $\lambda$ weakly inaccessible, $\langle
    2^\theta \colon \theta < \lambda \rangle$ is eventually constant.
\end{itemize}

\S2 \underline{Upgrading to larger cardinals}: 

\begin{itemize}
     \item Claim \ref{x.20}:  For suitable $\lambda_2 > 2^{\lambda_1}$, we
    construct $(\lambda_2,\kappa,\chi)$-\Bb-parameters $\bar C$ which   has the
    $(D_\lambda,\Upsilon,\theta)$-property from $(\lambda_1,\kappa,\chi)$-\Bb-parameters which  has the $(D_\lambda,\Upsilon,\theta)$-parameters.
    
    \item Claim \ref{x.21A}:  Getting $\bar C^i = \langle C^i_\delta \colon \delta
    \in S_i \rangle$ sets $S_i$ pairwise disjoint for $i < \lambda$ in Claim \ref{x.20}.
    
    \item Conclusion \ref{x.21}:  If $\mu$ is strong limit, $\lambda = \text{
    cf}(\lambda) > 2^{2^\mu}$ is successor of singular,   then for every large enough $\kappa < \mu$, some good $(\lambda,\kappa,\kappa^+)$-\Bb-parameter has the  $(D_\lambda,2^\mu,\theta)$-\Bb-property for every $\theta < \mu$.
    
    \item Claim \ref{x.23}:  A generalization of Enge\l king Kar\l owicz theorem. [no pf].
    
    \item Claim \ref{x.2.9}:  Dealing with the successor of singulars. [$\exists$ pf].
    
    \item Claim \ref{x.24}:  We advance the picture in \ref{x.21} dealing with weakly inaccessible   (i.e. regular limit) not strong limit cardinals; at a minor price - $\mu$  is of uncountable cofinality.
    
    \item Claim \ref{x.20.9}: This is a stepping-up lemma used in Claim \ref{x.24} above.
    
    \item Claim \ref{x.20.10}:  Sufficient conditions for the ``lifting" 
    used in Claim \ref{x.20.9}.
\end{itemize}

\nocite{ignore-this-bibtex-warning} 

\newpage

\bibliographystyle{amsalpha}
\bibliography{shlhetal}

\end{document}